\documentclass[10pt]{article}

\usepackage{latexsym}
\usepackage{amsfonts}
\usepackage{amsmath}
\usepackage{amsthm}
\usepackage{mathrsfs}
\usepackage{dsfont}    
\usepackage{bbold}     
\usepackage[english]{babel}
\usepackage{caption}
\usepackage{epsfig}
\usepackage{float}
\usepackage{subfigure}
\usepackage{psfrag}
\usepackage{graphicx}
\usepackage{epsfig}

\DeclareMathOperator{\Val}{\matV}

\newtheorem{theorem}{Theorem}
\newtheorem*{prop*}{Theorem}

\newtheorem{defi}[theorem]{Notation}
\newtheorem{lemma}[theorem]{Lemma}

\newtheorem{rmk}[theorem]{Remark}

\newtheorem{const}[theorem]{Constraint}

\newcommand{\zerarcounters}{\setcounter{equation}{0}\setcounter{theorem}{0}}

\newcommand{\ZZZ}{\mathds{Z}}
\newcommand{\CCC}{\mathds{C}}
\newcommand{\NNN}{\mathds{N}}

\newcommand{\RRR}{\mathds{R}}
\newcommand{\TTT}{\mathds{T}}
\newcommand{\uno}{\mathds{1}}

\newcommand{\calA}{{\mathcal A}}

\newcommand{\calF}{{\mathcal F}}
\newcommand{\calG}{{\mathcal G}}

\newcommand{\calL}{{\mathcal L}}

\newcommand{\calP}{{\mathcal P}}

\newcommand{\gotE}{{\mathfrak E}}

\newcommand{\gotN}{{\mathfrak N}}

\newcommand{\gotR}{{\mathfrak R}}
\newcommand{\gotS}{{\mathfrak S}}
\newcommand{\gotT}{{\mathfrak T}}

\newcommand{\matA}{{\mathscr A}}
\newcommand{\matB}{{\mathscr B}}

\newcommand{\matG}{{\mathscr G}}

\newcommand{\matL}{{\mathscr L}}

\newcommand{\matO}{{\mathscr O}}

\newcommand{\matR}{{\mathscr R}}

\newcommand{\matV}{{\mathscr V}}

\newcommand{\Fullbox}{{\rule{2.0mm}{2.0mm}}}

\newcommand{\EP}{\hfill\Fullbox\vspace{0.2cm}}
\newcommand{\prova}{\noindent{\it Proof. }}
\newcommand{\io}{\infty}
\newcommand{\e}{\varepsilon}

\newcommand{\de}{\delta}

\newcommand{\n}{\nu}

\newcommand{\g}{\gamma}
\newcommand{\om}{\omega}
\newcommand{\h}{\eta}

\newcommand{\s}{\sigma}

\newcommand{\oo}{\boldsymbol{\omega}}
\newcommand{\mm}{\boldsymbol{\mu}}
\newcommand{\aaa}{\boldsymbol{\alpha}}

\newcommand{\hhh}{\boldsymbol{\eta}}
\newcommand{\nn}{\boldsymbol{\nu}}

\newcommand{\vzero}{\boldsymbol{0}}
\newcommand{\xx}{\boldsymbol{x}}
\newcommand{\yy}{\boldsymbol{y}}
\newcommand{\zz}{\boldsymbol{z}}
\newcommand{\ww}{\boldsymbol{w}}
\newcommand{\ee}{\boldsymbol{e}}

\newcommand{\bs}{\boldsymbol{s}}

\newcommand{\ff}{\boldsymbol{f}}
\newcommand{\ggg}{\boldsymbol{g}}
\newcommand{\cc}{\boldsymbol{c}}

\newcommand{\AAA}{\boldsymbol{A}}

\newcommand{\ii}{{\rm i}}
\newcommand{\To}{{\mathring{T}}}
\newcommand{\Tpru}{{\breve{T}}}
\newcommand{\thetao}{{\mathring{\theta}}}

\newcommand{\thetapru}{{\breve{\theta}}}
\newcommand{\Ga}{\Gamma}
\newcommand{\Gao}{\mathring{\Ga}}

\def\ins#1#2#3{\vbox to0pt{\kern-#2 \hbox{\kern#1 #3}\vss}\nointerlineskip}

\begin{document}

\title{\bf KAM theory in configuration space\\
and cancellations in the Lindstedt series}
%for perturbations of isochronous systems}

\author
{\bf Livia Corsi$^{1}$, Guido Gentile$^{1}$, Michela Procesi$^{2}$
\vspace{2mm}
\\ \small 
$^{1}$ Dipartimento di Matematica, Universit\`a di Roma Tre, Roma,
I-00146, Italy
\\ \small
$^{2}$ Dipartimento di Matematica, Universit\`a di
Napoli ``Federico II'', Napoli, I-80126, Italy
\\ \small 
E-mail: lcorsi@mat.uniroma3.it, gentile@mat.uniroma3.it, procesi@unina.it}

\date{}

\maketitle

\begin{abstract}
The KAM theorem for analytic quasi-integrable anisochronous Hamiltonian
systems yields that the perturbation expansion (Lindstedt series) for
any quasi-periodic solution with Diophantine frequency vector converges. 
If one studies the Lindstedt series by following a perturbation
theory approach, one finds that convergence is ultimately related
to the presence of cancellations between contributions of the
same perturbation order. In turn, this is due to symmetries in the
problem. Such symmetries are easily visualised in action-angle
coordinates, where KAM theorem is usually formulated, by exploiting
the analogy between Lindstedt series and perturbation expansions
in quantum field theory and, in particular, the possibility of
expressing the solutions in terms of tree graphs, which are
the analogue of Feynman diagrams. If the
unperturbed system is isochronous, Moser's modifying terms
theorem ensures that an analytic quasi-periodic solution with the
same Diophantine frequency vector as the unperturbed Hamiltonian
exists for the system obtained by adding a suitable constant
(counterterm) to the vector field. Also in this case, one can
follow the alternative approach of studying the perturbation expansion
for both the solution and the counterterm, and again
convergence of the two series is obtained as a consequence of deep
cancellations between contributions of the same order.
In this paper, we revisit Moser's theorem, by studying the
perturbation expansion one obtains by working in Cartesian coordinates.
We investigate the symmetries giving rise to the
cancellations which makes possible the convergence of the series. 
We find that the cancellation mechanism works in a
completely different way in Cartesian coordinates, and the
interpretation of the underlying symmetries in terms of tree graphs
is much more subtle than in the case of action-angle coordinates.
\end{abstract}

%\newpage

%\tableofcontents

%\newpage

%%%%%%%%%%%%%%%%%%%%%%%%%%%%%%%%%%%%%%%%%%%%%%%%%%%%%%%%%%%%%%%%%%%%%%%%%
%%%%%%%%%%%%%%%%%%%%%%%%%%%%%%%%%%%%%%%%%%%%%%%%%%%%%%%%%%%%%%%%%%%%%%%%%
\zerarcounters
\section{Introduction}
\label{sec:1}
%%%%%%%%%%%%%%%%%%%%%%%%%%%%%%%%%%%%%%%%%%%%%%%%%%%%%%%%%%%%%%%%%%%%%%%%%
%%%%%%%%%%%%%%%%%%%%%%%%%%%%%%%%%%%%%%%%%%%%%%%%%%%%%%%%%%%%%%%%%%%%%%%%%

Consider an isochronous Hamiltonian system, described by the Hamiltonian
$H(\aaa,\AAA)=\oo\cdot\AAA+\e f(\aaa,\AAA)$, with $f$ real analytic in
$\TTT^{d}\times\calA$ and $\calA$ an open subset of $\RRR^{d}$.
The corresponding Hamilton equation are
\begin{equation}
\dot \aaa = \oo + \e \partial_{\AAA}f(\aaa,\AAA) , \qquad
\dot \AAA = - \e \partial_{\aaa}f(\aaa,\AAA) .
\label{eq:1.1} \end{equation}
Let $(\aaa_{0}(t),\AAA_{0}(t))=(\aaa_{0}+\oo t,\AAA_{0})$ be a
solution of (\ref{eq:1.1}) for $\e=0$.
For $\e\neq0$, in general, there is no quasi-periodic solution
to (\ref{eq:1.1}) with frequency vector $\oo$ which reduces to
$(\aaa_{0}(t),\AAA_{0}(t))$ as $\e\to0$. However, one can prove that,
if $\e$ is small enough and $\oo$ satisfies some Diophantine condition,
then there is a `correction' $\mm(\e,\AAA_{0})$, analytic in both
$\e$ and $\AAA_{0}$, such that the modified equations
\begin{equation}
\dot \aaa = \oo + \e \partial_{\AAA}f(\aaa,\AAA) + \mm(\e,\AAA_{0}) ,
\qquad \dot \AAA = - \e \partial_{\aaa}f(\aaa,\AAA) ,
\label{eq:1.2} \end{equation}
admit a quasi-periodic solution with frequency vector $\oo$
which reduces to $(\aaa_{0}(t),\AAA_{0}(t))$ as $\e\to0$.  This is a
well known result, called the \emph{modifying terms theorem},
or \emph{translated torus theorem},
first proved by Moser \cite{M}. By writing the solution
as a power series in $\e$ (Lindstedt series), the existence
of an analytic solution means that the series converges.
This is ultimately related to some deep cancellations in the series;
see \cite{BG} for a review.

Equations like (\ref{eq:1.1}) naturally arise when studying
the stability of an elliptic equilibrium point. For instance,
one can think of a mechanical system near a minimum point for the
potential energy, where the Hamiltonian describing the system looks like
\begin{equation}
H(x_{1},\ldots,x_{n},y_{1},\ldots,y_{n}) =
\frac{1}{2} \sum_{j=1}^{d} 
\left( y_{j}^{2} + \om_{j}^{2} x_{j}^{2} \right) +
\e F(x_{1},\ldots,x_{n},\e) ,
\label{eq:1.3} \end{equation}
where $F$ is a real analytic function at least of third order
in its arguments, the vector $\oo=(\om_{1},\ldots,\om_{d})$ satisfies
some Diophantine condition, and the factor $\e$ can be
assumed to be obtained after a rescaling of the original coordinates.
Indeed, the corresponding Hamilton equations, written
in action-angle variables, are of the form (\ref{eq:1.1}).
Unfortunately, the action-angle variables are singular near the
equilibrium, and hence there are problems in the region where
one of the actions is much smaller than the others.
Thus, it can be worthwhile to work directly in the original
Cartesian coordinates. In fact, there has been a lot of interest
for KAM theory in configuration space, that is, without
action-angle variables; see for instance \cite{SZ,LM,DGJV}.

In the light of Moser's theorem of the modifying terms, 
one expects that, by writing
\begin{equation}
\ddot x_{j} + \om_{j}^{2} x_{j} =
- \e \partial_{x_{i}} F(x_{1},\ldots,x_{n},\e) + \h_{j} x_{j} ,
\label{eq:1.4} \end{equation}
and taking the (arbitrary) unperturbed solution $x_{0,j}(t)=
C_{j} \cos \om_{j} t + S_{j} \sin\om_{j} t= c_{j} {\rm e}^{i \om_{j} t} +
c_{j}^{*} {\rm e}^{-i \om_{j} t}$, $j=1,\ldots,d$,
there exists a function $\hhh(\e,\cc)$, analytic
both in $\e$ and $\cc=(c_{1},\ldots,c_{d})$,
such that by fixing $\h_{j}=\h_{j}(\e,\cc)$, there exists
a quasi-periodic solution to (\ref{eq:1.4}) with frequency vector
$\oo$, which reduces to the unperturbed one as $\e\to0$.

One can try to write again the solution as a power series in $\e$, and
study directly the convergence of the series. In general, when considering
the Lindstedt series of some KAM problem, first of all one identifies
the terms of the series which are an obstruction to convergence:
such terms are usually called resonances (or self-energy clusters,
by analogy to what happens in quantum field theory). Crudely speaking,
the series is given by the sum of infinitely many terms
(finitely many for each perturbation order), and each term
looks like a product of `small divisors' $\de(\oo\cdot\nn_{i})$,
$\nn_{i}\in\ZZZ^{d}$, times some harmless factors: a resonance
is a particular structure in the product which allows a dangerous
accumulation of small divisors. This phenomenon is very easily
visualised when each term of the series is graphically represented
as a tree graph (tree \emph{tout court} in the following),
that is, a set of points and lines connecting them in such a way
that no loop arises. In terms of trees, each line $\ell$ carries
a label $\nn_{\ell}\in\ZZZ^{d}$ (that one call momentum,
again inspired by the terminology of quantum field theory)
and with each such line a small divisor $\de(\oo\cdot\nn_{\ell})$ is
associated: then a resonance becomes a subgraph which is between
two lines with the same $\nn$ and hence the same small divisor
$\de(\oo\cdot\nn)$. A tree with a chain of resonances represents
a term of the series containing a factor $\de(\oo\cdot\nn)$ to
a very large power, and this produces a factorial $k!$ to some
positive power when bounding some terms contributing to the $k$-th
order in $\e$ of the Lindstedt series, so preventing a proof of convergence.

However, a careful analysis of the resonances shows that there are
cancellations to all perturbation orders. This is what can be proved
in the case of the standard anisochronous KAM theorem,
as first pointed out by Eliasson \cite{E}; see also \cite{G0,GBG},
for a proof which more deeply exploits the similarity with
the techniques of quantum field theory.

More precisely the cancellation mechanism works in the following way.
Given two lines $\ell_{1}$ and $\ell_{2}$ with the same small divisor,
consider all possible resonances which can be inserted between
$\ell_{1}$ and $\ell_{2}$: when summing together the numerical values
corresponding to all such resonances, there are compensations
and the sum is in fact much smaller than each summand
(for more details we refer to \cite{GBG,G3}).

For the isochronous case, already in action-angle variables \cite{BG},
there are some kinds of resonances which do not cancel each other.
Nevertheless there are other kinds of resonances
for which the gain factor due to the cancellation is more than
what needed (that is, one has a second order instead of a
first order cancellation). Thus, the hope naturally arises that one
can use the extra gain factors to compensate the lack of gain factors
for the first kind of resonances, and in fact this happens.
Indeed, the resonances for which there is no cancellation cannot
accumulate too much without entailing the presence of as many
resonances with the extra gain factors, in such a way that
the overall number of gain factors is, in average, one per resonance
(this is essentially the meaning of Lemma 5.4 in \cite{BG}).

When working in Cartesian coordinates, one immediately meets a
difficulty. If one writes down the lowest order resonances,
there is no cancellation at all. This is slightly surprising
because a cancellation is expected somewhere: if the resonances
do not cancel each other, in principle one can construct chains of
arbitrarily many resonances, which produce factorials in the formal
power series expansion. However, we shall show that there are
cancellations, as soon as one has at least two resonances.
So, one has the curious phenomenon that resonances which do not cancel
each other are allowed, but they cannot accumulate too much.
Moreover, the cancellation mechanism is more involved than in other cases
(including the same problem in action-angle variables). First of all,
the resonances are no longer diagonal in the momenta, that is,
the lines $\ell_{1}$ and $\ell_{2}$ considered above can have
different momenta $\nn_{1}$ and $\nn_{2}$.
Second, the cancellation does not operate simply by collecting
together all resonances to a given order and then summing the
corresponding numerical values. As we mentioned,
in this way no cancellation is produced: to obtain a cancellation
one has to consider all possible ways to connect to each other
two resonances. Thus, there is a cancellation only if there is
a chain of at least two resonances.

What emerges eventually is that working in Cartesian coordinates
rather complicates the analysis. On the other hand,
as remarked above, it can be worthwhile to investigate the problem
in Cartesian coordinates. Moreover, the cancellations are
due to remarkable symmetries in the problem, which can be
of interest by their own; in this regard we mention
the problem of the reducibility of the skew-product flows
with Bryuno base \cite{G1}, where the convergence of the corresponding
`Lindstedt series' is also due to some cancellation mechanism and
hence, ultimately, to some deep symmetry of the system.

In this paper we shall assume the \emph{standard Diophantine
condition} on the frequency vector $\oo$; see (\ref{eq:2.3}) below.
Of course one could consider more general Diophantine
conditions than the standard one (for instance a Bryuno condition
\cite{B}; see also \cite{G2} for a discussion using the
Lindstedt series expansion). This would make the analysis slightly
more complicate, without shedding further light on the problem.
An important feature of the Lindstedt series
method is that, from a conceptual point of view, the general
strategy is exactly the same independently of the kind
of coordinates one uses (and independently of the fact that
the system is a discrete map or a continuous flow; see \cite{BG0,GBG}).
What is really important for the analysis is the form
of the unperturbed solution: the simpler is such a solution,
the easier is the analysis. Of course, an essential issue
is that the system one wants to study is a perturbation of
one which is exactly soluble. This is certainly true in the case of
quasi-integrable Hamiltonian systems, but of course the range of
applicability is much wider, and includes also non-Hamiltonian
systems; see for instance \cite{GBD,G4}. Moreover an assumption
of this kind is more or less always implicit in whatsoever method
one can envisage to deal with small divisor problems of this kind.

In the anisochronous case, the cancellations are due to
symmetry properties of the model, as first pointed out by
Eliasson \cite{E}. Ultimately, the cancellation mechanism for the
resonances is deeply related to that assuring the formal solubility
of the equations of motions, which in turn is due to a symmetry property,
as already shown by Poincar\'e \cite{P}. Subsequently,
stressing further the analogy with quantum field theory,
Bricmont et al. showed that the cancellations can be interpreted
as a consequence of suitable Ward identities of the corresponding
field theory \cite{BGK} (see also \cite{DK}). In the isochronous case,
in terms of Cartesian coordinates the cancellation mechanism works
in a completely different way with respect to action-angle coordinates.
However, as we shall see, the cancellation is still related to
underlying symmetry properties: it would be interesting to provide an
interpretation of the symmetry properties that we find in terms of
some invariance property of the corresponding quantum field model.

%%%%%%%%%%%%%%%%%%%%%%%%%%%%%%%%%%%%%%%%%%%%%%%%%%%%%%%%%%%%%%%%%%%%%%%%%
%%%%%%%%%%%%%%%%%%%%%%%%%%%%%%%%%%%%%%%%%%%%%%%%%%%%%%%%%%%%%%%%%%%%%%%%%
\zerarcounters
\section{Statement of the results}
\label{sec:2}
%%%%%%%%%%%%%%%%%%%%%%%%%%%%%%%%%%%%%%%%%%%%%%%%%%%%%%%%%%%%%%%%%%%%%%%%%
%%%%%%%%%%%%%%%%%%%%%%%%%%%%%%%%%%%%%%%%%%%%%%%%%%%%%%%%%%%%%%%%%%%%%%%%%

Consider the ordinary differential equations
\begin{equation}
\ddot x_{j} + \om_{j}^{2} x_{j} + f_{j}(x_{1},\ldots,x_{d},\e) +
\h_{j} x_{j} = 0 , \qquad j=1,\ldots,d,
\label{eq:2.1} \end{equation}
where $\xx=(x_{1},\ldots,x_{d})\in\RRR^{d}$, $\e$ is real parameter
(\emph{perturbation parameter}), the function $\ff(\xx,\e)=
(f_{1}(\xx,\e),\ldots,f_{d}(\xx,\e))$ is real analytic in $\xx$ and $\e$
at $(\xx,\e)=(\vzero,0)$ and at least quadratic in $\xx$,
\begin{equation}
f_{j}(\xx,\e) = \sum_{p=1}^{\io} \e^{p} \!\!\!\!\!\!
\sum_{\substack{
s_{1},\ldots,s_{d} \ge 0 \\ s_{1}+\ldots+s_{d}=p+1}}
\!\!\!\!\!\! f_{j,s_{1},\ldots,s_{d}} \,
x_{1}^{s_{1}} \ldots x_{d}^{s_{d}} ,
\label{eq:2.2} \end{equation}
$\hhh=(\h_{1},\ldots,\h_{d})$ is a vector
of parameters, and the \emph{frequency vector}
(or \emph{rotation vector})
$\oo=(\om_{1},\ldots,\om_{d})$ satisfies the Diophantine condition
\begin{equation}
\left| \oo \cdot \nn \right| > \g_{0} \, |\nn|^{-\tau}
\qquad \forall \nn\in\ZZZ^{d}_{*} , 
\label{eq:2.3} \end{equation}
with $\ZZZ^{d}_{*}=\ZZZ^{d}\setminus\{\vzero\}$,
$\tau>d-1$ and $\g_{0} > 0$. Here and henceforth
$\cdot$ denotes the standard scalar product in $\RRR^{d}$,
and $|\nn|=|\n_{1}|+\ldots+|\n_{d}|$.

The equations (\ref{eq:2.1}) naturally arise from
equations of the form
\begin{equation}
\ddot x_{j} + \om_{j}^{2} x_{j} + g_{j}(x_{1},\ldots,x_{d}) +
\h_{j} x_{j} = 0 , \qquad j=1,\ldots,d,
\nonumber \end{equation}
with $\ggg=(g_{1},\ldots,g_{d})$ analytic and at least quadratic in $\xx$,
after rescaling the coordinates: $\xx \to \e \xx$.
Such rescaling makes sense if one wants to study the behaviour
of the system near the origin.

We look for quasi-periodic solutions $\xx(t)$ of (\ref{eq:2.1}) with
frequency vector $\oo$. Therefore we expand the function $\xx(t)$
by writing
\begin{equation}
\xx(t) = \sum_{\nn\in\ZZZ^{d}} {\rm e}^{\ii\nn\cdot\oo t} \, \xx_{\nn} ,
\label{eq:2.4} \end{equation}
and we denote by $\ff_{\nn}(\xx,\e)$ the $\nn$-th Fourier coefficient
of the function that we obtain by Taylor-expanding $\ff(\xx,\e)$ in
powers of $\xx$ and Fourier-expanding $\xx$ according to (\ref{eq:2.4}).

Thus, in Fourier space (\ref{eq:2.1}) becomes
\begin{equation}
\left[ (\oo\cdot\nn)^{2} - \om_{j}^{2} \right] x_{j,\nn} =
f_{j,\nn}(\xx,\e) + \h_{j} \, x_{j,\nn} .
\label{eq:2.5} \end{equation}

For $\e=0$, $\hhh=\vzero$, the vector $\xx^{(0)}(t)$ with components
\begin{equation}
x^{(0)}_{j}(t) = c_{j} {\rm e}^{\ii\om_{j}t} + c_{j}^{*} {\rm e}^{-\ii\om_{j}t} ,
\qquad j=1,\ldots,d ,
\label{eq:2.6} \end{equation}
is a solution of (\ref{eq:2.1}) for any choice of the
complex constant $\cc=(c_{1},\ldots,c_{d})$.
Here and henceforth $*$ denotes complex conjugation.

Define $\ee_{j}$ as the vector with components $\de_{ij}$
(Kronecker delta).
Then we can split (\ref{eq:2.5}) into two sets of equations,
called respectively the \emph{bifurcation equation}
and the \emph{range equation},
\begin{subequations}
\begin{align}
& f_{j,\s\ee_{j}}(\xx,\e) + \h_{j} \, x_{j,\s\ee_{j}} = 0 ,
\qquad j=1,\ldots,d , \quad \s=\pm 1, 
\label{eq:2.7a} \\
& \left[ (\oo\cdot\nn)^{2} - \om_{j}^{2} \right] x_{j,\nn} =
f_{j,\nn}(\xx,\e) + \h_{j} \, x_{j,\nn} ,
\qquad j=1,\ldots,d , \quad \nn \neq \pm \ee_{j} .
\label{eq:2.7b}
\end{align}
\label{eq:2.7} \end{subequations}
\vskip-.3truecm
\noindent We shall study both equations (\ref{eq:2.7}) simultaneously,
by showing that for all choices of the parameters $\cc$
there exist suitable \emph{counterterms} $\hhh$, depending
analytically on $\e$ and $\cc$, such that (\ref{eq:2.7})
admits a quasi-periodic solution with frequency vector $\oo$,
which is analytic in $\e$, $\cc$, and $t$. Moreover,
with the choice $x_{j,\ee_{j}}=c_{j}$ for all $j=1,\ldots,d$,
the counterterms are uniquely determined.

We formulate the following result.

%%%%%%%%%%%%%%%%%%%%%%%%%%%%%%%%%%%%%%%%%%%%%%%%%%%%%%%%%%%%%%%%%%%%%%%%%%
\begin{theorem} \label{thm:1}
Consider the system described by the equations (\ref{eq:2.1})
and let (\ref{eq:2.6}) be a solution at $\e=0$, $\hhh=\vzero$.
Set $\Gamma(\cc)=\max\{|c_{1}|,\ldots,|c_{d}|,1\}$. There exist
a positive constant $\eta_{0}$, small enough and independent
of $\e,\cc$, and a unique function $\hhh(\e,\cc)$,
holomorphic in the domain $|\e|\Gamma^{3}(\cc)\le \eta_{0}$ and
real for real $\e$, such that the system
\begin{equation}
\ddot x_{j} + \om_{j}^{2} x_{j} + f_{j}(x_{1},\ldots,x_{d},\e) +
\h_{j}(\e,\cc)\, x_{j} = 0 , \qquad j=1,\ldots,d,
\nonumber \end{equation}
admits a solution $\xx(t)=\xx(t,\e,\cc)$ of the form (\ref{eq:2.4}),
holomorphic in the domain
$|\e| \Gamma^{3}(\cc) {\rm e}^{3|\oo|\,|{\rm Im}\, t|}\le \eta_{0}$
and real for real $\e,t$,
with Fourier coefficients $x_{j,\ee_{j}}=c_{j}$
and $x_{j,\nn}=O(\e)$ if $\nn\neq\pm\ee_{j}$ for $j=1,\ldots,d$.
\end{theorem}
%%%%%%%%%%%%%%%%%%%%%%%%%%%%%%%%%%%%%%%%%%%%%%%%%%%%%%%%%%%%%%%%%%%%%%%%%%

The proof is organised as follows. After introducing the small
divisors and proving some simple preliminary properties
in Section \ref{sec:3}, we develop in Section \ref{sec:4}
a graphical representation for the power series of the counterterms
and the solution (tree expansion). In particular we perform
a multiscale analysis which allows us to single out the
contributions (self-energy clusters) which give problems when
trying to bound the coefficients of the series. In Section \ref{sec:5}
we show that, as far as such contributions are neglected, there
is no difficulty in obtaining power-like estimates on the
coefficients: these estimates, which are generalisations of the
Siegel-Bryuno bounds holding for anisochronous systems \cite{G0,GBG},
would imply the convergence of the series and hence analyticity.
In Section \ref{sec:6} we discuss how to deal with the
self-energy clusters: in particular we single out the leading
part of their contributions (localised values), which are proved
in Section \ref{sec:7} to satisfy some deep symmetry properties.
Finally, in Section \ref{sec:8} we show how the symmetry properties
can be exploited in order to obtain cancellations involving the
localised parts, in such a way that the remaining contributions
can still bounded in a summable way. This will yield the
convergence of the full series and hence the analyticity of
both the solution and the counterterms.

Note that the system dealt with in Theorem \ref{thm:1} can be
non-Hamiltonian. On the other hand the most general case for
a Hamiltonian system near a stable equilibrium allows for
Hamiltonians of the form
\begin{equation}
H(x_{1},\ldots,x_{n},y_{1},\ldots,y_{n}) =
\frac{1}{2} \sum_{j=1}^{d} 
\left( y_{j}^{2} + \om_{j}^{2} x_{j}^{2} \right) +
\e F(x_{1},\ldots,x_{n},y_{1},\ldots,y_{n},\e) ,
\label{eq:2.8} \end{equation}
which lead to the equations
\begin{equation}
\begin{cases}
\dot x_{j} = y_{j} + \e \partial_{y_{i}} F(\xx,\yy,\e) , \\
\dot y_{j} = - \om_{j}^{2} x_{j} -\e \partial_{x_{i}} F(\xx,\yy,\e) .
\end{cases}
\label{eq:2.9}
\end{equation}
Also in this case one can consider the modified equations
\begin{equation}
\begin{cases}
\dot x_{j} = y_{j} + \e \partial_{y_{i}} F(\xx,\yy,\e) , \\
\dot y_{j} = - \om_{j}^{2} x_{j} -\e \partial_{x_{i}} F(\xx,\yy,\e) +
\h_{j} x_{j} ,
\end{cases}
\label{eq:2.10}
\end{equation}
which are not of the form considered in Theorem \ref{thm:1}.
However, a result in the same spirit as Theorem \ref{thm:1} still holds.

%%%%%%%%%%%%%%%%%%%%%%%%%%%%%%%%%%%%%%%%%%%%%%%%%%%%%%%%%%%%%%%%%%%%%%%%%%
\begin{theorem} \label{thm:2}
Consider the system described by the equations (\ref{eq:2.10})
and let $(\xx^{(0)}(t),\yy^{(0)}(t))$ be a solution at $\e=0$,
$\hhh=\vzero$, with $\xx^{(0)}(t)$ given by (\ref{eq:2.6})
and $\yy^{(0)}(t)=\dot\xx^{(0)}(t)$.
Set $\Gamma(\cc)=\max\{|c_{1}|,\ldots,|c_{d}|,1\}$.
Then there exist a positive constant $\eta_{0}$,
small enough and independent of $\e,\cc$,
and a unique function $\hhh(\e,\cc)$, holomorphic in the domain
$|\e|\Gamma^{3}(\cc)\le \eta_{0}$ and
real for real $\e$, such that the system
\begin{equation}
\begin{cases}
\dot x_{j} = y_{j} + \e \partial_{y_{i}} F(\xx,\yy,\e) , \\
\dot y_{j} = - \om_{j}^{2} x_{j} -\e \partial_{x_{i}} F(\xx,\yy,\e) +
\h_{j}(\e,\cc)\, x_{j} 
\end{cases}
\nonumber
\end{equation}
admits a solution $(\xx(t,\e,\cc),\yy(t,\e,\cc))$,
holomorphic in the domain $|\e| \Gamma^{3}(\cc) {\rm e}^{3|\oo|
\,|{\rm Im}\, t|}\le \eta_{0}$ and real for real $\e,t$,
with Fourier coefficients $x_{j,\ee_{j}}=y_{j,\ee_{j}}/\ii
\om_{j}=c_{j}$ and $x_{j,\nn}=y_{j,\nn}=O(\e)$ if $\nn\neq\pm\ee_{j}$
for $j=1,\ldots,d$.
\end{theorem}
%%%%%%%%%%%%%%%%%%%%%%%%%%%%%%%%%%%%%%%%%%%%%%%%%%%%%%%%%%%%%%%%%%%%%%%%%%

The proof follows the same lines as that of Theorem \ref{thm:1}, and it
is discussed in Appendices \ref{app:a} and \ref{app:b}.
Finally in Appendix \ref{app:c} we briefly sketch an alternative
approach based on the resummation of the perturbation series.

%%%%%%%%%%%%%%%%%%%%%%%%%%%%%%%%%%%%%%%%%%%%%%%%%%%%%%%%%%%%%%%%%%%%%%%%%%
%%%%%%%%%%%%%%%%%%%%%%%%%%%%%%%%%%%%%%%%%%%%%%%%%%%%%%%%%%%%%%%%%%%%%%%%%%
\zerarcounters
\section{Preliminary results}
\label{sec:3}
%%%%%%%%%%%%%%%%%%%%%%%%%%%%%%%%%%%%%%%%%%%%%%%%%%%%%%%%%%%%%%%%%%%%%%%%%%
%%%%%%%%%%%%%%%%%%%%%%%%%%%%%%%%%%%%%%%%%%%%%%%%%%%%%%%%%%%%%%%%%%%%%%%%%%

We shall denote by $\NNN$ the set of (strictly) positive integers,
and set $\ZZZ_{+}=\NNN\cup\{0\}$.
For any $j=1,\ldots,d$ and $\nn\in\ZZZ^{d}$ define
the \emph{small divisors}
\begin{equation}
\de_{j}(\oo\cdot\nn) := \min\{
\left| \oo\cdot\nn - \om_{j} \right|,
 \left| \oo\cdot\nn + \om_{j} \right|\}=
|\oo\cdot(\nn-\s(\nn,j)\,\ee_{j})|,
\label{eq:3.1} \end{equation}
where $\s(\nn,j)$ is the minimizer.
Note that the Diophantine condition (\ref{eq:2.3}) implies that
\begin{subequations}
\begin{align}
& \de_{j}(\oo\cdot\nn) \ge \g |\nn|^{-\tau}
\qquad \forall j =1\ldots,d, \quad
\forall \nn\neq \vzero,\s(\nn,j)\,\ee_{j},
\label{eq:3.2a} \\
& \de_{j}(\oo\cdot\nn) \! + \! \de_{j'}(\oo\cdot\nn')
\ge \g |\nn-\nn'|^{-\tau}
\quad \forall j,j' \!\! = \!\! 1,\ldots,d, \quad
\forall \nn \! \neq \! \nn',\,
\nn \!\! - \!\!  \nn' \! \neq \! \s(\nn,j)\,\ee_{j} \! - \! 
\s(\nn',j')\,\ee_{j'} ,
\label{eq:3.2b}
\end{align}
\label{eq:3.2} \end{subequations}
\vskip-.3truecm
\noindent for a suitable positive $\g>0$. We can (and shall)
assume that $\g$ is sufficiently smaller than $\g_{0}$,
and hence than $\de(0)=\min\{|\om_{1}|,\ldots,|\om_{d}|\}$
and $\om := \min\{||\om_{i}|-|\om_{j}|| : 1 \le i<j \le d\}$.

%%%%%%%%%%%%%%%%%%%%%%%%%%%%%%%%%%%%%%%%%%%%%%%%%%%%%%%%%%%%%%%%%%%%%%%%%%
\begin{lemma} \label{lem:3.1}
Given $\nn,\nn'\in \ZZZ^{d}$, with $\nn\neq\nn'$, and
$\de_{j}(\oo\cdot\nn)=\de_{j'}(\oo\cdot\nn')$ for some
$j,j'\in\{1,\ldots,d\}$, then either
$|\nn-\nn'| \ge |\nn|+|\nn'|-2$ or $|\nn-\nn'|=2$.
\end{lemma}
%%%%%%%%%%%%%%%%%%%%%%%%%%%%%%%%%%%%%%%%%%%%%%%%%%%%%%%%%%%%%%%%%%%%%%%%%%

%%%%%%%%%%%%%%%%%%%%%%%%%%%%%%%%%%%%%%%%%%%%%%%%%%%%%%%%%%%%%%%%%%%%%%%%%%
\prova One has $\de_{j}(\oo\cdot\nn)=|\oo\cdot\nn-\s\om_{j}|$ and
$\de_{j'}(\oo\cdot\nn')=|\oo\cdot\nn'-\s'\om_{j'}|$,
with $\s=\s(\nn,j)$ and $\s'=\s(\nn',j')$.
Set $\bar\nn=\nn-\s\ee_{j}$ and $\bar\nn'=\nn'-\s'\ee_{j'}$.
By the Diophantine condition (\ref{eq:2.3}) one can have
$\de_{j}(\oo\cdot\nn)=\de_{j'}(\oo\cdot\nn')$, and hence
$|\oo\cdot\bar\nn|=|\oo\cdot\bar\nn'|$, if and only if
$\bar\nn=\pm\bar\nn'$.

If $\bar\nn=-\bar\nn'$ then for $\s=-\s'$ one has
$|\nn-\nn'|=|\nn|+|\nn'|$, while for $\s=\s'$ one obtains
$|\nn-\nn'| \ge |\nn|+|\nn'|-2$.
If $\bar\nn=\bar\nn'$ and $j=j'$ one has $\n_{i}=\n_{i}'$ for all
$i\neq j$ and $\n_{j}-\s=\n_{j}'-\s'$, and hence $|\n_{j}-\n_{j}'|=2$.
If $\bar\nn=\bar\nn'$ and $j\neq j'$ then $\n_{i}=\n_{i}'$ for all
$i\neq j,j'$, while $\n_{j}-\s=\n_{j}'$ and $\n_{j'}=\n_{j'}'-\s'$,
and hence $|\n_{j}-\n_{j}'|=|\n_{j'}-\n_{j'}'|=1$.\EP
%%%%%%%%%%%%%%%%%%%%%%%%%%%%%%%%%%%%%%%%%%%%%%%%%%%%%%%%%%%%%%%%%%%%%%%%%%

%%%%%%%%%%%%%%%%%%%%%%%%%%%%%%%%%%%%%%%%%%%%%%%%%%%%%%%%%%%%%%%%%%%%%%%%%%
\begin{lemma} \label{lem:3.2} 
Let $\nn,\nn'\in \ZZZ^{d}$ be such that $\nn\neq\nn'$ and,
for some $n\in\ZZZ_{+}$, $j,j'\in\{1,\ldots,d\}$, both $\de_{j}
(\oo\cdot\nn) \le 2^{-n}\g$ and $\de_{j'}(\oo\cdot\nn') \le 2^{-n}\g$
hold.
Then either $|\nn-\nn'|>2^{(n-2)/\tau}$
or $|\nn-\nn'| = 2$ and $\de_{j}(\oo\cdot\nn)=\de_{j'}(\oo\cdot\nn')$.
\end{lemma}
%%%%%%%%%%%%%%%%%%%%%%%%%%%%%%%%%%%%%%%%%%%%%%%%%%%%%%%%%%%%%%%%%%%%%%%%%%

%%%%%%%%%%%%%%%%%%%%%%%%%%%%%%%%%%%%%%%%%%%%%%%%%%%%%%%%%%%%%%%%%%%%%%%%%%
\prova Write $\de_{j}(\oo\cdot\nn)=|\oo\cdot\nn-\s\om_{j}|$ and
$\de_{j'}(\oo\cdot\nn')=|\oo\cdot\nn'-\s'\om_{j'}|$, with
$\s=\s(\nn,j)$ and $\s'=\s(\nn',j')$, and set $\bar\nn=\nn-\s\ee_{j}$
and $\bar\nn'=\nn'-\s'\ee_{j'}$ as above.

If $\bar\nn\neq\bar\nn'$, by the Diophantine condition (\ref{eq:3.2b}),
one has
\begin{equation}
\g \left|\bar\nn-\bar\nn' \right|^{-\tau} <
\left| \oo\cdot (\bar\nn-\bar\nn') \right| \le
\left| \oo\cdot\bar\nn \right| +
\left| \oo\cdot\bar\nn'\right| < 2^{-(n-1)}\g ,
\nonumber \end{equation}
which implies $|\bar\nn-\bar\nn'|>2^{(n-1)/\tau}$, and hence
we have $|\nn-\nn'|>2^{(n-2)/\tau}$ in such a case.

If $\bar\nn=\bar\nn'$ then, as in Lemma \ref{lem:3.1}, one has $|\nn-\nn'|
=2$ and $\de_{j}(\oo\cdot\nn)=\de_{j'}(\oo\cdot\nn')$.\EP
%%%%%%%%%%%%%%%%%%%%%%%%%%%%%%%%%%%%%%%%%%%%%%%%%%%%%%%%%%%%%%%%%%%%%%%%%%

%%%%%%%%%%%%%%%%%%%%%%%%%%%%%%%%%%%%%%%%%%%%%%%%%%%%%%%%%%%%%%%%%%%%%%%%%%
\begin{rmk}\label{rmk:3.3}
Note that $|\nn-\nn'| \le 2$ and $\de_{j}(\oo\cdot\nn)=\de_{j'}
(\oo\cdot\nn')$ if and only if $\nn-\nn'=\s(\nn,j)\ee_{j}-\s(\nn',j')
\ee_{j'}$.
\end{rmk}
%%%%%%%%%%%%%%%%%%%%%%%%%%%%%%%%%%%%%%%%%%%%%%%%%%%%%%%%%%%%%%%%%%%%%%%%%%

%%%%%%%%%%%%%%%%%%%%%%%%%%%%%%%%%%%%%%%%%%%%%%%%%%%%%%%%%%%%%%%%%%%%%%%%%%
\begin{lemma} \label{lem:3.4}
Let $\nn_{1},\ldots,\nn_{p}\in \ZZZ^{d}$ and $j_{1},\ldots,j_{p}
\in \{1,\ldots,d\}$, with $p \ge 2$, be such that
$|\nn_{i}-\nn_{i-1}| \le 2$ and $\de_{j_{i}}(\oo\cdot\nn_{i})=
\de_{j_{1}}(\oo\cdot\nn_{1}) \le \g$ for $i=2,\ldots,p$.
Then $|\nn_{1}-\nn_{p}| \le 2$.
\end{lemma}
%%%%%%%%%%%%%%%%%%%%%%%%%%%%%%%%%%%%%%%%%%%%%%%%%%%%%%%%%%%%%%%%%%%%%%%%%%

%%%%%%%%%%%%%%%%%%%%%%%%%%%%%%%%%%%%%%%%%%%%%%%%%%%%%%%%%%%%%%%%%%%%%%%%%%
\prova Set $\s_{i}=\s(\nn_{i},j_{i})$ and
$\bar\nn_{i}=\nn_{i}-\s_{i}\ee_{j_{i}}$ for $i=1,\ldots,p$.
For all $i=2,\ldots,p$, the assumption $\de_{j_{i}}(\oo\cdot\nn_{i})=
\de_{j_{i-1}}(\oo\cdot\nn_{i-1})$ implies $\bar\nn_{i}=\pm\bar\nn_{i-1}$, 
which in turn yields $\bar\nn_{i}=\bar\nn_{i-1}$ since
$|\nn_{i}-\nn_{i-1}| \le 2$. In particular $\bar\nn_{1}=\bar\nn_{p}$,
and hence $|\nn_{1}-\nn_{p}| \le 2$.\EP
%%%%%%%%%%%%%%%%%%%%%%%%%%%%%%%%%%%%%%%%%%%%%%%%%%%%%%%%%%%%%%%%%%%%%%%%%%

%%%%%%%%%%%%%%%%%%%%%%%%%%%%%%%%%%%%%%%%%%%%%%%%%%%%%%%%%%%%%%%%%%%%%%%%%%
%%%%%%%%%%%%%%%%%%%%%%%%%%%%%%%%%%%%%%%%%%%%%%%%%%%%%%%%%%%%%%%%%%%%%%%%%%
\zerarcounters
\section{Multiscale analysis and diagrammatic rules}
\label{sec:4}
%%%%%%%%%%%%%%%%%%%%%%%%%%%%%%%%%%%%%%%%%%%%%%%%%%%%%%%%%%%%%%%%%%%%%%%%%%
%%%%%%%%%%%%%%%%%%%%%%%%%%%%%%%%%%%%%%%%%%%%%%%%%%%%%%%%%%%%%%%%%%%%%%%%%%

As we are looking for $\xx(t,\e,\cc)$ and $\hhh(\e,\cc)$ analytic in $\e$,
we formally write
\begin{equation}
x_{j,\nn}=\sum_{k=0}^{\io}\e^{k}x_{j,\nn}^{(k)},\qquad
\h_{j}=\sum_{k=1}^{\io}\e^{k}\h_{j}^{(k)}.
\label{eq:4.1} \end{equation}
It is not difficult to see that using (\ref{eq:4.1}) in (\ref{eq:2.7})
one can recursively compute (at least formally) the coefficients
$x_{j,\nn}^{(k)}$, $\h_{j}^{(k)}$ to all orders. Here we introduce a
graphical representation for each contribution to $x_{j,\nn}^{(k)}$,
$\h_{j}^{(k)}$, which will allow us to study the convergence of the series.

%%%%%%%%%%%%%%%%%%%%%%%%%%%%%%%%%%%%%%%%%%%%%%%%%%%%%%%%%%%%%%%%%%%%%%%%%%
\subsection{Trees} \label{sec:4.1}
%%%%%%%%%%%%%%%%%%%%%%%%%%%%%%%%%%%%%%%%%%%%%%%%%%%%%%%%%%%%%%%%%%%%%%%%%%

\begin{figure}[!ht]
\begin{center}{
\psfrag{th}{$\theta\!=$}
\includegraphics[width=12cm]{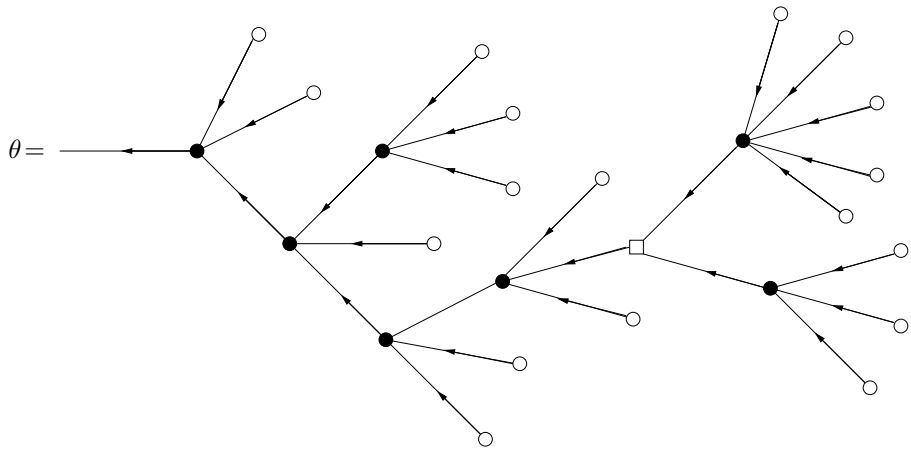}
}
\caption{\label{fig:1}
\footnotesize{An unlabelled tree: the arrows on the lines
all point toward the root, according to the tree partial ordering.
}}
\end{center}
\end{figure}

A graph is a connected set of points and lines.
A \emph{tree} $\theta$ is a graph with no cycle,
such that all the lines are oriented toward a unique point
(\emph{root}) which has only one incident line (\emph{root line}).
All the points in a tree except the root are called \emph{nodes}.
The orientation of the lines in a tree induces a partial ordering 
relation ($\preceq$) between the nodes and the lines: we can imagine
that each line carries an arrow pointing toward the root;
see Figure \ref{fig:1}. Given two nodes $v$ and $w$,
we shall write $w \prec v$ every time $v$ is along the path
(of lines) which connects $w$ to the root.

%%%%%%%%%%%%%%%%%%%%%%%%%%%%%%%%%%%%%%%%%%%%%%%%%%%%%%%%%%%%%%%%%%%%%%%%%%
% figure 1
%%%%%%%%%%%%%%%%%%%%%%%%%%%%%%%%%%%%%%%%%%%%%%%%%%%%%%%%%%%%%%%%%%%%%%%%%%
%%%%%%%%%%%%%%%%%%%%%%%%%%%%%%%%%%%%%%%%%%%%%%%%%%%%%%%%%%%%%%%%%%%%%%%%%%

We call $E(\theta)$ the set of \emph{end nodes} in $\theta$,
that is, the nodes which have no entering line, and $V(\theta)$
the set of \emph{internal nodes} in $\theta$, that is, the set of
nodes which have at least one entering line.
Set $N(\theta)=E(\theta) \amalg V(\theta)$. For all $v\in N(\theta)$
denote by $s_{v}$ the number of lines entering the node $v$.

%%%%%%%%%%%%%%%%%%%%%%%%%%%%%%%%%%%%%%%%%%%%%%%%%%%%%%%%%%%%%%%%%%%%%%%%%%
\begin{rmk} \label{rmk:4.1}
One has ${\displaystyle \sum_{v\in V(\theta)} s_{v}= |N(\theta)|-1}$.
\end{rmk}
%%%%%%%%%%%%%%%%%%%%%%%%%%%%%%%%%%%%%%%%%%%%%%%%%%%%%%%%%%%%%%%%%%%%%%%%%%

We denote by $L(\theta)$ the set of lines in $\theta$.
We call \emph{internal line} a line exiting
an internal node and \emph{end line} a line exiting an end node.
Since a line $\ell\in L(\theta)$ is uniquely identified
with the node $v$ which it leaves, we may write $\ell = \ell_{v}$.
We write $\ell_{w} \prec \ell_{v}$ if $w\prec v$; we say that a node
$w$ precedes a line $\ell$, and write $w\prec \ell$, if $\ell_{w}
\preceq \ell$.

%%%%%%%%%%%%%%%%%%%%%%%%%%%%%%%%%%%%%%%%%%%%%%%%%%%%%%%%%%%%%%%%%%%%%%%%%%
\begin{defi}\label{def:4.2}
$\null$\\ 
(1) If $\ell$ and $\ell'$ are two comparable lines, i.e.,
$\ell' \prec \ell$, we denote by $\calP(\ell,\ell')$ the
(unique) path of lines connecting $\ell'$ to $\ell$.\\
(2) Each internal line $\ell\in L(\theta)$ can be seen as the root line
of the tree $\theta_{\ell}$ whose nodes and lines are those of $\theta$
which precede $\ell$, that is,
$N(\theta_{\ell})=\{v'\in N(\theta) : v' \prec \ell\}$ and
$L(\theta_{\ell})=\{\ell'\in L(\theta) : \ell' \preceq \ell\}$.
\end{defi}
%%%%%%%%%%%%%%%%%%%%%%%%%%%%%%%%%%%%%%%%%%%%%%%%%%%%%%%%%%%%%%%%%%%%%%%%%%

%%%%%%%%%%%%%%%%%%%%%%%%%%%%%%%%%%%%%%%%%%%%%%%%%%%%%%%%%%%%%%%%%%%%%%%%%%
\subsection{Tree labels} \label{sec:4.2}
%%%%%%%%%%%%%%%%%%%%%%%%%%%%%%%%%%%%%%%%%%%%%%%%%%%%%%%%%%%%%%%%%%%%%%%%%%

\begin{figure}[!ht]
\begin{center}{
\psfrag{a}{(a)}
\psfrag{b}{(b)}
\psfrag{c}{(c)}
\psfrag{v}{$v$}
\psfrag{njs}{$\nn_{v}\,j_{v}\,\s_{v}$}
\psfrag{jk}{$j_{v}\,k_{v}$}
\psfrag{jks}{$j_{v}\,k_{v}\,\s_{v}$}
\includegraphics[width=16cm]{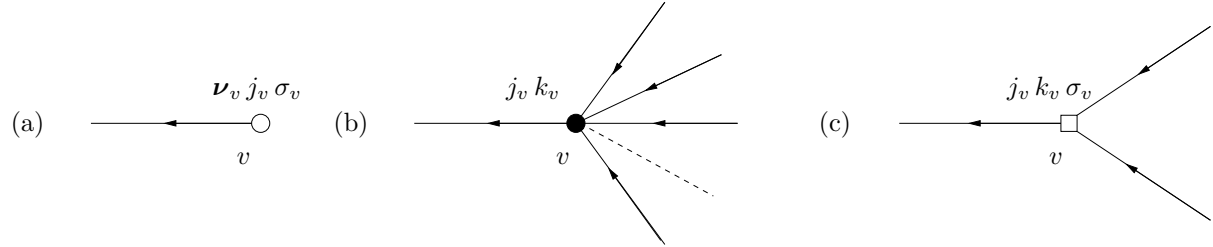}
}
\caption{\label{fig:2}
\footnotesize{Nodes and labels associated with the nodes:
(a) end node $v$ with $s_{v}=0$, $j_{v}\in\{1,\ldots,d\}$,
$\s_{v}\in\{\pm\}$, and $\nn_{v}=\s_{v}\ee_{j_{v}}$
(cf. Section \ref{sec:4.3});
(b) internal node $v$ with $s_{v} \ge 2$, $j_{v}\in\{1,\ldots,d\}$,
and $k_{v}=s_{v}-1$ (cf. Section \ref{sec:4.3});
(c) internal node $v$ with $s_{v} = 2$, $j_{v}\in\{1,\ldots,d\}$
$k_{v}=0$, $\s_{v}\in\{\pm\}$ (cf. Section \ref{sec:4.3}).}}
\end{center}
\end{figure}

With each end node $v\in E(\theta)$ we associate
a \emph{mode} label $\nn_{v}\in \ZZZ^{d}$, a \emph{component} label
$j_{v}\in\{1,\ldots,d\}$, and a \emph{sign} label $\s_{v}\in\{\pm\}$;
see Figure \ref{fig:2}.
We call $E_{j}^{\s}(\theta)$ the set of end nodes $v\in E(\theta)$
such that $j_{v}=j$ and $\s_{v}=\s$.

With each internal node $v\in V(\theta)$ we associate a
\emph{component} label $j_{v}\in\{1,\ldots,d\}$,
and an \emph{order} label $k_{v}\in\ZZZ_{+}$.
Set $V_{0}(\theta)=\{v\in V(\theta):k_{v}=0\}$
and $N_{0}(\theta)=E(\theta) \amalg V_{0}(\theta)$.
We also associate a sign label $\s_{v}\in\{\pm\}$ with each
$v\in V_{0}(\theta)$. The internal nodes $v$ with $k_{v}\ge 1$ will be
drawn as black bullets, while the end nodes and the internal nodes
with $k_{v}=0$ will be drawn as white bullets and white squares,
respectively; see Figure \ref{fig:2}.

%%%%%%%%%%%%%%%%%%%%%%%%%%%%%%%%%%%%%%%%%%%%%%%%%%%%%%%%%%%%%%%%%%%%%%%%%%
% figure 2
%%%%%%%%%%%%%%%%%%%%%%%%%%%%%%%%%%%%%%%%%%%%%%%%%%%%%%%%%%%%%%%%%%%%%%%%%%
%%%%%%%%%%%%%%%%%%%%%%%%%%%%%%%%%%%%%%%%%%%%%%%%%%%%%%%%%%%%%%%%%%%%%%%%%%

\begin{figure}[!ht]
\begin{center}{
\psfrag{njsn}{$\nn_{\ell}\,j_{\ell}\,\s_{\ell}\,n_{\ell}$}
\psfrag{l}{$\ell$}
\includegraphics[width=4cm]{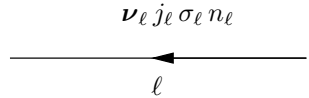}
}
\caption{\label{fig:3}
\footnotesize{Labels associated with a line.
One has $\s_{\ell}=\s(\nn_{\ell},j_{\ell})$ (cf. Section \ref{sec:4.3})
Moreover if $\ell=\ell_{v}$ then $j_{\ell}=j_{v}$;
if $v\in V_{0}(\theta)$ one has also $\s_{\ell}=\s_{v}$;
if $\nn_{\ell}=\s_{\ell}\ee_{j_{\ell}}$ then
$n_{\ell}=-1$, otherwise $n_{\ell}\ge 0$
(cf. Section \ref{sec:4.3}).}}
\end{center}
\end{figure}

With each line $\ell$ we associate a \emph{momentum} label
$\nn_{\ell} \in \ZZZ^{d}$, a \emph{component} label
$j_{\ell}\in\{1,\ldots,d\}$, a \emph{sign} label $\s_{\ell}\in\{\pm\}$,
and \emph{scale} label $n_{\ell}\in\ZZZ_{+} \cup \{-1\}$;
see Figure \ref{fig:3}.

%%%%%%%%%%%%%%%%%%%%%%%%%%%%%%%%%%%%%%%%%%%%%%%%%%%%%%%%%%%%%%%%%%%%%%%%%%
% figure 3
%%%%%%%%%%%%%%%%%%%%%%%%%%%%%%%%%%%%%%%%%%%%%%%%%%%%%%%%%%%%%%%%%%%%%%%%%%
%%%%%%%%%%%%%%%%%%%%%%%%%%%%%%%%%%%%%%%%%%%%%%%%%%%%%%%%%%%%%%%%%%%%%%%%%%

Denote by $s_{v,j}$ the number of lines
$\ell$ with component label $j_{\ell}=j$ entering the node $v$,
and with $r_{v,j,\s}$ the number of end lines with component label $j$
and sign label $\s$ which enter the node $v$.
Of course $s_{v} = s_{v,1}+\ldots+s_{v,d}$ and
$s_{v,j} \ge r_{v,j,+}+r_{v,j,-}$ for all $j=1,\ldots,d$.

Finally call
\begin{equation}
k(\theta) := \sum_{v\in V(\theta)} k_{v} 
\nonumber \end{equation}
the \emph{order} of the tree $\theta$.

In the following we shall call \emph{trees} tout court
the trees with labels, and we shall use the
term \emph{unlabelled trees} for the trees without labels.

%%%%%%%%%%%%%%%%%%%%%%%%%%%%%%%%%%%%%%%%%%%%%%%%%%%%%%%%%%%%%%%%%%%%%%%%%%
\subsection{Constraints on the tree labels} \label{sec:4.3}
%%%%%%%%%%%%%%%%%%%%%%%%%%%%%%%%%%%%%%%%%%%%%%%%%%%%%%%%%%%%%%%%%%%%%%%%%%

%%%%%%%%%%%%%%%%%%%%%%%%%%%%%%%%%%%%%%%%%%%%%%%%%%%%%%%%%%%%%%%%%%%%%%%%%%
\begin{const}\label{const:4.3}
We have the following constraints on the labels of the nodes
(see Figure \ref{fig:2}):\\
(1) if $v\in V(\theta)$ one has $s_{v} \ge 2$;\\
(2) if $v\in E(\theta)$ one has $\nn_{v}=\s_{v}\ee_{j_{v}}$;\\
(3) if $v \in V(\theta)$ then $k_{v}=s_{v}-1$,
except for $s_{v}=2$, where both $k_{v}=1$ and $k_{v}=0$ are allowed.
\end{const}
%%%%%%%%%%%%%%%%%%%%%%%%%%%%%%%%%%%%%%%%%%%%%%%%%%%%%%%%%%%%%%%%%%%%%%%%%%

\begin{figure}[!ht]
\begin{center}{
\psfrag{a}{(a)}
\psfrag{b}{(b)}
\psfrag{c}{(c)}
\psfrag{v}{$v$}
\psfrag{lv}{$\ell_{v}$}
\psfrag{l'}{$\ell_{1}$}
\psfrag{l}{$\ell_{2}$}
\psfrag{j0s}{$j_{v} \, 0 \, \s_{v}$}
\psfrag{nsjv}{$\nn_{\ell_{v}} \, \s_{v} \, j_{v}$}
\psfrag{nsj'}{$\s_{v}\ee_{j_{v}} \, \s_{v} \, j_{v}$}
\psfrag{nsj}{$\nn_{\ell_{2}} \, \s_{v} \, j_{v}$}
\includegraphics[width=6cm]{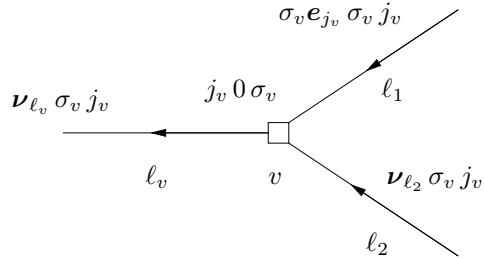}
}
\caption{\label{fig:4}
\footnotesize{If there is an internal node $v$ with $k_{v}=0$
then $s_{v}=2$ and the following constraints are imposed on the
other labels:
$\s_{\ell_{v}}=\s_{\ell_{1}}=\s_{\ell_{2}}=\s_{v}$;
$j_{\ell_{v}}=j_{\ell_{1}}=j_{\ell_{2}}=j_{v}$;
either $\nn_{\ell_{1}}=\s_{v}\ee_{j_{v}}$ and
$\nn_{\ell_{2}} \neq \s_{v}\ee_{j_{v}}$ (as in the figure) or
$\nn_{\ell_{2}}=\s_{v}\ee_{j_{v}}$ and
$\nn_{\ell_{1}} \neq \s_{v}\ee_{j_{v}}$. (The scale labels are not shown).}}
\end{center}
\end{figure}

%%%%%%%%%%%%%%%%%%%%%%%%%%%%%%%%%%%%%%%%%%%%%%%%%%%%%%%%%%%%%%%%%%%%%%%%%%
\begin{const}\label{const:4.4}
The following constraints will be imposed on the labels of the lines:\\
(1) $j_{\ell}=j_{v}$, $\nn_{\ell}=\nn_{v}$, and
$\s_{\ell}=\s_{v}$ if $\ell$ exits $v\in E(\theta)$;\\
(2) $j_{\ell}=j_{v}$ if $\ell$ exits $v\in V(\theta)$;\\
(3) if $\ell$ is an internal line
then $\s_{\ell}=\s(\nn_{\ell},j_{\ell})$, i.e.,
$\de_{j_{\ell}}(\oo\cdot\nn_{\ell})=|\oo\cdot\nn-\s_{\ell}\om_{j_{\ell}}|$
(see (\ref{eq:3.1}) for notations);\\
(4) if $v\in V_{0}(\theta)$ then (see Figure \ref{fig:4}) \\
\null\hskip1.truecm 1. $s_{v}=2$;\\
\null\hskip1.truecm 2. both lines $\ell_{1}$ and $\ell_{2}$ entering
$v$ are internal and have $\s_{\ell_{1}}=\s_{\ell_{2}}=\s_{v}$
and $j_{\ell_{1}}=j_{\ell_{2}}=j_{v}$;\\
\null\hskip1.truecm 3. either $\nn_{\ell_{1}}=\s_{v}\ee_{j_{v}}$ 
and $\nn_{\ell_{2}} \neq \s_{v}\ee_{j_{v}}$ or
$\nn_{\ell_{1}} \neq \s_{v}\ee_{j_{v}}$ 
and $\nn_{\ell_{2}} = \s_{v}\ee_{j_{v}}$;\\
\null\hskip1.truecm 4. $\s_{\ell_{v}}=\s_{v}$;\\
(5) if $\ell$ is an internal line and $\nn_{\ell}=\s_{\ell}\ee_{j_{\ell}}$,
then $\ell$ enters a node $v\in V_{0}(\theta)$;\\
(6) $n_{\ell} \ge 0$ if $\nn_{\ell} \neq \s_{\ell}\ee_{j_{\ell}}$ and
$n_{\ell}=-1$ otherwise.
\end{const}
%%%%%%%%%%%%%%%%%%%%%%%%%%%%%%%%%%%%%%%%%%%%%%%%%%%%%%%%%%%%%%%%%%%%%%%%%%

%%%%%%%%%%%%%%%%%%%%%%%%%%%%%%%%%%%%%%%%%%%%%%%%%%%%%%%%%%%%%%%%%%%%%%%%%%
% figure 4
%%%%%%%%%%%%%%%%%%%%%%%%%%%%%%%%%%%%%%%%%%%%%%%%%%%%%%%%%%%%%%%%%%%%%%%%%%
%%%%%%%%%%%%%%%%%%%%%%%%%%%%%%%%%%%%%%%%%%%%%%%%%%%%%%%%%%%%%%%%%%%%%%%%%%

%%%%%%%%%%%%%%%%%%%%%%%%%%%%%%%%%%%%%%%%%%%%%%%%%%%%%%%%%%%%%%%%%%%%%%%%%%
\begin{defi}\label{def:4.5}
Given a tree $\theta$, call $\ell_{0}$ its root line and
consider the internal lines $\ell_{1},\ldots,\ell_{p} \in
L(\theta)$ on scale $-1$ (if any) such that one has
$n_{\ell}\ge 0$ for all $\ell\in\calP(\ell_{0},\ell_{i})$,
$i=1,\ldots,p$; we shall say that $\ell_{1},\ldots,\ell_{p}$
are the lines on scale $-1$ which are closest to the root of $\theta$.
For each such line $\ell_{i}$, call $\theta_{i}=\theta_{\ell{i}}$.
Then we call \emph{pruned tree} $\thetapru$
the subgraph with set of nodes and set of lines
\begin{equation}
N(\thetapru)=N(\theta)\setminus \bigcup_{i=1}^{p} N(\theta_{i}) ,
\qquad L(\thetapru)=L(\theta)\setminus \bigcup_{i=1}^{p} L(\theta_{i}) ,
\nonumber \end{equation}
respectively.
\end{defi}
%%%%%%%%%%%%%%%%%%%%%%%%%%%%%%%%%%%%%%%%%%%%%%%%%%%%%%%%%%%%%%%%%%%%%%%%%%

By construction, $\thetapru$ is a tree,
except that, with respect to the constraints listed above,
one has $s_{v}=1$ whenever $k_{v}=0$; moreover one has
$\nn_{\ell} \neq \s_{\ell}\ee_{j_{\ell}}$ (and hence
$n_{\ell} \ge 0$) for all internal lines $\ell\in L(\thetapru)$
except possibly the root line.

%%%%%%%%%%%%%%%%%%%%%%%%%%%%%%%%%%%%%%%%%%%%%%%%%%%%%%%%%%%%%%%%%%%%%%%%%%
\begin{const}\label{const:4.6}
The modes of the end nodes and the momenta of the lines
are related as follows: if $\ell = \ell_{v}$ one has the
\emph{conservation law}
\begin{equation}
\nn_{\ell} = \sum_{\substack{w\in E(\theta) \\ w \preceq v}} \nn_{w} -
\sum_{\substack{ w \in V_{0}(\theta) \\ w \preceq v}} \s_{w} \ee_{j_{w}}
=\sum_{\substack{w\in E(\thetapru) \\ w \preceq v}} \nn_{w} .
\nonumber\end{equation}
\end{const}
%%%%%%%%%%%%%%%%%%%%%%%%%%%%%%%%%%%%%%%%%%%%%%%%%%%%%%%%%%%%%%%%%%%%%%%%%%

\begin{figure}[!ht]
\begin{center}{
\psfrag{a}{(a)}
\psfrag{b}{(b)}
\psfrag{v}{$v$}
\psfrag{l}{$\ell$}
\psfrag{l1}{$\ell_{1}$}
\psfrag{l2}{$\ell_{2}$}
\psfrag{l3}{$\ell_{3}$}
\psfrag{ls}{$\ell_{s_{v}}$}
\psfrag{jk}{$j_{\ell} \! \, k_{v}$}
\psfrag{j0s}{$j_{\ell} \! \, 0 \, \!\s_{\ell}$}
\psfrag{nsj}{$\nn_{\ell} \, \s_{\ell} \, j_{\ell}$}
\psfrag{nsj1}{$\nn_{\ell_{1}} \, \s_{\ell_{1}} \, j_{\ell_{1}}$}
\psfrag{nsj2}{$\nn_{\ell_{2}} \, \s_{\ell_{2}} \, j_{\ell_{2}}$}
\psfrag{nsj3}{$\nn_{\ell_{3}} \, \s_{\ell_{3}} \, j_{\ell_{3}}$}
\psfrag{nsjs}{$\nn_{\ell_{s_{v}}} \, \s_{\ell_{s_{v}}} \, j_{\ell_{s_{v}}}$}
\psfrag{esj}{$\s_{\ell}\ee_{j_{\ell}} \, \s_{\ell} \, j_{\ell}$}
\includegraphics[width=12cm]{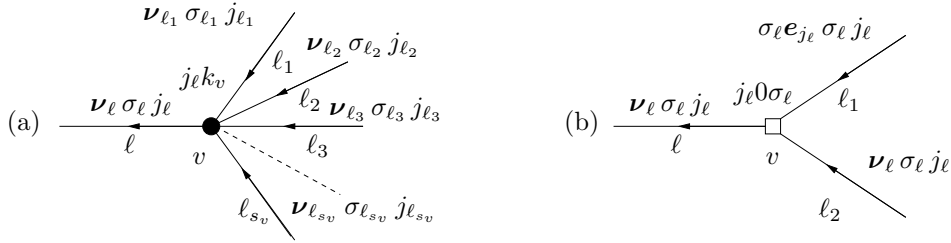}
}
\caption{\label{fig:5}
\footnotesize{Conservation law: (a) $v$ with $k_{v} =s_{v}-1 \ge 1$,
so that $\nn_{\ell}=\nn_{\ell_{1}}+\ldots+\nn_{\ell_{s_{v}}}$,
(b) $v$ with $s_{v}=2$ and $k_{v}=0$.
(The scale labels are not shown).}}
\end{center}
\end{figure}

Note that by Constraint \ref{const:4.6} one has
$\nn_{\ell} = \nn_{v}$ if $v\in E(\theta)$, and
$\nn_{\ell}=\nn_{\ell_{1}}+\ldots+\nn_{\ell_{s_{v}}}$ if
$v \in V(\theta)$, $k_{v}\ge1$, and $\ell_{1},\ldots,\ell_{s_{v}}$
are the lines entering $v$; see Figure \ref{fig:5}.
Moreover for any line $\ell\in L(\theta)$
one has $|\nn_{\ell}|\le |E(\thetapru)|$.

%%%%%%%%%%%%%%%%%%%%%%%%%%%%%%%%%%%%%%%%%%%%%%%%%%%%%%%%%%%%%%%%%%%%%%%%%%
% figure 5
%%%%%%%%%%%%%%%%%%%%%%%%%%%%%%%%%%%%%%%%%%%%%%%%%%%%%%%%%%%%%%%%%%%%%%%%%%
%%%%%%%%%%%%%%%%%%%%%%%%%%%%%%%%%%%%%%%%%%%%%%%%%%%%%%%%%%%%%%%%%%%%%%%%%%

%%%%%%%%%%%%%%%%%%%%%%%%%%%%%%%%%%%%%%%%%%%%%%%%%%%%%%%%%%%%%%%%%%%%%%%%%%
\begin{rmk} \label{rmk:4.7}
In the following we shall repeatedly consider the operation of
changing the sign label of the nodes. Of course this changes
produces the change of other labels, consistently with the
constraints mentioned above: for instance, if we change
the label $\s_{v}$ of an end node $v$ into $-\s_{v}$, then also
$\nn_{v}$ is changed into $-\nn_{v}$; if we change the sign labels
of all the end nodes, then also the momenta of all the lines are
changed, according to the conservation law (Constraint \ref{const:4.6});
and so on.
\end{rmk}
%%%%%%%%%%%%%%%%%%%%%%%%%%%%%%%%%%%%%%%%%%%%%%%%%%%%%%%%%%%%%%%%%%%%%%%%%%

Two unlabelled trees are called \emph{equivalent} if they can be
transformed into each other by continuously deforming the lines
in such a way that they do not cross each other.
We shall call equivalent two trees if the same happens in such a way
that all labels match.

%%%%%%%%%%%%%%%%%%%%%%%%%%%%%%%%%%%%%%%%%%%%%%%%%%%%%%%%%%%%%%%%%%%%%%%%%%
\begin{defi}\label{def:4.8}
We denote by $\gotT^{k}_{j,\nn}$ the set of
inequivalent trees of order $k$ with \emph{tree component} $j$ and
\emph{tree momentum} $\nn$, that is, such that the component label
and the momentum of the root line are $j$ and $\nn$, respectively. 
Finally for $n\ge -1$ define $\gotT^{k}_{j,\nn}(n)$ the set
of trees $\theta\in\gotT^{k}_{j,\nn}$ such that $n_{\ell}\le n$
for all $\ell\in L(\theta)$.
\end{defi}
%%%%%%%%%%%%%%%%%%%%%%%%%%%%%%%%%%%%%%%%%%%%%%%%%%%%%%%%%%%%%%%%%%%%%%%%%%

%%%%%%%%%%%%%%%%%%%%%%%%%%%%%%%%%%%%%%%%%%%%%%%%%%%%%%%%%%%%%%%%%%%%%%%%%%
\begin{rmk} \label{rmk:4.9}
For $\theta\in \gotT^{k}_{j,\nn}$, by writing $\nn=(\n_{1},\ldots,\n_{d})$,
one has $\n_{i}=|E^{+}_{i}(\thetapru)|-|E^{-}_{i}(\thetapru)|$
for $i=1,\ldots,d$. In particular for $\nn=\s\ee_{j}$,
one has $|E^{\s}_{j}(\thetapru)|=|E^{-\s}_{j}(\thetapru)|+1 \ge 1$,
and $|E^{\s}_{j'}(\thetapru)|=|E^{-\s}_{j'}(\thetapru)|$ for all $j'\neq j$.
\end{rmk}
%%%%%%%%%%%%%%%%%%%%%%%%%%%%%%%%%%%%%%%%%%%%%%%%%%%%%%%%%%%%%%%%%%%%%%%%%%

%%%%%%%%%%%%%%%%%%%%%%%%%%%%%%%%%%%%%%%%%%%%%%%%%%%%%%%%%%%%%%%%%%%%%%%%%%
\begin{lemma} \label{lem:4.10}
The number of unlabelled trees $\theta$ with $N$ nodes is
bounded by $4^{N}$. If $k(\theta)=k$ then $|E(\theta)|\le E_{0}k$
and $|V(\theta)|\le V_{0}k$, 
for suitable positive constants $E_{0}$ and $V_{0}$.
\end{lemma}
%%%%%%%%%%%%%%%%%%%%%%%%%%%%%%%%%%%%%%%%%%%%%%%%%%%%%%%%%%%%%%%%%%%%%%%%%%

%%%%%%%%%%%%%%%%%%%%%%%%%%%%%%%%%%%%%%%%%%%%%%%%%%%%%%%%%%%%%%%%%%%%%%%%%%
\prova The bound $|V(\theta)| \le |E(\theta)|-1$ is easily proved
by induction using that $s_{v}\ge2$ for all $v\in V(\theta)$.
So it is enough to bound $|E(\theta)|$. The definition of order
and Remark \ref{rmk:4.1} yield $|E(\theta)|=1+k(\theta)+|V_{0}(\theta)|$,
and the bound $|V_{0}(\theta)| \le 2k(\theta)-1$ immediately
follows by induction on the order of the tree, simply using that
$s_{v}\ge 2$ for $v\in V(\theta)$. Thus, the assertions are proved
with $E_{0}=V_{0}=3$.\EP
%%%%%%%%%%%%%%%%%%%%%%%%%%%%%%%%%%%%%%%%%%%%%%%%%%%%%%%%%%%%%%%%%%%%%%%%%%

%%%%%%%%%%%%%%%%%%%%%%%%%%%%%%%%%%%%%%%%%%%%%%%%%%%%%%%%%%%%%%%%%%%%%%%%%%
\subsection{Tree expansion} \label{sec:4.4}
%%%%%%%%%%%%%%%%%%%%%%%%%%%%%%%%%%%%%%%%%%%%%%%%%%%%%%%%%%%%%%%%%%%%%%%%%%

Now we shall see how to associate with each tree $\theta \in
\gotT^{k}_{j,\nn}$ a contribution to the coefficients $x^{(k)}_{j,\nn}$
and $\h^{(k)}_{j}$ of the power series in (\ref{eq:4.1}).

For all $j=1,\ldots,d$ set $c_{j}^{+}=c_{j}$ and $c_{j}^{-}=c_{j}^{*}$.
We associate with each end node $v\in E(\theta)$ a \emph{node factor}
\begin{equation}
F_{v} := c_{j_{v}}^{\s_{v}} ,
\label{eq:4.2} \end{equation}
and with each internal node $v \in V(\theta)$ a \emph{node factor}
\begin{equation}
F_{v} := \begin{cases}
\displaystyle{
\frac{s_{v,1}!\ldots s_{v,d}!}{s_{v}!} \,
f_{j_{v},s_{v,1},\ldots,s_{v,d}} } ,
& \qquad k_{v} \ge 1 , \\
& \\
\displaystyle{ - \frac{1}{2c_{j_{v}}^{\s_{v}}} } ,
& \qquad k_{v} = 0 ,
\end{cases}
\label{eq:4.3} \end{equation}
where the coefficients $f_{j,s_{1},\ldots,s_{d}}$ are defined
in (\ref{eq:2.2}).

\begin{figure}[!ht]
\begin{center}{
\psfrag{u}{$u$}
\psfrag{psi}{$\psi(u)$}
\psfrag{psin}{$\Psi_{n}(u)$}
\psfrag{g2}{$\frac{\g}{2}$}
\psfrag{58g}{$\frac{5}{8}\g$}
\psfrag{78g}{$\frac{7}{8}\g$}
\psfrag{g}{$\g$}
\psfrag{2n1m}{$2^{-n-1}\g$}
\psfrag{2n}{$2^{\!-\!n}\g$}
\psfrag{2n1p}{$2^{\!-\!n\!+\!1}\g$}
\includegraphics[width=15cm]{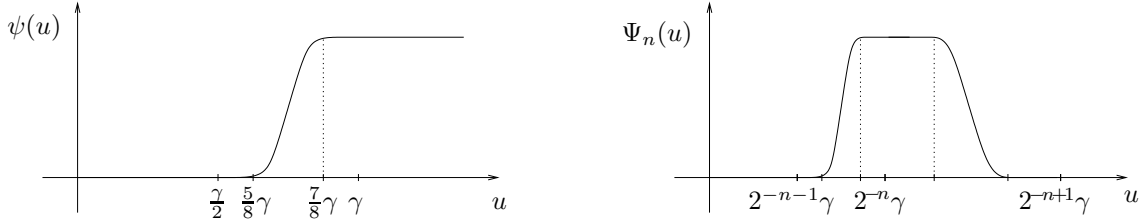}
}
\caption{\label{fig:6}
\footnotesize{The functions $\psi$ and $\Psi_{n}$.}}
\end{center}
\end{figure}

Let $\psi$ be a non-decreasing $C^{\infty}$ function defined
in $\RRR_{+}$, such that (see Figure \ref{fig:6})
\begin{equation}
\psi(u) = \left\{
\begin{array}{ll}
1 , & \text{for } u \geq 7\g/8 , \\
0 , & \text{for } u \leq 5\g/8 ,
\end{array} \right.
\label{eq:4.4} \end{equation}
and set $\chi(u) := 1-\psi(u)$. For all $n \in \ZZZ_{+}$ define
$\chi_{n}(u) := \chi(2^{n}u)$ and $\psi_{n}(u) := \psi(2^{n}u)$,
and set (see Figure \ref{fig:6})
\begin{equation}
\Psi_{n}(u) = \chi_{n-1}(u) \, \psi_{n}(u) ,
\label{eq:4.5} \end{equation}
where $\chi_{-1}(u)=1$. Note that $\chi_{n-1}(u)\chi_{n}(u)=\chi_{n}(u)$,
and hence $\{\Psi_{n}(u)\}_{n\in\ZZZ_{+}}$ is a partition of unity.

%%%%%%%%%%%%%%%%%%%%%%%%%%%%%%%%%%%%%%%%%%%%%%%%%%%%%%%%%%%%%%%%%%%%%%%%%%
% figure 6
%%%%%%%%%%%%%%%%%%%%%%%%%%%%%%%%%%%%%%%%%%%%%%%%%%%%%%%%%%%%%%%%%%%%%%%%%%
%%%%%%%%%%%%%%%%%%%%%%%%%%%%%%%%%%%%%%%%%%%%%%%%%%%%%%%%%%%%%%%%%%%%%%%%%%

We associate with each line $\ell$ a \emph{propagator}
$G_{\ell} := G^{[n_{\ell}]}_{j_{\ell}}(\oo\cdot\nn_{\ell})$, where
\begin{equation}
G^{[n]}_{j}(u) := \begin{cases}\displaystyle{
\frac{\Psi_{n}(\de_{j}(u))}{u^{2}-\om_{j}^{2}},}  & n\ge 0, \\
1, & n=-1.
\end{cases}
\label{eq:4.6} \end{equation}
%

%%%%%%%%%%%%%%%%%%%%%%%%%%%%%%%%%%%%%%%%%%%%%%%%%%%%%%%%%%%%%%%%%%%%%%%%%%
\begin{rmk} \label{rmk:4.11}
The number of scale labels which can be associated with a line
$\ell$ in such a way that $G_{\ell} \neq 0$ is at most $2$.
In particular, given a line $\ell$ with momentum $\nn_{\ell}=\nn$ and
scale $n_{\ell}=n$, such that $\Psi_{n}(\de_{j_{\ell}}(\oo\cdot\nn)) \neq
0$, then (see Figure \ref{fig:6})
\begin{equation}
2^{-(n+1)}\g \le \frac{5}{8} 2^{-n} \g \le
\de_{j_{\ell}}(\oo\cdot\nn) \le
\frac{7}{8} 2^{-(n-1)} \g \le 2^{-(n-1)}\g .
\label{eq:4.7} \end{equation}
and if $\Psi_{n}(\de_{j_{\ell}}(\oo\cdot\nn))
\Psi_{n+1}(\de_{j_{\ell}}(\oo\cdot\nn)) \neq 0$, then
\begin{equation}
\frac{5}{8} \, 2^{-n} \g \le \de_{j_{\ell}}(\oo\cdot\nn)
\le \frac{7}{8} 2^{-n }\g.
\label{eq:4.8} \end{equation}
\end{rmk}
%%%%%%%%%%%%%%%%%%%%%%%%%%%%%%%%%%%%%%%%%%%%%%%%%%%%%%%%%%%%%%%%%%%%%%%%%%

We define
\begin{equation}
\Val(\theta) := \Big(\prod_{\ell \in L(\theta)} G_{\ell} \Big) 
\Big( \prod_{v \in N(\theta)} F_{v} \Big) ,
\label{eq:4.9} 
\end{equation}
and call $\Val(\theta)$ the \emph{value of the tree} $\theta$.

%%%%%%%%%%%%%%%%%%%%%%%%%%%%%%%%%%%%%%%%%%%%%%%%%%%%%%%%%%%%%%%%%%%%%%%%%%
\begin{rmk} \label{rmk:4.12}
The number of trees $\theta\in\gotT^{k}_{j,\nn}$ with $\Val(\theta)\neq 0$
is bounded proportionally to $C^{k}$, for some positive constant $C$.
This immediately follows from Lemma \ref{lem:4.10} and the
observation that the number of trees obtained from
a given unlabelled tree by assigning the labels to the nodes
and the lines is also bounded by a constant to the power $k$
(use Remark \ref{rmk:4.11} to bound the number of allowed scale labels).
\end{rmk}
%%%%%%%%%%%%%%%%%%%%%%%%%%%%%%%%%%%%%%%%%%%%%%%%%%%%%%%%%%%%%%%%%%%%%%%%%%

%%%%%%%%%%%%%%%%%%%%%%%%%%%%%%%%%%%%%%%%%%%%%%%%%%%%%%%%%%%%%%%%%%%%%%%%%%
\begin{rmk} \label{rmk:4.13}
In any tree $\theta$ there is at least one end node with
node factor factor $c_{j}^{\s}$ for each internal node $v$
with $k_{v}=0$, $\s_{v}=\s$ and $j_{v}=j$ (this is easily proved
by induction on the order of the pruned tree):
the node factors $-1/2c_{j}^{\s}$ do not introduce any singularity at
$c_{j}^{\s}=0$. Therefore for any tree $\theta$ the corresponding value
$\Val(\theta)$ is well defined because both propagators and node factors
are finite quantities. Remark \ref{rmk:4.12} implies that also
$$ \sum_{\theta\in\gotT^{k}_{j,\nn}} \Val(\theta) $$
is well defined for all $k \in \NNN$, all $j\in\{1,\ldots,d\}$,
and all $\nn\in\ZZZ^{d}$.
\end{rmk}
%%%%%%%%%%%%%%%%%%%%%%%%%%%%%%%%%%%%%%%%%%%%%%%%%%%%%%%%%%%%%%%%%%%%%%%%%%

%%%%%%%%%%%%%%%%%%%%%%%%%%%%%%%%%%%%%%%%%%%%%%%%%%%%%%%%%%%%%%%%%%%%%%%%%%
\begin{lemma} \label{lem:4.14}
For all $k\in\NNN$, all $j=1,\ldots,d$,
and any $\theta\in\gotT^{k}_{j,\s \ee_{j}}$, there exists
$\theta'\in\gotT^{k}_{j,-\s\ee_{j}}$ such that
$c_{j}^{-\s}\Val(\theta) = c_{j}^{\s}\Val(\theta')$.
The tree $\theta'$ is obtained from $\theta$ by changing the
sign labels of all the nodes $v\in N_{0}(\theta)$.
\end{lemma}
%%%%%%%%%%%%%%%%%%%%%%%%%%%%%%%%%%%%%%%%%%%%%%%%%%%%%%%%%%%%%%%%%%%%%%%%%%

%%%%%%%%%%%%%%%%%%%%%%%%%%%%%%%%%%%%%%%%%%%%%%%%%%%%%%%%%%%%%%%%%%%%%%%%%%
\prova The proof is by induction on the order of the tree.
For any tree $\theta\in\gotT^{k}_{j,\ee_{j}}$ consider the tree
$\theta'\in\gotT^{k}_{j,-\ee_{j}}$ obtained from $\theta$ by replacing
all the labels $\s_{v}$ of all nodes $v \in N_{0}(\theta)$ with
$-\s_{v}$, so that the mode labels $\nn_{v}$ are replaced with $-\nn_{v}$
and the momenta $\nn_{\ell}$ with $-\nn_{\ell}$ (see Remark \ref{rmk:4.7}).
Call $\ell_{1},\ldots,\ell_{p}$ the lines on scale $-1$ (if any)
closest to the root of $\theta$, and for $i=1,\ldots,p$ denote
by $v_{i}$ the node $\ell_{i}$ enters and $\theta_{i}=\theta_{\ell_{i}}$
(recall (2) in Notation \ref{def:4.2}).
As an effect of the change of the sign labels,
each tree $\theta_{i}$ is replaced with a tree $\theta_{i}'$ such that
$c_{j_{v_{i}}}^{-\s}\Val(\theta_{i}) =c_{j_{v_{i}}}^{\s}\Val(\theta_{i}')$,
by the inductive hypothesis. Thus, for each node $v_{i}$ the quantity
$F_{v_{i}}\Val(\theta_{i})$ is not changed.
Moreover, neither the propagators of the lines $\ell\in L(\thetapru)$
nor the node factors corresponding to the internal nodes
$v\in V(\thetapru)$ with $k_{v}\neq0$ change,
while the node factors $c_{j_{v}}^{\s_{v}}$ of the nodes
$v\in E(\thetapru)$ are changed into $c_{j_{v}}^{-\s_{v}}$.
On the other hand one has $|E_{i}^{+}(\thetapru)|=|E_{i}^{-}(\thetapru)|$
for all $i\neq j$, whereas $|E_{j}^{+}(\thetapru)|=
|E_{j}^{-}(\thetapru)|+1$ and $|E_{j}^{+}(\thetapru')|+1=
|E_{j}^{-}(\thetapru')|$. Therefore one obtains
$c_{j}^{-\s} \Val(\theta) = c_{j}^{\s} \Val(\theta')$,
and the assertion follows.\EP
%%%%%%%%%%%%%%%%%%%%%%%%%%%%%%%%%%%%%%%%%%%%%%%%%%%%%%%%%%%%%%%%%%%%%%%%%%

For $k\in\NNN$, $j\in\{1,\ldots,d\}$, and $\s\in\{\pm\}$, define
\begin{equation}
\h^{(k)}_{j,\s} = -\frac{1}{c^{\s}_{j}}
\sum_{\theta \in \gotT^{k}_{j,\s\ee_{j}}} \Val(\theta) .
\nonumber \end{equation}
%

%%%%%%%%%%%%%%%%%%%%%%%%%%%%%%%%%%%%%%%%%%%%%%%%%%%%%%%%%%%%%%%%%%%%%%%%%%
\begin{lemma} \label{lem:4.15}
For all $k\in\NNN$ and all $j=1,\ldots,d$
one has $\h^{(k)}_{j,+}=\h^{(k)}_{j,-}$.
\end{lemma}
%%%%%%%%%%%%%%%%%%%%%%%%%%%%%%%%%%%%%%%%%%%%%%%%%%%%%%%%%%%%%%%%%%%%%%%%%%

%%%%%%%%%%%%%%%%%%%%%%%%%%%%%%%%%%%%%%%%%%%%%%%%%%%%%%%%%%%%%%%%%%%%%%%%%%
\prova Lemma \ref{lem:4.14} implies
\begin{equation}
c_{j}^{-} \sum_{\theta\in\gotT^{k}_{j,\ee_{j}}} \Val(\theta) =
c_{j}^{+} \sum_{\theta\in\gotT^{k}_{j,-\ee_{j}}} \Val(\theta)
\nonumber \end{equation}
for all $k\in\NNN$ and all $j=1,\ldots,d$, so that the assertion follows
from the definition of $\eta^{(k)}_{j,\s}$.\EP
%%%%%%%%%%%%%%%%%%%%%%%%%%%%%%%%%%%%%%%%%%%%%%%%%%%%%%%%%%%%%%%%%%%%%%%%%%

%%%%%%%%%%%%%%%%%%%%%%%%%%%%%%%%%%%%%%%%%%%%%%%%%%%%%%%%%%%%%%%%%%%%%%%%%%
\begin{lemma} \label{lem:4.16}
Equations (\ref{eq:2.7}) formally hold, i.e., they hold
to all perturbation orders, provided that
for all $k\in\NNN$ and $j=1,\ldots,d$ we set formally
\begin{eqnarray}
\displaystyle{  x_{j,\nn}} &=&\displaystyle{\sum_{k=1}^{\io} 
\e^{k} x^{(k)}_{j,\nn}\, , \qquad x^{(k)}_{j,\nn} =
 \sum_{\theta \in \gotT^{k}_{j,\nn}} \!\!\! \Val(\theta)
\quad \forall\nn\in\ZZZ^{d}\setminus\{\pm\ee_{j}\}\,,
\qquad x^{(k)}_{j,\pm\ee_{j}}=0
\,,}
\label{eq:4.10} \\ 
\displaystyle{  \h_{j}}&=&
\displaystyle{\sum_{k=1}^{\io} \e^{k} \h^{(k)}_{j} , \qquad
\h^{(k)}_{j} = -\frac{1}{c_{j}}
\sum_{\theta \in \gotT^{k}_{j,\ee_{j}}} \!\!\! \Val(\theta) .}
\label{eq:4.11} \end{eqnarray}
\end{lemma}
%%%%%%%%%%%%%%%%%%%%%%%%%%%%%%%%%%%%%%%%%%%%%%%%%%%%%%%%%%%%%%%%%%%%%%%%%%

%%%%%%%%%%%%%%%%%%%%%%%%%%%%%%%%%%%%%%%%%%%%%%%%%%%%%%%%%%%%%%%%%%%%%%%%%%
\prova The proof is a direct check.\EP
%%%%%%%%%%%%%%%%%%%%%%%%%%%%%%%%%%%%%%%%%%%%%%%%%%%%%%%%%%%%%%%%%%%%%%%%%%

%%%%%%%%%%%%%%%%%%%%%%%%%%%%%%%%%%%%%%%%%%%%%%%%%%%%%%%%%%%%%%%%%%%%%%%%%%
\begin{rmk} \label{rmk:4.17}
In $\eta^{(k)}_{j}$, defined as (\ref{eq:4.11}), there is no
singularity in $c_{j}=0$ because $\Val(\thetapru)$ contains
at least one factor $c_{j}^{+}=c_{j}$ by Remark \ref{rmk:4.9}.
\end{rmk}
%%%%%%%%%%%%%%%%%%%%%%%%%%%%%%%%%%%%%%%%%%%%%%%%%%%%%%%%%%%%%%%%%%%%%%%%%%

In the light of Lemma \ref{lem:4.16} one can wonder why the definition
of the propagators for $\nn_{\ell}\neq \s_{\ell} \ee_{j_{\ell}}$ is so
involved; as a matter of fact one could define
$$ G_{\ell}=\frac{1}{(\oo\cdot\nn_{\ell})^{2}-\om_{j}^{2}}. $$
However, since $\sum_{n\ge 0}\Psi_{n}(u)\equiv 1$, the two definitions
are equivalent. We use the definition (\ref{eq:4.6}) so that we can
immediately identify the factors $O(2^{n})$ which could prevent
the convergence of the power series (\ref{eq:4.1}).
In what follows we shall make this idea more precise.

%%%%%%%%%%%%%%%%%%%%%%%%%%%%%%%%%%%%%%%%%%%%%%%%%%%%%%%%%%%%%%%%%%%%%%%%%%
\subsection{Clusters} \label{sec:4.5}
%%%%%%%%%%%%%%%%%%%%%%%%%%%%%%%%%%%%%%%%%%%%%%%%%%%%%%%%%%%%%%%%%%%%%%%%%%

\begin{figure}[!ht]
\begin{center}{
\psfrag{1}{$-\!1$}
\psfrag{2}{$2$}
\psfrag{3}{$3$}
\psfrag{5}{$5$}
\psfrag{a}{(a)}
\psfrag{b}{(b)}
\includegraphics[width=16cm]{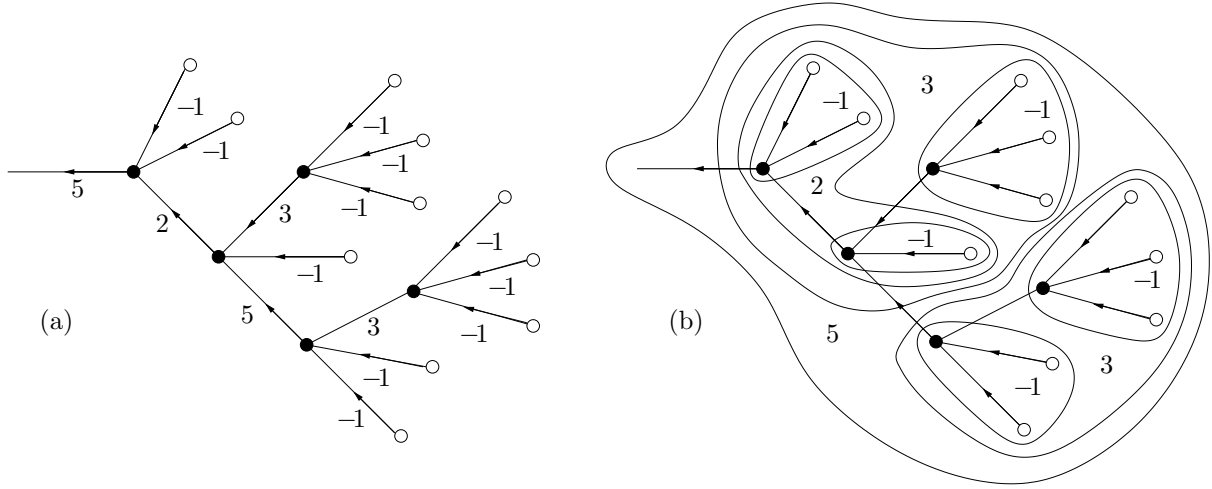}
}
\caption{\label{fig:7}
\footnotesize{Example of tree and the corresponding clusters:
once the scale labels have been assigned to the lines of the tree
as in (a), one obtains the cluster structure depicted in (b).}}
\end{center}
\end{figure}

A \emph{cluster} $T$ on scale $n$ is a maximal set of nodes and lines
connecting them such that all the lines have scales $n'\le n$
and there is at least one line with scale $n$; see Figure \ref{fig:7}.
The lines entering
the cluster $T$ and the line coming out from it (unique if existing
at all) are called the \emph{external} lines of the cluster $T$.
We call $V(T)$, $E(T)$, and $L(T)$ the set of
internal nodes, of end nodes, and of lines of $T$, respectively;
note that the external lines of $T$ do not belong to $L(T)$.
Define also $E^{\s}_{j}(T)$ as the set of end nodes $v\in E(T)$
such that $\s_{v}=\s$ and $j_{v}=j$. By setting
\begin{equation}
k(T) := \sum_{v\in V(T)} k_{v},
\nonumber \end{equation}
we say that the cluster $T$ has \emph{order} $k$ if $k(T)=k$.

%%%%%%%%%%%%%%%%%%%%%%%%%%%%%%%%%%%%%%%%%%%%%%%%%%%%%%%%%%%%%%%%%%%%%%%%%%
% figure 7
%%%%%%%%%%%%%%%%%%%%%%%%%%%%%%%%%%%%%%%%%%%%%%%%%%%%%%%%%%%%%%%%%%%%%%%%%%
%%%%%%%%%%%%%%%%%%%%%%%%%%%%%%%%%%%%%%%%%%%%%%%%%%%%%%%%%%%%%%%%%%%%%%%%%%

%%%%%%%%%%%%%%%%%%%%%%%%%%%%%%%%%%%%%%%%%%%%%%%%%%%%%%%%%%%%%%%%%%%%%%%%%%
\subsection{Self-energy clusters} \label{sec:4.6}
%%%%%%%%%%%%%%%%%%%%%%%%%%%%%%%%%%%%%%%%%%%%%%%%%%%%%%%%%%%%%%%%%%%%%%%%%%

\begin{figure}[!ht]
\begin{center}{
\psfrag{1}{$-\!1$}
\psfrag{2}{$2$}
\psfrag{3}{$3$}
\psfrag{n1}{$\nn_{1}$}
\psfrag{n2}{$\nn_{2}$}
\psfrag{n3}{$\nn_{3}$}
\psfrag{n4}{$\nn_{4}$}
\psfrag{n5}{$\nn_{5}$}
\psfrag{n6}{$\nn_{6}$}
\psfrag{l}{$\ell_{T}$}
\psfrag{l'}{$\ell_{T}'$}
\psfrag{ell}{$\ell$}
\includegraphics[width=7cm]{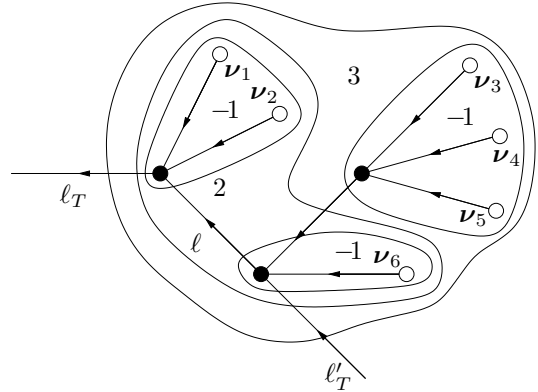}
}
\caption{\label{fig:8}
\footnotesize{Example of self-energy cluster: consider the
cluster $T$ on scale $3$ in Figure \ref{fig:7}, and suppose that the
mode labels of the end nodes are such that $|\nu_{1}+\nu_{2}+\nu_{3}+
\nu_{4}+\nu_{5}+\nu_{6}|\le 2$ and $\de_{j_{\ell_{T}}}(\oo\cdot
\nn_{\ell_{T}})=\de_{j_{\ell_{T}'}}(\oo\cdot\nn_{\ell_{T}'})$.
Then $T$ is a self-energy cluster with external lines
$\ell_{T}'$ (entering line) and $\ell_{T}$ (exiting line).
The path $\calP_{T}$ is such that $\calP_{T}=\{\ell\}$.}}
\end{center}
\end{figure}

We call \emph{self-energy cluster} any cluster $T$ such that
(see Figure \ref{fig:8})\\
(1) $T$ has only one entering line and
one exiting line,\\
(2) one has $n_{\ell} \le\min\{n_{\ell_{T}},n_{\ell_{T}'}\}-2$ for any
$\ell \in L(T)$,\\
(3) one has $|\nn_{\ell_{T}}-\nn_{\ell_{T}'}| \le 2$ and
$\de_{j_{\ell_{T}}}(\oo\cdot\nn_{\ell_{T}})=
\de_{j_{\ell_{T}'}}(\oo\cdot\nn_{\ell_{T}'})$.

%%%%%%%%%%%%%%%%%%%%%%%%%%%%%%%%%%%%%%%%%%%%%%%%%%%%%%%%%%%%%%%%%%%%%%%%%%
\begin{defi}\label{def:4.18}
For any self-energy cluster $T$ we denote by $\ell_{T}$ and $\ell_{T}'$
the exiting and the entering line of $T$ respectively.
We call $\calP_{T}$ the path of lines $\ell\in L(T)$ connecting
$\ell_{T}'$ to $\ell_{T}$, i.e., $\calP_{T}=\calP(\ell_{T},\ell_{T}')$
(recall (1) in Notation \ref{def:4.2}),
and set $n_{T}=\min\{n_{\ell_{T}},n_{\ell_{T}'}\}$.
\end{defi}
%%%%%%%%%%%%%%%%%%%%%%%%%%%%%%%%%%%%%%%%%%%%%%%%%%%%%%%%%%%%%%%%%%%%%%%%%%

%%%%%%%%%%%%%%%%%%%%%%%%%%%%%%%%%%%%%%%%%%%%%%%%%%%%%%%%%%%%%%%%%%%%%%%%%%
\begin{rmk} \label{rmk:4.19}
Notice that, by Remark \ref{rmk:3.3}, for any self-energy
cluster the label $\nn_{\ell_T}$ is  uniquely fixed by the labels
$j_{\ell_{T}},\s_{\ell_{T}},j_{\ell_{T}'},\s_{\ell_{T}'},\nn_{\ell_T'}$.
In particular, fixed $\nn$ and $j$ such that $\de_{j}(\oo\cdot\nn) \le \g$,
there are only $2d-1$ momenta $\nn' \neq \nn$ such that
$|\nn'-\nn|\le 2$ and $\de_{j'}(\oo\cdot\nn')=\de_{j}(\oo\cdot\nn)$
for some $j'$ and $\s'$, depending on $\nn'$. All the other $\nn''$
with small divisor equal to $\de_{j}(\oo\cdot\nn)$ are far away
from $\nn$, according to Lemma \ref{lem:3.1}.
%
% for fixed $\nn_{\ell_T'},j_{\ell_{T}'},\s_{\ell_{T}'}$ there is
%only a finite number of possible choices for $j_{\ell_{T}},\s_{\ell_{T}}$
%and consequently of $\nn_{\ell_T}=\nn_{\ell_T'}+ \s_{\ell_{T}} 
%\ee_{j_{\ell_{T}}} 
%-\s_{\ell_{T}'} \ee_{j_{\ell_{T}'}} $.
\end{rmk}
%%%%%%%%%%%%%%%%%%%%%%%%%%%%%%%%%%%%%%%%%%%%%%%%%%%%%%%%%%%%%%%%%%%%%%%%%%

%%%%%%%%%%%%%%%%%%%%%%%%%%%%%%%%%%%%%%%%%%%%%%%%%%%%%%%%%%%%%%%%%%%%%%%%%%
% figure 8
%%%%%%%%%%%%%%%%%%%%%%%%%%%%%%%%%%%%%%%%%%%%%%%%%%%%%%%%%%%%%%%%%%%%%%%%%%
%%%%%%%%%%%%%%%%%%%%%%%%%%%%%%%%%%%%%%%%%%%%%%%%%%%%%%%%%%%%%%%%%%%%%%%%%%

\begin{figure}[!ht]
\begin{center}{
\psfrag{l}{$\ell$}
\psfrag{T1}{$T_{1}$}
\psfrag{T2}{$T_{2}$}
\includegraphics[width=8cm]{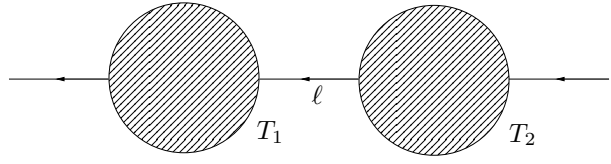}
}
\caption{\label{fig:9}
\footnotesize{Example of resonant line: $\ell$ is resonant if
both $T_{1}$ and $T_{2}$ are self-energy clusters.}}
\end{center}
\end{figure}

We say that a line $\ell$ is a \emph{resonant line} if it is
both the exiting line of a self-energy cluster and the entering line
of another self-energy cluster, that is, $\ell$ is resonant
if there exist two self-energy clusters $T_{1}$ and $T_{2}$ such that
$\ell=\ell_{T_{1}}'=\ell_{T_{2}}$; see Figure \ref{fig:9}.

%%%%%%%%%%%%%%%%%%%%%%%%%%%%%%%%%%%%%%%%%%%%%%%%%%%%%%%%%%%%%%%%%%%%%%%%%%
\begin{rmk} \label{rmk:4.20}
The notion of self-energy cluster was first introduced by Eliasson,
in the context of KAM theorem, in \cite{E}, where it was called
resonance. We prefer the term self-energy cluster to stress further
the analogy with quantum field theory.
\end{rmk}
%%%%%%%%%%%%%%%%%%%%%%%%%%%%%%%%%%%%%%%%%%%%%%%%%%%%%%%%%%%%%%%%%%%%%%%%%%

%%%%%%%%%%%%%%%%%%%%%%%%%%%%%%%%%%%%%%%%%%%%%%%%%%%%%%%%%%%%%%%%%%%%%%%%%%
% figure 9
%%%%%%%%%%%%%%%%%%%%%%%%%%%%%%%%%%%%%%%%%%%%%%%%%%%%%%%%%%%%%%%%%%%%%%%%%%
%%%%%%%%%%%%%%%%%%%%%%%%%%%%%%%%%%%%%%%%%%%%%%%%%%%%%%%%%%%%%%%%%%%%%%%%%%

The notion of equivalence given for trees can be extended in the
obvious way to self-energy clusters. 

\begin{figure}[!ht]
\centering
%\begin{minipage}[b]{8cm}
%\centering
{
\psfrag{njs}{$\nn\,j\,\s\,n'$}
\psfrag{njs'}{$\nn\,j\,\s\,n''$}
\psfrag{se}{$\s\ee_{j}$}
\psfrag{v}{$v_{T}$}
\psfrag{T}{$T$}
\includegraphics[width=6cm]{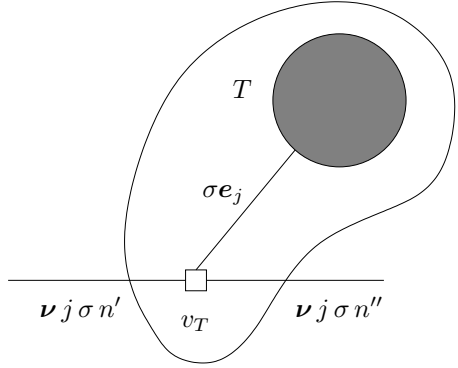}
}
\caption{\label{fig:10}
\footnotesize{A self-energy cluster in
$\gotE^{k}_{j,\s,j,\s}(n)$; $T$ contains at least one line
on scale $\le n$ and $n$ such that $\min\{n',n''\}\ge n+2$.}}
%\end{minipage}
\end{figure}

%%%%%%%%%%%%%%%%%%%%%%%%%%%%%%%%%%%%%%%%%%%%%%%%%%%%%%%%%%%%%%%%%%%%%%%%%%
\begin{defi}\label{def:4.21}
We denote by
$\gotR^{k}_{j,\s,j',\s'}(\oo\cdot\nn',n)$ the set of inequivalent
self-energy clusters $T$ on scale $\le n$ of order $k$, such that
$\nn_{\ell_{T}'}=\nn'$, $j_{\ell_{T}}=j$, $\s_{\ell_{T}}=\s$,
$j_{\ell_{T}'}=j'$ and $\s_{\ell_{T}'}=\s'$. By definition of cluster
for $T\in\gotR^{k}_{j,\s,j',\s'}(\oo\cdot\nn',n)$ one must have 
$n\le n_{T}-2$. For $j=j'$ and $\s=\s'$ define also $\gotE^{k}_{j,\s,j,\s}
(\oo\cdot\nn',n)$ the set of self-energy clusters $T\in
\gotR^{k}_{j,\s,j,\s}(\oo\cdot\nn',n)$ such that
(1) $\ell_{T}'$ enters the same node $v$ which $\ell_{T}$ exits and
(2) $k_{v}=0$. We call $v_{T}$ such a special node and set
$\overline\gotR^{k}_{j,\s,j,\s}(\oo\cdot\nn',n)= \gotR^{k}_{j,\s,j,\s}
(\oo\cdot\nn',n)\setminus \gotE^{k}_{j,\s,j,\s}(\oo\cdot\nn',n) $;
see Figure \ref{fig:10}.
\end{defi}
%%%%%%%%%%%%%%%%%%%%%%%%%%%%%%%%%%%%%%%%%%%%%%%%%%%%%%%%%%%%%%%%%%%%%%%%%%

%%%%%%%%%%%%%%%%%%%%%%%%%%%%%%%%%%%%%%%%%%%%%%%%%%%%%%%%%%%%%%%%%%%%%%%%%%
% figure 10
%%%%%%%%%%%%%%%%%%%%%%%%%%%%%%%%%%%%%%%%%%%%%%%%%%%%%%%%%%%%%%%%%%%%%%%%%%
%%%%%%%%%%%%%%%%%%%%%%%%%%%%%%%%%%%%%%%%%%%%%%%%%%%%%%%%%%%%%%%%%%%%%%%%%%

\begin{figure}[!ht]
\centering
{
\psfrag{njsn}{$\nn\,j\,\s\,n'$}
\psfrag{njsn'}{$\nn\,j\,\s\,n''$}
\psfrag{se}{$\s\ee_{j}$}
\psfrag{se'}{$\s'\ee_{j'}$}
\psfrag{se-}{$-\s\ee_{j}$}
\psfrag{se'-}{$-\s'\ee_{j'}$}
\psfrag{v}{$v_{T}$}
\psfrag{T}{$T$}
\psfrag{th}{$\theta_{T}=$}
\includegraphics[width=14cm]{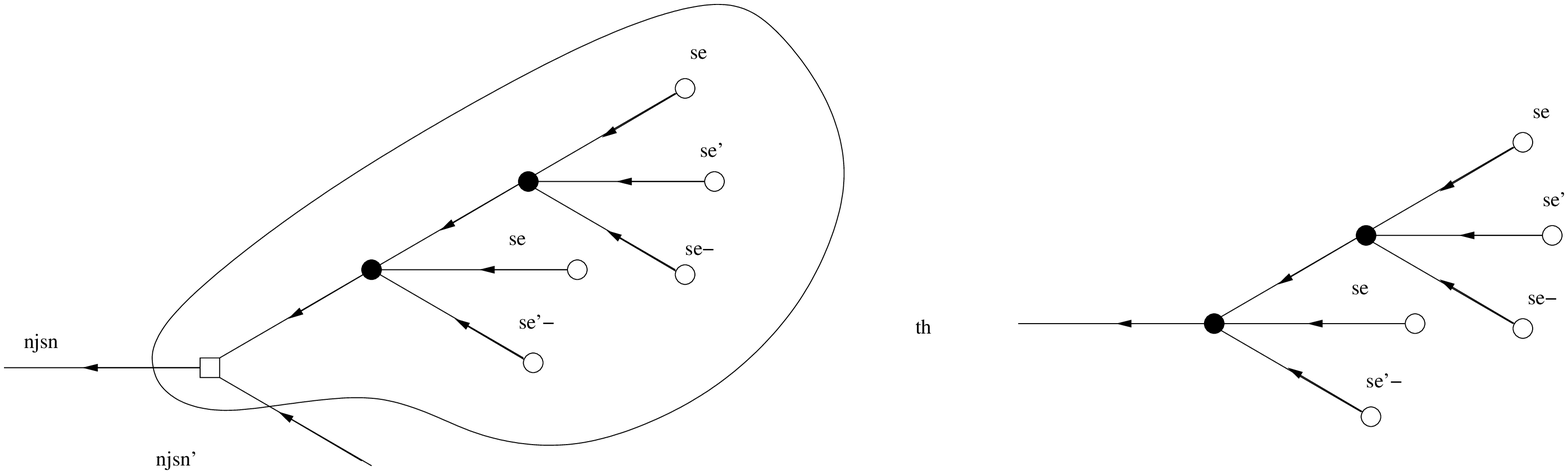}
}
\caption{\label{fig:11}
\footnotesize{An example of self-energy cluster
$T \in \gotE^{k}_{j,\s,j,\s}(n)$ and the corresponding
tree $\theta_{T}$. (Only the mode labels of the end nodes are shown
in $T$ and $\theta_{T}$)}}
\end{figure}

%%%%%%%%%%%%%%%%%%%%%%%%%%%%%%%%%%%%%%%%%%%%%%%%%%%%%%%%%%%%%%%%%%%%%%%%%%
\begin{defi}\label{def:4.22}
For any $T\in\gotE^{k}_{j,\s,j,\s}(\oo\cdot\nn',n)$ we call $\theta_{T}$
the tree which has as root line the line $\ell\in L(T)$
entering $v_{T}$ (one can imagine to obtain $\theta_{T}$ from $T$
by `removing' the node $v_{T}$); see Figure \ref{fig:11}.
Note that  $\theta_{T}\in\gotT^{k}_{j,\s\ee_{j}}(n)$.
\end{defi}
%%%%%%%%%%%%%%%%%%%%%%%%%%%%%%%%%%%%%%%%%%%%%%%%%%%%%%%%%%%%%%%%%%%%%%%%%%

%%%%%%%%%%%%%%%%%%%%%%%%%%%%%%%%%%%%%%%%%%%%%%%%%%%%%%%%%%%%%%%%%%%%%%%%%%
% figure 11
%%%%%%%%%%%%%%%%%%%%%%%%%%%%%%%%%%%%%%%%%%%%%%%%%%%%%%%%%%%%%%%%%%%%%%%%%%
%%%%%%%%%%%%%%%%%%%%%%%%%%%%%%%%%%%%%%%%%%%%%%%%%%%%%%%%%%%%%%%%%%%%%%%%%%

%%%%%%%%%%%%%%%%%%%%%%%%%%%%%%%%%%%%%%%%%%%%%%%%%%%%%%%%%%%%%%%%%%%%%%%%%%
\begin{defi}\label{def:4.23}
Consider a self-energy cluster $T$ such that $n_{\ell}\ge 0$ for all lines
$\ell\in\calP_{T}$.
Call $\ell_{1},\ldots,\ell_{p} \in L(T)$ the internal lines
on scale $-1$ (if any) which are closest
to the exiting line of $T$, that is, such that $n_{\ell}\ge 0$
for all lines $\ell\in\calP(\ell_{T},\ell_{i})$, $i=1,\ldots,p$.
For each line $\ell_{i}$ set $\theta_{i}=\theta_{\ell_{i}}$.
Then the \emph{pruned self-energy cluster}
$\Tpru$ is the subgraph with set of
nodes and set of lines
$$ N(\Tpru)=N(T)\setminus \bigcup_{i=1}^{p} N(\theta_{i}) ,
\qquad L(\Tpru)=L(T)\setminus \bigcup_{i=1}^{p} L(\theta_{i}) , $$
respectively. 
\end{defi}
%%%%%%%%%%%%%%%%%%%%%%%%%%%%%%%%%%%%%%%%%%%%%%%%%%%%%%%%%%%%%%%%%%%%%%%%%%

%%%%%%%%%%%%%%%%%%%%%%%%%%%%%%%%%%%%%%%%%%%%%%%%%%%%%%%%%%%%%%%%%%%%%%%%%%
\begin{rmk} \label{rmk:4.24}
For $T\in\gotR^{k}_{j,\s,j',\s'}(\oo\cdot\nn',n)$ such that
$n_{\ell}\ge 0$ for all $\ell\in\calP_{T}$, one  has
$|E^{+}_{i}(\Tpru)|=|E^{-}_{i}(\Tpru)|$ for all $i\neq j,j'$. 
If $j\neq j'$ then $|E_{j'}^{-\s'}(\Tpru)|=|E^{\s'}_{j'}(\Tpru)|+1$
and $|E_{j}^{\s}(\Tpru)|=|E^{-\s}_{j}(\Tpru)|+1$; if $j=j'$,
$\s=\s'$ and $T\in\overline\gotR^{k}_{j,\s,j,\s}(\oo\cdot\nn',n)$
then $|E_{j}^{\s}(\Tpru)|=|E^{-\s}_{j}(\Tpru)|$, while
if $j=j'$ and $\s=-\s'$ then
$|E_{j}^{\s}(\Tpru)|=|E^{-\s}_{j}(\Tpru)|+2$. 
Finally, for any $T\in\gotE^{k}_{j,\s,j,\s}(\oo\cdot\nn',n)$ one has
$|E_{j}^{\s}(\Tpru)|=|E_{j}^{-\s}(\Tpru)|+1\ge 1$.
\end{rmk}
%%%%%%%%%%%%%%%%%%%%%%%%%%%%%%%%%%%%%%%%%%%%%%%%%%%%%%%%%%%%%%%%%%%%%%%%%%

We shall define
\begin{equation}
\Val(T,\oo\cdot\nn_{\ell_{T}'}) :=
\Big(\prod_{\ell \in L(T)} G_{\ell} \Big)
\Big( \prod_{v \in N(T)} F_{v} \Big) ,
\label{eq:4.12} \end{equation}
\noindent where $\Val(T,\oo\cdot\nn_{\ell_{T}'})$
will be called
the \emph{value of the self-energy cluster} $T$.

The value $\Val(T,\oo\cdot\nn_{\ell_{T}'})$ depends on
$\oo\cdot\nn_{\ell_{T}'}$ through the propagators of the
lines $\ell\in\calP_{T}$. 

%%%%%%%%%%%%%%%%%%%%%%%%%%%%%%%%%%%%%%%%%%%%%%%%%%%%%%%%%%%%%%%%%%%%%%%%%%
\begin{rmk}\label{rmk:4.25}
The value of a self-energy cluster $T\in\gotE^{k}_{j,\s,j,\s}(u,n)$ does
not depend on $u$ so that we shall write 
$$\Val(T,u)=\Val(T)=
-\frac{1}{2c_{j}^{\s}}\Val(\theta_{T}).
$$
\end{rmk}
%%%%%%%%%%%%%%%%%%%%%%%%%%%%%%%%%%%%%%%%%%%%%%%%%%%%%%%%%%%%%%%%%%%%%%%%%%

We define also for future convenience
\begin{equation}
M^{(k)}_{j,\s,j',\s'}(\oo\cdot\nn',n)  := \!\!\!\!\!\!\!\!\!\!\!
\sum_{T \in \gotR^{k}_{j,\s,j',\s'}(\oo\cdot\nn',n)} \!\!\!\!\!\!\!\!\!\!\!
\Val(T,\oo\cdot\nn') ,
\label{eq:4.13} \end{equation}
\vskip-.3truecm
\noindent Note that $M^{(k)}_{j,\s,j,\s}(\oo\cdot\nn',n)=
\widetilde M^{(k)}_{j,\s,j,\s}(n)+
\overline M^{(k)}_{j,\s,j,\s}(\oo\cdot\nn',n)$, where
$\widetilde M^{(k)}_{j,\s,j,\s}(n)$ and
$\overline M^{(k)}_{j,\s,j,\s}(\oo\cdot\nn',n)$ are defined as in
(\ref{eq:4.13}) but for the sum restricted to the set
$\gotE^{k}_{j,\s,j,\s}(\oo\cdot\nn',n)$ and
$\overline\gotR^{k}_{j,\s,j,\s}(\oo\cdot\nn',n)$ respectively.

%%%%%%%%%%%%%%%%%%%%%%%%%%%%%%%%%%%%%%%%%%%%%%%%%%%%%%%%%%%%%%%%%%%%%%%%%%
\begin{rmk} \label{rmk:4.26}
Both the quantities $M^{(k)}_{j,\s,j',\s'}(\oo\cdot\nn',n)$
and the coefficients $x^{(k)}_{j,\nn}$ and $\h^{(k)}_{j}$
are well defined to all orders because the number of terms
which one sums over is finite (by the same argument 
in Remark \ref{rmk:4.12}). At least formally, we can define
\begin{equation}
M_{j,\s,j',\s'}(\oo\cdot\nn')=\sum_{k=1}^{\io}\e^{k}\sum_{n\geq -1}
 M^{(k)}_{j,\s,j',\s'}(\oo\cdot\nn',n) .
\nonumber \end{equation}
\end{rmk}
%%%%%%%%%%%%%%%%%%%%%%%%%%%%%%%%%%%%%%%%%%%%%%%%%%%%%%%%%%%%%%%%%%%%%%%%%%

We define the \emph{depth} $D(T)$ of a self-energy cluster $T$
recursively as follows: we set $D(T)=1$ if there is no
self-energy cluster containing $T$, and set $D(T)=D(T')+1$ if
$T$ is contained inside a self-energy cluster $T'$ and no
other self-energy clusters inside $T'$ (if any) contain $T$.
We denote by $\gotS_{D}(\theta)$ the set of self-energy clusters
of depth $D$ in $\theta$, and by $\gotS_{D}(\theta,T)$ the set of
self-energy clusters of depth $D$ in $\theta$ contained inside $T$.

%%%%%%%%%%%%%%%%%%%%%%%%%%%%%%%%%%%%%%%%%%%%%%%%%%%%%%%%%%%%%%%%%%%%%%%%%%
\begin{defi}\label{def:4.27}
Call $\thetao=\theta\setminus\gotS_{1}(\theta)$ the
subgraph of $\theta$ formed by the set of nodes and lines of $\theta$
which are outside the set $\gotS_{1}(\theta)$
(the external lines of the self-energy clusters $T\in\gotS_{1}(\theta)$
being included in $\thetao$), and, analogously, for
$T\in \gotS_{D}(\theta)$ call $\To=T\setminus\gotS_{D+1}(\theta,T)$ the
subgraph of $T$ formed by the set of nodes
and lines of $T$ which are outside the set $\gotS_{D+1}(\theta,T)$.
We denote by $V(\To)$, $E(\To)$, and $L(\To)$ the set of internal nodes,
of end nodes, and of lines of $\To$, and by
$k(\To)$ the order of $\To$, that is, the
sum of the labels $k_{v}$ of all the internal nodes $v\in V(\To)$.
\end{defi}
%%%%%%%%%%%%%%%%%%%%%%%%%%%%%%%%%%%%%%%%%%%%%%%%%%%%%%%%%%%%%%%%%%%%%%%%%%

%%%%%%%%%%%%%%%%%%%%%%%%%%%%%%%%%%%%%%%%%%%%%%%%%%%%%%%%%%%%%%%%%%%%%%%%%%
\begin{lemma} \label{lem:4.28}
Given a line $\ell\in L(\theta)$, if $T$ is the self-energy cluster
with largest depth containing $\ell$ (if any), $\ell\in\calP_{T}$
and there is no line $\ell' \in \calP_{T}$ preceding $\ell$
with $n_{\ell'}=-1$, one can write
$\nn_{\ell}=\nn_{\ell}^{0}+\nn_{\ell_{T}'}$. Then one has
$|\nn_{\ell}^{0}|\le E_{1}k(\To)$, for a suitable positive constant $E_{1}$,
if $k(\To)\ge 1$, and $|\nn_{\ell}^{0}|\le 2$ if $k(\To)=0$.
\end{lemma}
%%%%%%%%%%%%%%%%%%%%%%%%%%%%%%%%%%%%%%%%%%%%%%%%%%%%%%%%%%%%%%%%%%%%%%%%%%

%%%%%%%%%%%%%%%%%%%%%%%%%%%%%%%%%%%%%%%%%%%%%%%%%%%%%%%%%%%%%%%%%%%%%%%%%%
\prova We first prove that for any tree $\theta$, if we
denote by $\ell_{0}$ its root line, one has
\begin{equation}
|\nn_{\ell_{0}}| \le \begin{cases}
E_{1} k(\thetao)-2 , & \text{if } \ell_{0}
\text{ does not exit a self-energy cluster} , \\
E_{1} k(\thetao) , & \text{if } \ell_{0}
\text{ exits a self-energy cluster} ,
\end{cases}
\label{eq:4.14} \end{equation}
for a suitable constant $E_{1}\ge 4$. The proof is by induction on the
order of the tree $\theta$.
If $k(\theta)=1$ (and hence $\thetao=\theta$) then the only
internal line of $\theta$ is $\ell_{0}$ and $|\nn_{\ell_{0}}| \le 2$,
so that the assertion trivially holds provided $E_{1}\ge 4$.
If $k(\theta)>1$ let $v_{0}$ be the node which $\ell_{0}$ exits.
If $v_{0}$ is not contained inside a self-energy cluster let
$\ell_{1},\ldots,\ell_{m}$, $m\ge 0$, be the internal lines
entering $v_{0}$ and $\theta_{i}=\theta_{\ell_{i}}$ for all $i=1,\ldots,m$.
Finally let $\ell_{m+1},\ldots,\ell_{m+m'}$ be the end-lines entering $v_0$.
By definition we have $k(\thetao)=k_{v_{0}}+k(\thetao_{1})+\ldots+
k(\thetao_{m})$. If $k_{v_{0}}>0$,
we have  $\nn_{\ell_{0}}=\nn_{\ell_{1}}+\ldots+\nn_{\ell_{m+m'}}$.
This implies in turn
$$
|\nn_{\ell_{0}} |\leq |\nn_{\ell_{1}}|+\ldots+|\nn_{\ell_{m}}|+m' 
\leq E_1\big(k(\thetao_{1})+\ldots+
k(\thetao_{m})\big)+m' \leq E_1(k(\thetao)-m-m'+1)+m'
$$
The assertion follows for $E_1\ge 4$ by the inductive hypothesis
(the worst possible case is $m=0$, $m'=2$).

If  $k_{v_{0}}=0$ then  $s_v= 2$ and $m'=0$. moreover one of the lines,
say $\ell_1$, is on scale $n=-1$ while for the other line one has
$\nn_{\ell_{0}}=\nn_{\ell_{2}}$. Once more the bound follows
from the inductive hypothesis
since $|\nn_{\ell_{2}}|\le E_1k(\thetao_{2})$ $\le E_1(k(\thetao)-1)$. 

Finally, if $v_{0}$ is contained inside a self-energy cluster, then
$\ell_{0}$ exits a self-energy cluster $T_{1}$. There will be
$p$ self-energy clusters $T_{1},\ldots,T_{p}$, $p\ge 1$,
such that the exiting line of $T_{i}$ is the entering line
of $T_{i-1}$, for $i=2,\ldots,p$, while the entering line
$\ell'$ of $T_{p}$ does not exit any self-energy cluster.
By Lemma \ref{lem:3.4}, one has $|\nn_{\ell_{0}}-\nn_{\ell'}| \le 2$
and $k(\thetao)=k(\thetao_{\ell'})$. Then,
by the inductive hypothesis, one finds $|\nn_{\ell_{0}}| \le 2 +
E_{1}k(\thetao_{\ell'})-2 =E_{1}k(\thetao)$.

\begin{figure}[!ht]
\centering
{
\psfrag{lT}{$\ell_{T}$}
\psfrag{lT'}{$\ell_{T}'$}
\psfrag{T1}{$T_{1}$}
\psfrag{Tp}{$T_{p}$}
\psfrag{v}{$v$}
\psfrag{l0}{$\ell_{0}$}
\psfrag{l}{$\ell$}
\psfrag{l1}{$\ell_{1}$}
\psfrag{l2}{$\ell_{2}$}
\psfrag{l1'}{$\ell_{1}'$}
\psfrag{l2'}{$\ell_{2}'$}
\psfrag{l3'}{$\ell_{3}'$}
\psfrag{th1}{$\theta_{1}$}
\psfrag{th2}{$\theta_{2}$}
\includegraphics[width=14cm]{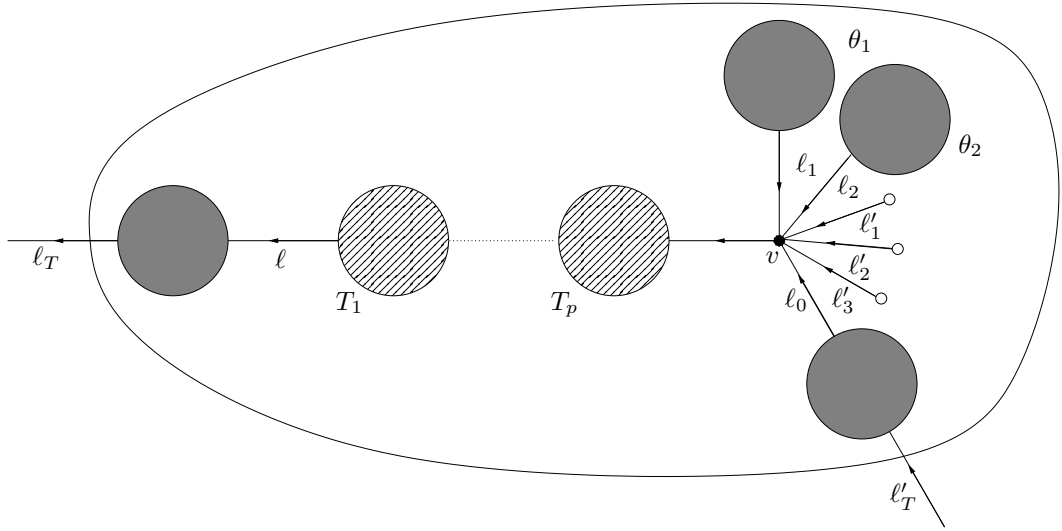}
}
\caption{\label{fig:12}
\footnotesize{The self-energy cluster $T$
considered in the proof of Lemma \ref{lem:4.28},
with $m=2$, $m'=3$, and a chain of $p$ self-energy clusters
between $\ell$ and $\ell_{v}$ (one has $p\ge 0$,
and $\ell=\ell_{v}$ if $p=0$).}}
\end{figure}

Now for $\ell$ and $T$ as in the statement we prove,
by induction on the order of the self-energy cluster, the bound
\begin{equation}
|\nn_{\ell}^{0}| \le \begin{cases}
E_{1} k(\To_{\ell})-2 , & \text{if } k(\To_{\ell}) \ge 1 , \\
2 & \text{if } k(\To_{\ell}) = 0 ,
\end{cases}
\label{eq:4.15} \end{equation}
where $\To_{\ell}$ is the set of nodes and lines of $\To$ which precede
$\ell$. The bound is trivially satisfied when $k(\To_{\ell})=0$.
Otherwise let $v$ be the node in $V(\To)$ between
$\ell$ and $\ell_{T}'$ which is closest to $\ell$.
If $k_{v}=0$ the bound follows trivially
by using the bound (\ref{eq:4.14}). If $k_{v}\ge 1$ call $\ell_{1},
\ldots,\ell_{m}$, $m\ge 0$, the internal lines entering $v$
which are not along the path $\calP_{T}$, and
$\ell_{m+1},\ldots,\ell_{m+m'}$ the end lines entering $v$;
one has $m+m' \ge 1$. There is a further line
$\ell_{0}\in\calP_{T}$ entering $v$ such that
$\nn_{\ell_{0}}=\nn_{\ell_{0}}^{0}+\nn_{\ell_{T}'}$; see Figure \ref{fig:12}.
Using also  Lemma \ref{lem:3.4} one has
$|\nn_{\ell}^{0}| \le 2 + |\nn_{\ell_{0}}^{0}| +
|\nn_{\ell_{1}}|+\ldots+|\nn_{\ell_{m}}|+m'$. As $n_{\ell_{0}}\le
n_{\ell_{T}'}-2$ one has $k(\To_{\ell_{0}})\ge 1$ and hence,
by (\ref{eq:4.14}) and the inductive hypothesis, one has
\begin{equation}
|\nn_{\ell}^{0}| \le 2 + \left(E_{1} k(\To_{\ell_{0}})-2\right) + 
E_{1}\left( k(\thetao_{1}) + \ldots k(\thetao_{m})\right)+m' ,
\nonumber \end{equation}
where $\theta_{i}=\theta_{\ell_{i}}$ for all $i=1,\ldots,m$.
Thus, since  $k(\To_{\ell_{0}})+k(\thetao_{1})+\ldots+
k(\thetao_{m})+(m+m')=k(\To_{\ell})$ and $m+m'\ge 1$, one finds
\begin{equation}
|\nn_{\ell}^{0}| \le E_{1} \left(k(\To_{\ell}) - m-m'\right) +m' \le
E_{2} k(\To_{\ell}) - 2 ,
\nonumber \end{equation}
provided $E_{1} \ge 4$.
Therefore, the assertion follows with, say, $E_{1}=4$.\EP
%%%%%%%%%%%%%%%%%%%%%%%%%%%%%%%%%%%%%%%%%%%%%%%%%%%%%%%%%%%%%%%%%%%%%%%%%%

%%%%%%%%%%%%%%%%%%%%%%%%%%%%%%%%%%%%%%%%%%%%%%%%%%%%%%%%%%%%%%%%%%%%%%%%%%
% figure 12
%%%%%%%%%%%%%%%%%%%%%%%%%%%%%%%%%%%%%%%%%%%%%%%%%%%%%%%%%%%%%%%%%%%%%%%%%%
%%%%%%%%%%%%%%%%%%%%%%%%%%%%%%%%%%%%%%%%%%%%%%%%%%%%%%%%%%%%%%%%%%%%%%%%%%

%%%%%%%%%%%%%%%%%%%%%%%%%%%%%%%%%%%%%%%%%%%%%%%%%%%%%%%%%%%%%%%%%%%%%%%%%%
\begin{defi}\label{def:4.29}
Given a tree $\theta$ and a line $\ell \in L(\theta)$, call
$\Ga_{\ell}=\Ga_{\ell}(\theta)$
the subgraph formed by the set of nodes and lines
which do not precede $\ell$; see figure.
Let us call $\Gao_{\ell}$ the set of nodes and lines of $\Ga_{\ell}$
which are outside any self-energy cluster contained inside $\Ga_{\ell}$.
\end{defi}
%%%%%%%%%%%%%%%%%%%%%%%%%%%%%%%%%%%%%%%%%%%%%%%%%%%%%%%%%%%%%%%%%%%%%%%%%%

%%%%%%%%%%%%%%%%%%%%%%%%%%%%%%%%%%%%%%%%%%%%%%%%%%%%%%%%%%%%%%%%%%%%%%%%%%
% figure 12
%%%%%%%%%%%%%%%%%%%%%%%%%%%%%%%%%%%%%%%%%%%%%%%%%%%%%%%%%%%%%%%%%%%%%%%%%%
\begin{figure}[!ht]
\centering
{
\psfrag{th}{$\theta=$}
\psfrag{l}{$\ell$}
\psfrag{G}{$\Ga_{\ell}$}
\psfrag{H}{$\theta_{\ell}$}
\includegraphics[width=8cm]{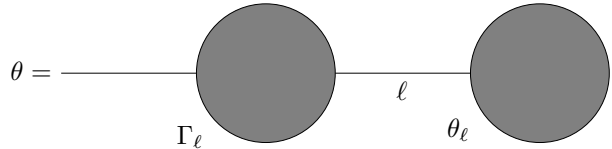}
}
\caption{\label{fig:13}
\footnotesize{The set $\Ga_{\ell}=\Ga_{\ell}(\theta)$ and the subtree
$\theta_{\ell}$ determined by the line $\ell\in L(\theta)$.
If $\ell$ is the root line then $\Ga_{\ell}=\emptyset$.}}
\end{figure}
%%%%%%%%%%%%%%%%%%%%%%%%%%%%%%%%%%%%%%%%%%%%%%%%%%%%%%%%%%%%%%%%%%%%%%%%%%

%%%%%%%%%%%%%%%%%%%%%%%%%%%%%%%%%%%%%%%%%%%%%%%%%%%%%%%%%%%%%%%%%%%%%%%%%%
\begin{lemma} \label{lem:4.30}
Given a tree $\theta$ let $\ell_{0}$  and $\ell$ be the root line and an
arbitrary internal line preceding $\ell_{0}$. If $k(\Gao_{\ell})\ge 1$
one has $|\nn_{\ell_{0}}-\nn_{\ell}| \le E_{2} k(\Gao_{\ell})$,
for a suitable positive constant $E_{2}$.
\end{lemma}
%%%%%%%%%%%%%%%%%%%%%%%%%%%%%%%%%%%%%%%%%%%%%%%%%%%%%%%%%%%%%%%%%%%%%%%%%%

%%%%%%%%%%%%%%%%%%%%%%%%%%%%%%%%%%%%%%%%%%%%%%%%%%%%%%%%%%%%%%%%%%%%%%%%%%
\prova We prove by induction on the order of $\Ga_{\ell}$ the bound
\begin{equation}
|\nn_{\ell_{0}}-\nn_{\ell}| \le \begin{cases}
E_{2} k(\Gao_{\ell})- 2 , & \text{if } \ell_{0}
\text{ does not exit a self-energy cluster} , \\
E_{2} k(\Gao_{\ell}) , & \text{if } \ell_{0}
\text{ exits a self-energy cluster} .
\end{cases}
\label{eq:4.16} \end{equation}
We mimic the proof of (\ref{eq:4.14}) in Lemma \ref{lem:4.28}.
The case $k(\Gao_{\ell})=1$ is trivial provided $E_{2}\ge 3$, so let us
consider $k(\Gao_{\ell})>1$ and call $v_{0}$ the node which $\ell_{0}$
exits. If $v_{0}$ is not contained inside a self-energy cluster and
$k_{v_{0}}\ge 1$ then $\nn_{\ell_{0}}=\nn_{\ell_{1}}+\ldots+\nn_{\ell_{m+m'}}$,
where $\ell_{1},\ldots,\ell_{m}$ are the internal lines
entering $v_{0}$, with (say) $\ell_{m}\in \calP(\ell_{0},\ell)\cup\{\ell\}$,
and $\ell_{m+1},\ldots,\ell_{m+m'}$ are the end lines
entering $v_{0}$.
Hence $k(\Gao_{\ell})=k_{v_{0}}+k(\thetao_{1})+\ldots+
k(\thetao_{m-1})+k(\Gao_{m})$, where  $\theta_{i}=\theta_{\ell_{i}}$
and $\Ga_{m}=\Ga_{\ell}(\theta_{\ell_{m}})$
($\Ga_{\ell_{m}}=\emptyset$ if $\ell_{m}=\ell$). Thus, the assertion follows
by (\ref{eq:4.14}) and the inductive hypothesis.
If $v_{0}$ is not contained inside a self-energy cluster
and $k_{v_{0}}=0$ then two lines $\ell_{1}$ and $\ell_{2}$ enter
$v_{0}$, and one of them, say $\ell_{1}$, is such that
$|\nn_{\ell_{1}}|=1$. If $\ell=\ell_{2}$ the result is trivial.
If $\ell_{2} \in \calP(\ell_{0},\ell)$
the bound follows once more from the inductive hypothesis.
If $\ell=\ell_{1}$ one has
\begin{equation}
|\nn_{\ell_{0}} - \nn_{\ell}| \le |\nn_{\ell_{0}}| + 1 \le
E_{1} k(\thetao_{2}) + 1 \le E_{2} k(\Gao_{\ell})-2 ,
\nonumber \end{equation}
where $\theta_{2}=\theta_{\ell_{2}}$, provided $E_{2} \ge E_{1}+3$,
if $E_{1}$ is the constant defined in Lemma \ref{lem:4.28}.
If $\ell_{1}\in \calP(\ell_{0},\ell)$ denote by $\ell_{1}'$ the line
on scale $-1$ along the path $\{\ell_{1}\}\cup\calP(\ell_{1},\ell)$
which is closest to $\ell$. Again call $\theta_{2}=\theta_{\ell_{2}}$ 
and $J_{1}$ the subgraph formed by the set of nodes
and lines preceding $\ell_{1}'$ (with $\ell_{1}'$ included)
but not $\ell$; define also $\theta_{1}$ as the tree obtained
from $J_{1}$ by (1) reverting the arrows of all lines
along $\{\ell_{1}',\ell\}\cup\calP(\ell_{1}',\ell)$, (2) replacing
$\ell_{1}'$ with an end line carrying the same sign and component labels
as $\ell_{1}'$, and (3) replacing all the labels $\s_{v}$, $v\in 
N_{0}(J_{1})$ with $-\s_{v}$. One has, by using also (\ref{eq:4.14}),
\begin{equation}
|\nn_{\ell_{0}} - \nn_{\ell}| \le |\nn_{\ell_{0}}| + |\nn_{\ell}| \le
E_{1} k(\thetao_{1}) + E_{1} k(\thetao_{2}) \le E_{2} k(\Gao_{\ell})-2 ,
\nonumber \end{equation}
provided $E_{2} \ge E_{1}+2$ so that the bound follows once more.
Finally, if $v_{0}$ is contained inside a self-energy cluster, then
$\ell_{0}$ exits a self-energy cluster $T_{1}$. There will be
$p$ self-energy clusters $T_{1},\ldots,T_{p}$, $p\ge 1$,
such that the exiting line of $T_{i}$ is the entering line
of $T_{i-1}$, for $i=2,\ldots,p$, while the entering line
$\ell'$ of $T_{p}$ does not exit any self-energy cluster.
By Lemma \ref{lem:3.4}, one has $|\nn_{\ell_{0}}-\nn_{\ell'}| \le 2$
and $k(\Gao_{\ell})=k(\Gao')$,
where $\Ga'=\Ga_{\ell}(\theta_{\ell'})$. Then,
the inductive hypothesis yields $|\nn_{\ell_{0}}-\nn_{\ell}|
\le 2 + |\nn_{\ell'}-\nn_{\ell}| \le 2 + E_{2}k(\Gao')-2 =E_{2}k(\Gao)$.
Therefore the assertion follows with, say, $E_{2}=E_{1}+3$
(and hence $E_{2}=7$ if $E_{1}=4$).\EP
%%%%%%%%%%%%%%%%%%%%%%%%%%%%%%%%%%%%%%%%%%%%%%%%%%%%%%%%%%%%%%%%%%%%%%%%%%

%%%%%%%%%%%%%%%%%%%%%%%%%%%%%%%%%%%%%%%%%%%%%%%%%%%%%%%%%%%%%%%%%%%%%%%%%%
\begin{rmk} \label{rmk:4.31}
Lemma \ref{lem:4.28} will be used in Section \ref{sec:6} to control
the change of the momenta as an effect of the regularisation
procedure (to be defined). Furthermore, both Lemmas \ref{lem:4.28}
and \ref{lem:4.30} will be used in Section \ref{sec:8} to show that
the resonant lines which are not regularised cannot accumulate too much.
\end{rmk}
%%%%%%%%%%%%%%%%%%%%%%%%%%%%%%%%%%%%%%%%%%%%%%%%%%%%%%%%%%%%%%%%%%%%%%%%%%

%%%%%%%%%%%%%%%%%%%%%%%%%%%%%%%%%%%%%%%%%%%%%%%%%%%%%%%%%%%%%%%%%%%%%%%%%%
%%%%%%%%%%%%%%%%%%%%%%%%%%%%%%%%%%%%%%%%%%%%%%%%%%%%%%%%%%%%%%%%%%%%%%%%%%
\zerarcounters
\section{Dimensional bounds}
\label{sec:5}
%%%%%%%%%%%%%%%%%%%%%%%%%%%%%%%%%%%%%%%%%%%%%%%%%%%%%%%%%%%%%%%%%%%%%%%%%%
%%%%%%%%%%%%%%%%%%%%%%%%%%%%%%%%%%%%%%%%%%%%%%%%%%%%%%%%%%%%%%%%%%%%%%%%%%

In this section we discuss how to prove that the series (\ref{eq:4.10})
and (\ref{eq:4.11}) converge if the resonant lines are excluded.
We shall see in the following sections
how to take into account the presence of the resonant lines.

Call $\gotN_{n}(\theta)$ the number of non-resonant lines
$\ell\in L(\theta)$ such that $n_{\ell}\ge n$,
and $\gotN_{n}(T)$ the number of non-resonant lines
$\ell\in L(T)$ such that $n_{\ell}\ge n$.

The analyticity assumption on $f$ yields that one has
\begin{equation}
\left| F_{v} \right| \le \Phi^{s_{v}+k_{v}}
\qquad \forall v \in V(\theta) \setminus V_{0}(\theta) ,
\label{eq:5.1} \end{equation}
for a suitable positive constant $\Phi$.

%%%%%%%%%%%%%%%%%%%%%%%%%%%%%%%%%%%%%%%%%%%%%%%%%%%%%%%%%%%%%%%%%%%%%%%%%%
\begin{lemma} \label{lem:5.1}
Assume that $2^{-(n_{\ell}+2)}\g \le \de_{j_{\ell}}(\oo\cdot\nn_{\ell})\le
2^{-(n_{\ell}-2)}\g$ for all trees $\theta$ and all lines $\ell\in L(\theta)$.
Then there exists a positive constant $c$ such that for any tree $\theta$
one has $\gotN_{n}(\theta) \le c \, 2^{-n/\tau}k(\theta)$.
\end{lemma}
%%%%%%%%%%%%%%%%%%%%%%%%%%%%%%%%%%%%%%%%%%%%%%%%%%%%%%%%%%%%%%%%%%%%%%%%%%

%%%%%%%%%%%%%%%%%%%%%%%%%%%%%%%%%%%%%%%%%%%%%%%%%%%%%%%%%%%%%%%%%%%%%%%%%%
\prova We prove that $\gotN_{n}(\theta) \le \max \{ 0, c \, 2^{-n/\tau}
k(\theta)-2\}$ by induction on the order of $\theta$.\\
1. First of all note that for a tree $\theta$ to have a line $\ell$
on scale $n_{\ell}\ge n$ one needs $k(\theta)\ge k_{n}=
E_{0}^{-1}2^{(n-2)/\tau}$, as it follows from the
Diophantine condition (\ref{eq:3.2a}) and Lemma \ref{lem:4.10}.
Hence the bound is trivially true for $k<k_{n}$.\\
2. For $k(\theta)\ge k_{n}$, let $\ell_{0}$ be the root line of $\theta$
and set $\nn=\nn_{\ell_{0}}$ and $j=j_{\ell_{0}}$.
If $n_{\ell_{0}} < n$ the assertion follows
from the inductive hypothesis. If $n_{\ell_{0}}\ge n$,
call $\ell_{1},\ldots,\ell_{m}$ the lines with scale $\ge n-1$
which are closest to $\ell_{0}$ (that is, such that $n_{\ell}\le n-2$
for all $p=1,\ldots,m$ and all lines $\ell\in\calP(\ell_{1},\ell_{p})$).
The case $m=0$ is trivial.
If $m\ge 2$ the bound follows once more from the inductive hypothesis.\\
3. If $m=1$, then $\ell_{1}$ is the only entering line of a cluster $T$.
Set $\nn'=\nn_{\ell_{1}}$, $j'=j_{\ell_{1}}$ and $n'=n_{\ell_{1}}$.
By hypothesis one has $\de_{j}(\oo\cdot\nn)
\le 2^{-(n-2)}\g$ and $\de_{j'}(\oo\cdot\nn')\le 2^{-(n-3)}\g$,
so that, by Lemma \ref{lem:3.2}, either $|\nn-\nn'|>2^{(n-5)/\tau}$
or $|\nn-\nn'|\le 2$ and $\de_{j}(\oo\cdot\nn)=\de_{j'}(\oo\cdot\nn')$.
In the first case, since
\begin{equation}
\nn-\nn' = \sum_{w\in E(T)} \nn_{w} - 
\sum_{\substack{ w \in V(T) \\ k_{w}=0}} \s_{w} \ee_{j_{w}} =
\sum_{w\in E(\Tpru)} \nn_{w} ,
\nonumber \end{equation}
the same argument used to prove Lemma \ref{lem:4.10} yields
$|\nn-\nn'| \le |E(T)| \le E_{0} k(T)$, and hence $k(T) \ge
E_{0}^{-1}2^{(n-5)/\tau}$. Thus, if $\theta_{1}=\theta_{\ell_{1}}$,
one has $k(\theta) = k(T)+k(\theta_{1})$,
so that
\begin{equation}
\gotN_{n}(\theta) = 1 + \gotN_{n}(\theta_{1}) \le
c \, 2^{-n/\tau}k(\theta_{1})-1 \le
c \, 2^{-n/\tau}k(\theta) - c \, 2^{-n/\tau}k(T)-1
\le c \, 2^{-n/\tau}k(\theta) - 2 ,
\nonumber \end{equation}
provided $c\ge E_{0}2^{5/\tau}$.\\
4. If instead $|\nn-\nn'|\le 2$ and $\de_{j}(\oo\cdot\nn)=
\de_{j'}(\oo\cdot\nn')$, then the only way for $T$ not to be a
self-energy cluster is that $n_{\ell_{1}}=n_{\ell_{0}}-1=n-1$ and there
is at least a line $\ell\in T$ with $n_{\ell}=n-2$.
But then $\de_{j}(\oo\cdot\nn)\neq\de_{j_{\ell}}(\oo\cdot\nn_{\ell})$
so that $|\nn-\nn_{\ell}|>2^{(n-6)/\tau}$ and we can reason as in
the previous case provided $c\ge E_{0}2^{6/\tau}$.
Otherwise $T$ is a self-energy cluster and $\ell_{1}$
can be either resonant or not-resonant.
Call $\ell_{1}',\ldots,\ell_{m'}'$ the lines with scale $\ge n-1$ which
are closest to $\ell_{1}$. Once more the cases $m'=0$ and $m'\ge 2$
are trivial.\\
5. If $m'=1$, then $\ell_{1}'$ is the only entering line of a cluster
$T'$. If $\theta_{1}'=\theta_{\ell_{1}'}$, then
$\gotN_{n}(\theta)=1+\gotN_{n}(\theta_{1}')$ if $\ell_{1}$ is resonant
and $\gotN_{n}(\theta)\le 2+\gotN_{n}(\theta_{1}')$ if $\ell_{1}$ is
non-resonant. Consider first the case of $\ell_{1}$ being non-resonant.
Set $\nn''=\nn_{\ell_{1}'}$, $j''=j_{\ell_{1}'}$ and
$n''=n_{\ell_{1}'}$. By reasoning as before we find that one has either
$|\nn'-\nn''|>2^{(n-5)/\tau}$ or $|\nn'-\nn''|\le 2$ and
$\de_{j'}(\oo\cdot\nn')=\de_{j''}(\oo\cdot\nn'')$.
If $|\nn'-\nn''|>2^{(n-5)/\tau}$ then
$k(T') \ge E_{0}^{-1}2^{(n-5)/\tau}$; thus, by using that
$k(\theta)=k(T)+k(T')+k(\theta_{1}')$, we obtain
\begin{eqnarray}
\gotN_{n}(\theta) & \!\! \le \!\! & 2 + \gotN_{n}(\theta_{1}') \le
c \, 2^{-n/\tau}k(\theta) - c \, 2^{-n/\tau}k(T) -
c \, 2^{-n/\tau}k(T') \nonumber \\
& \!\! \le \!\! & 
c \, 2^{-n/\tau}k(\theta) - c \, 2^{-n/\tau}k(T') \le
c \, 2^{-n/\tau}k(\theta) - 2 ,
\nonumber \end{eqnarray}
provided $c\ge 2 E_{0}2^{5/\tau}$.\\
6. Otherwise one has $|\nn-\nn'|\le 2$, $|\nn'-\nn''| \le 2$, and
$\de_{j}(\oo\cdot\nn) = \de_{j'}(\oo\cdot\nn') =
\de_{j''}(\oo\cdot\nn'')$. Since we are assuming $\ell_{1}$
to be non-resonant then, $T'$ is not a self-energy cluster.
But then there is at least a line $\ell'\in T$ with $n_{\ell'}=n-2$ and
we can reason as in item 4.\\
7. So we are left with the case in which $\ell_{1}$ is resonant and
hence $T'$ is a self-energy cluster. Let $\ell_{1}'$ be the entering
line of $T'$. Once more $\ell_{1}'$ is either resonant or non-resonant.
If it is non-resonant we repeat the same argument as done before for
$\ell_{1}$. If it is resonant, we iterate the construction, and so on.
Therefore we proceed until either we find a non-resonant line on
scale $\ge n$, for which we can reason as before, or we reach a tree
$\theta'$ of order so small that it cannot contain any line on scale
$\ge n$ (i.e., $k(\theta')<k_{n}$).\\
8. Therefore the assertion follows with, say, $c=2E_{0}2^{6/\tau}$.\EP
%%%%%%%%%%%%%%%%%%%%%%%%%%%%%%%%%%%%%%%%%%%%%%%%%%%%%%%%%%%%%%%%%%%%%%%%%%

%%%%%%%%%%%%%%%%%%%%%%%%%%%%%%%%%%%%%%%%%%%%%%%%%%%%%%%%%%%%%%%%%%%%%%%%%%
\begin{rmk} \label{rmk:5.2}
One can wonder why in Lemma \ref{lem:5.1} did we assume
$2^{-(n_{\ell}+2)}\g \le \de_{j_{\ell}}(\oo\cdot\nn_{\ell})\le
2^{-(n_{\ell}-2)}\g$ when Remark \ref{rmk:4.11} assures the stronger
condition $2^{-(n_{\ell}+1)}\g \le \de_{j_{\ell}}(\oo\cdot\nn_{\ell})\le
2^{-(n_{\ell}-1)}\g$. The reason is that later on we shall need
to slightly change the momenta of the lines, in such a way
that the scales in general no longer satisfy the condition
(\ref{eq:4.7}) noted in Remark \ref{rmk:4.11}. However the condition
assumed for proving Lemma \ref{lem:5.1} will still be satisfied.
\end{rmk}
%%%%%%%%%%%%%%%%%%%%%%%%%%%%%%%%%%%%%%%%%%%%%%%%%%%%%%%%%%%%%%%%%%%%%%%%%%

For any tree $\theta$ we call $L_{\rm R}(\theta)$ and $L_{\rm NR}(\theta)$
the sets of resonant lines and of non-resonant line, respectively,
in $L(\theta)$. Then we can write
\begin{equation}
\Val(\theta) = \Big(\prod_{\ell \in L_{\rm R}(\theta)} G_{\ell} \Big) 
\Val_{\rm NR}(\theta) , \qquad
\Val_{\rm NR}(\theta) :=
\Big(\prod_{\ell \in L_{\rm NR}(\theta)} G_{\ell} \Big) 
\Big( \prod_{v \in N(\theta)} F_{v} \Big) ,
\label{eq:5.2} \end{equation}
where each propagator $G_{\ell}$ can be bounded as $C_{0}2^{n_{\ell}}$,
for some constant $C_{0}$.

%%%%%%%%%%%%%%%%%%%%%%%%%%%%%%%%%%%%%%%%%%%%%%%%%%%%%%%%%%%%%%%%%%%%%%%%%%
\begin{lemma} \label{lem:5.3}
For all tree $\theta$ with $k(\theta)=k$ one has
$|\Val_{\rm NR}(\theta)| \le C^{k}\Gamma^{3k}(\cc)$, where
$\Gamma(\cc):=$ $\max\{|c_{1}|,\ldots,|c_{d}|,1\}$ and $C$
is a suitable positive constant.
\end{lemma}
%%%%%%%%%%%%%%%%%%%%%%%%%%%%%%%%%%%%%%%%%%%%%%%%%%%%%%%%%%%%%%%%%%%%%%%%%%

%%%%%%%%%%%%%%%%%%%%%%%%%%%%%%%%%%%%%%%%%%%%%%%%%%%%%%%%%%%%%%%%%%%%%%%%%%
\prova One has
\begin{eqnarray}
\left| \Val_{\rm NR}(\theta) \right|
& \!\!\! \le \!\!\! &
C_{0}^{k} \Gamma^{3k}(\cc) \Phi^{k} 
\Big(\prod_{\ell \in L_{NR}(\theta)} 2^{n_{\ell}} \Big) 
\le C_{0}^{k} \Gamma^{3k}(\cc) \Phi^{k}
\prod_{n=0}^{\io} 2^{n \gotN_{n}(\theta)}
\nonumber
\\
& \!\!\! \le \!\!\! &
C_{0}^{k} \Gamma^{3k}(\cc) \Phi^{k}
\exp \left( c \, \log 2 \, k \sum_{n=1}^{\io}
2^{-n/\tau} n \right) .
\nonumber \end{eqnarray}
The last sum converges: this is enough to prove the lemma.\EP
%%%%%%%%%%%%%%%%%%%%%%%%%%%%%%%%%%%%%%%%%%%%%%%%%%%%%%%%%%%%%%%%%%%%%%%%%%

So far the only bound that we have on the propagators of the
resonant lines is $|G_{\ell}| \le 1/\om_{j_{\ell}}\de_{j_{\ell}}
(\oo\cdot\nn_{\ell})\le C_{0}2^{n_{\ell}}$. What we need is to obtain a
gain factor proportional to $2^{-n_{\ell}}$ for each resonant line $\ell$
with $n_{\ell}\ge 1$.

%%%%%%%%%%%%%%%%%%%%%%%%%%%%%%%%%%%%%%%%%%%%%%%%%%%%%%%%%%%%%%%%%%%%%%%%%%
\begin{lemma} \label{lem:5.4}
Given $\theta$ such that $\Val(\theta)\neq 0$, let $\ell\in L(\theta)$ be
a resonant line and let $T$ be the
self-energy cluster of largest depth containing $\ell$ (if any). Then
there is at least one non-resonant line in $T$ on scale $\ge n_{\ell}-1$.
\end{lemma}
%%%%%%%%%%%%%%%%%%%%%%%%%%%%%%%%%%%%%%%%%%%%%%%%%%%%%%%%%%%%%%%%%%%%%%%%%%

%%%%%%%%%%%%%%%%%%%%%%%%%%%%%%%%%%%%%%%%%%%%%%%%%%%%%%%%%%%%%%%%%%%%%%%%%%
\prova Set $n=n_{\ell}$. There are in general $p\ge 2$ self-energy
clusters $T_{1},\ldots,T_{p}$, contained inside $T$, connected
by resonant lines $\ell_{1},\ldots,\ell_{p-1}$, and $\ell$ is one
of such lines, while the entering line $\ell_{p}$ of $T_{p}$
and the exiting line $\ell_{0}$ of $T_{1}$ are non-resonant.
Moreover $\de(\oo\cdot\nn_{\ell_{i}})=\de(\oo\cdot\nn_{\ell})$
for all $i=0,\ldots,p$, so that all the lines
$\ell_{0},\ldots,\ell_{p}$ have scales either $n,n-1$ or $n,n+1$,
by Remark \ref{rmk:4.11}. In any case the lines $\ell_{0},\ell_{p}$ must
be in $T$ by definition of self-energy cluster. \EP
%%%%%%%%%%%%%%%%%%%%%%%%%%%%%%%%%%%%%%%%%%%%%%%%%%%%%%%%%%%%%%%%%%%%%%%%%%

%%%%%%%%%%%%%%%%%%%%%%%%%%%%%%%%%%%%%%%%%%%%%%%%%%%%%%%%%%%%%%%%%%%%%%%%%%
%%%%%%%%%%%%%%%%%%%%%%%%%%%%%%%%%%%%%%%%%%%%%%%%%%%%%%%%%%%%%%%%%%%%%%%%%%
\zerarcounters
\section{Renormalisation}
\label{sec:6}
%%%%%%%%%%%%%%%%%%%%%%%%%%%%%%%%%%%%%%%%%%%%%%%%%%%%%%%%%%%%%%%%%%%%%%%%%%
%%%%%%%%%%%%%%%%%%%%%%%%%%%%%%%%%%%%%%%%%%%%%%%%%%%%%%%%%%%%%%%%%%%%%%%%%%

\begin{figure}[!ht]
\begin{center}{
\psfrag{T1}{$T_{1}$}
\psfrag{T2}{$T_{2}$}
\psfrag{T3}{$T_{3}$}
\psfrag{Tp}{$T_{p}$}
\includegraphics[width=14cm]{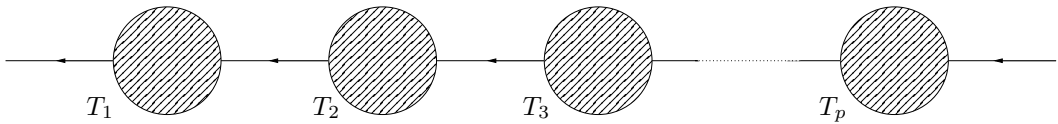}
}
\caption{\label{fig:14}
\footnotesize{A chain of self-energy clusters.}}
\end{center}
\end{figure}

Now we shall see how to deal with the resonant lines. In principle,
one can have trees containing chains of arbitrarily many
self-energy clusters (see Figure \ref{fig:14}), and this produces
an accumulation of small divisors, and hence a bound proportional
to $k!$ to some positive power for the corresponding values. 

%%%%%%%%%%%%%%%%%%%%%%%%%%%%%%%%%%%%%%%%%%%%%%%%%%%%%%%%%%%%%%%%%%%%%%%%%%
% figure 14
%%%%%%%%%%%%%%%%%%%%%%%%%%%%%%%%%%%%%%%%%%%%%%%%%%%%%%%%%%%%%%%%%%%%%%%%%%
%%%%%%%%%%%%%%%%%%%%%%%%%%%%%%%%%%%%%%%%%%%%%%%%%%%%%%%%%%%%%%%%%%%%%%%%%%

Let $K_{0}$ be such that $E_{1}K_{0}=2^{-8/\tau}$. For $T\in
\gotR^{k}_{j,\s,j',\s'}(u,n)$,
define the \emph{localisation operator} $\matL$ by setting
\begin{equation}
\matL \Val(T,u) :=
\begin{cases}
\Val(T,\s'\om_{j'}) , &
k(\To) \le K_{0}2^{n_{T}/\tau} , \;
n_{\ell} \ge 0 \; \forall \ell \in \calP_{T} , \\
& \\
0 , & \mbox{otherwise}
\end{cases}
\label{eq:6.1} \end{equation}
which will be called the \emph{localised value} of the
self-energy cluster $T$. Define also $\matR:=\uno-\matL$, by setting,
for $T\in\gotR^{k}_{j,\s,j',\s'}(u,n)$,
\begin{equation}
\matR \Val(T,u) =
\begin{cases}
\displaystyle{ \left( u-\s'\om_{j'} \right)
\!\!\int_{0}^{1} \!\!\!\!{\rm d}t \, \partial_{u}
\Val(T,\s'\om_{j'}+t(u \! - \! \s'\om_{j'})) } , &
k(\To) \le K_{0}2^{n_{T}/\tau} , \;
n_{\ell} \ge 0 \; \forall \ell \in \calP_{T} , \\
& \\
\Val(T,u) , & \mbox{otherwise}
\end{cases}
\label{eq:6.2} \end{equation}
so that
\begin{subequations}
\begin{align}
\matL M^{(k)}_{j,\s,j',\s'}(u,n) & = \!\!\!\!\!\!\!
\sum_{T \in \gotR^{k}_{j,\s,j',\s'}(u,n)}
\!\!\!\!\!\!\! \calL\Val(T,u)  ,
\label{eq:6.3a} \\
\matR M^{(k)}_{j,\s,j',\s'}(u,n) & = \!\!\!\!\!\!\!
\sum_{T \in \gotR^{k}_{j,\s,j',\s'}(u,n)}
\!\!\!\!\!\!\! \matR\Val(T,u) .
\label{eq:6.3b}
\end{align}
\label{eq:6.3} \end{subequations}
\vskip-.3truecm
\noindent We shall call $\matR$ the \emph{regularisation} operator and
$\matR \Val(T,u)$ the \emph{regularised value} of $T$.

%%%%%%%%%%%%%%%%%%%%%%%%%%%%%%%%%%%%%%%%%%%%%%%%%%%%%%%%%%%%%%%%%%%%%%%%%%
\begin{rmk} \label{rmk:6.1}
If $T\in\gotE^{k}_{j,\s,j,\s}(u,n)$ the localisation operator
acts as
\begin{equation}
\matL \Val(T) =
\begin{cases}
\Val(T) , & k(\To) \le K_{0}2^{n_{T}/\tau} , \\
& \\
0 , & k(\To) > K_{0}2^{n_{T}/\tau} .
\end{cases}
\nonumber \end{equation}
\end{rmk}
%%%%%%%%%%%%%%%%%%%%%%%%%%%%%%%%%%%%%%%%%%%%%%%%%%%%%%%%%%%%%%%%%%%%%%%%%%

%%%%%%%%%%%%%%%%%%%%%%%%%%%%%%%%%%%%%%%%%%%%%%%%%%%%%%%%%%%%%%%%%%%%%%%%%%
\begin{rmk} \label{rmk:6.2}
If in a self-energy cluster $T$ there is a line $\ell\in\calP_{T}$ such
that $\nn_{\ell}=\s_{\ell}\ee_{j_{\ell}}$ (and hence $n_{\ell}=-1$) then
$\matL\Val(T',u)=0$ for all self-energy clusters containing $T$
such that $\ell\in\calP_{T'}$.
\end{rmk}
%%%%%%%%%%%%%%%%%%%%%%%%%%%%%%%%%%%%%%%%%%%%%%%%%%%%%%%%%%%%%%%%%%%%%%%%%%

Recall the definition of the sets $\gotS_{D}(\theta)$ and
$\gotS_{D}(\theta,T)$ after Remark \ref{rmk:4.26}.
For any tree $\theta$ we can write its value as
\begin{equation}
\Val(\theta) = \Big( \prod_{T\in\gotS_{1}(\theta)}
\Val(T,\oo\cdot\nn_{\ell_{T}'}) \Big)
\Big(\prod_{\ell \in L(\theta\setminus\gotS_{1}(\theta))} G_{\ell} \Big) 
\Big( \prod_{v \in N(\theta\setminus\gotS_{1}(\theta))} F_{v} \Big) ,
\label{eq:6.4} \end{equation}
and, recursively, for any self-energy cluster $T$ of depth $D$ we have
\begin{equation}
\Val(T,\oo\cdot\nn_{\ell_{T}'}) =
\Big( \prod_{T'\in\gotS_{D+1}(\theta,T)}
\Val(T',\oo\cdot\nn_{\ell_{T'}'}) \Big)
\Big(\prod_{\ell \in L(T\setminus\gotS_{D+1}(\theta,T))} G_{\ell} \Big) 
\Big( \prod_{v \in N(T\setminus\gotS_{D+1}(\theta,T))} F_{v} \Big) .
\label{eq:6.5} \end{equation}

Then we modify the diagrammatic rules given in Section \ref{sec:4}
by assigning a further label $\matO_{T}\in\{\matR,\matL\}$,
which will  be called the \emph{operator} label,
to each self-energy cluster $T$. Then, by writing $\Val(\theta)$
according to (\ref{eq:6.4}) and (\ref{eq:6.5}), one replaces
$\Val(T,\oo\cdot\nn_{\ell_{T}'})$ with
$\matL\Val(T,\oo\cdot\nn_{\ell_{T}'})$ if $\matO_{T}=\matL$ and with
$\matR\Val(T,\oo\cdot\nn_{\ell_{T}'})$ if $\matO_{T}=\matR$.
When considering the regularised value of a self-energy cluster
$T\in \gotR^{k}_{j,\s,j',\s'}(u,n)$ with $k(\To) \le
K_{0}2^{n_{T}/\tau}$ and $n_{\ell}\ge0$ for all $\ell\in\calP_{T}$,
then we have also an interpolation parameter $t$ to
consider: we shall denote it by $t_{T}$ to keep trace of the self-energy
cluster which it is associated with. We set $t_{T}=1$ for a regularised
self-energy cluster $T$ with either $k(\To)>K_{0}2^{n_{T}/\tau}$ or
$\calP_{T}$ containing at least one line $\ell$ with $n_{\ell}=-1$.

We call \emph{renormalised trees} the trees $\theta$ carrying the
further labels $\matO_{T}$, associated with the self-energy clusters
$T$ of $\theta$. As an effect of the localisation and regularisation
operators the arguments of the propagators of some lines are changed.

%%%%%%%%%%%%%%%%%%%%%%%%%%%%%%%%%%%%%%%%%%%%%%%%%%%%%%%%%%%%%%%%%%%%%%%%%%
\begin{rmk} \label{rmk:6.3}
For any self-energy cluster $T$ the localised value $\matL\Val(T,u)$
does not depend on the operator labels of the self-energy clusters
containing $T$.
\end{rmk}
%%%%%%%%%%%%%%%%%%%%%%%%%%%%%%%%%%%%%%%%%%%%%%%%%%%%%%%%%%%%%%%%%%%%%%%%%%

Given a self-energy cluster $T\in\gotR^{k}_{j,\s,j',\s'}(u,n)$ such
that no line along $\calP_{T}$ is on scale $-1$, let $\ell$
be a line such that (1) $\ell\in\calP_{T}$, and (2) $T$ is the
self-energy cluster with largest depth containing $\ell$. If one has
$\matO_{T}=\matR$, then the quantity $\oo\cdot\nn_{\ell}$ is changed
according to the operator labels of all the self-energy clusters $T'$
such that (1) $T'$ contains $T$, (2) no line along $\calP_{T'}$
has scale $-1$, and (3) $\ell\in\calP_{T'}$.
Call $T_{p} \subset T_{p-1} \subset \ldots
\subset T_{1}$ such self-energy clusters, with $T_{p}=T$. If
$\matO_{T_{i}}=\matR$ for all $i=1,\ldots,p$ then $\oo\cdot\nn_{\ell}$
is replaced with
\begin{eqnarray}
\oo\cdot\nn_{\ell}(\underline{t}_{\ell})
& \!\!\! = \!\!\! &
\oo\cdot\nn_{\ell}^{0} + \s_{\ell_{p}}\om_{j_{p}} +
t_{p} \left( \oo\cdot\nn_{\ell_{p}}^{0} + \s_{\ell_{p-1}}\om_{j_{p-1}} -
\s_{\ell_{p}}\om_{j_{p}} \right)
\nonumber \\
& \!\!\! + \!\!\! & \sum_{i=2}^{p-1} t_{p} \ldots t_{i} 
\left( \oo\cdot\nn_{\ell_{i}}^{0} + \s_{\ell_{i-1}}\om_{j_{i-1}} -
\s_{\ell_{i}}\om_{j_{i}} \right) +
t_{p} \ldots t_{1} \left( \oo\cdot\nn_{\ell_{1}} -
\s_{\ell_{1}}\om_{j_{1}} \right) ,
\label{eq:6.6} \end{eqnarray}
where we have set $\underline{t_{\ell}}=(t_{1},\ldots,t_{p})$,
$\ell_{T_{i}}'=\ell_{i}$ and $t_{T_{i}}=t_{i}$ for simplicity.

Otherwise let $T_{q}$ be the self-energy cluster of highest
depth, among $T_{1},\ldots,T_{p-1}$, with $\matO_{T_{q}}=\matL$
(so that $\matO_{T_{i}}=\matR$ for $i\ge q+1$). In that case,
instead of (\ref{eq:6.6}), one has
\begin{eqnarray}
\oo\cdot\nn_{\ell}(\underline{t}_{\ell})
& \!\!\! = \!\!\! &
\oo\cdot\nn_{\ell}^{0} +
\s_{\ell_{p}}\om_{j_{p}} +
t_{p} \left( \oo\cdot\nn_{\ell_{p}}^{0} + \s_{\ell_{p-1}}\om_{j_{p-1}} -
\s_{\ell_{p}}\om_{j_{p}} \right)
\nonumber \\
& \!\!\! + \!\!\! & 
\sum_{i=q+1}^{p-1} t_{p} \ldots t_{i} 
\left( \oo\cdot\nn_{\ell_{i}}^{0} + \s_{\ell_{i-1}}\om_{j_{i-1}} -
\s_{\ell_{i}}\om_{j_{i}} \right) ,
\label{eq:6.7} \end{eqnarray}
with the same notations used in (\ref{eq:6.6}).

If $\matO_{T_{p}}=\matL$, since $\oo\cdot\nn_{\ell}$ is replaced
with $\oo\cdot\nn_{\ell}^{0} + \s_{\ell_{T}'}\om_{j_{\ell_{T}'}}$
for $\ell\in\calP_{T}$, we can write
$\oo\cdot\nn_{\ell}^{0}+\s_{\ell_{T}'}\om_{j_{\ell_{T}'}}$
as in (\ref{eq:6.6}) by setting $t_{p}=0$.
More generally, if we set $t_{T}=0$ whenever $\matO_{T}=\matL$,
we see that we can always claim that, under the action of the
localisation and regularisation operators, the momentum $\nn_{\ell}$
of any line $\ell \in \calP_{T}$ is changed to $\nn_{\ell}
(\underline{t}_{\ell})$, in such a way that $\oo\cdot\nn_{\ell}
(\underline{t}_{\ell})$ is given by (\ref{eq:6.6}).

%%%%%%%%%%%%%%%%%%%%%%%%%%%%%%%%%%%%%%%%%%%%%%%%%%%%%%%%%%%%%%%%%%%%%%%%%%
\begin{lemma} \label{lem:6.4}
Given $\theta$ such that $\Val(\theta)\neq 0$,
for all $\ell\in L(\theta)$ one has
$4 \, \de_{j_{\ell}}(\oo\cdot\nn_{\ell})\le 5 \,
\de_{j_{\ell}}(\oo\cdot\nn_{\ell}(\underline{t}_{\ell}))
\le 6\,\de_{j_{\ell}}(\oo\cdot\nn_{\ell})$.
\end{lemma}
%%%%%%%%%%%%%%%%%%%%%%%%%%%%%%%%%%%%%%%%%%%%%%%%%%%%%%%%%%%%%%%%%%%%%%%%%%

%%%%%%%%%%%%%%%%%%%%%%%%%%%%%%%%%%%%%%%%%%%%%%%%%%%%%%%%%%%%%%%%%%%%%%%%%%
\prova The proof is by induction on the depth
of the self-energy cluster.\\
1. Consider first the case that $\ell\in \calP_{T}$, with $\matO_{T}=
\matL$. Set $n=n_{\ell_{T}'}$, $\nn'=\nn_{\ell_{T}'}$, $\s'=\s_{\ell_{T}'}$,
and $j'=j_{\ell_{T}'}$. Then $\oo\cdot\nn'$ is replaced with
$\s'\om_{j'}$, and, as a consequence, $\oo\cdot\nn_{\ell}$
is replaced with $\oo\cdot\nn_{\ell}(\underline{t}_{\ell})=
\oo\cdot\nn_{\ell}^{0}+\s'\om_{j'}$. Define $\tilde n_{\ell}$ such that 
\begin{equation}
2^{-(\tilde n_{\ell}+1)}\g \le
\de_{j_{\ell}}(\oo\cdot\nn_{\ell}^{0}+\s'\om_{j'})
\le 2^{-(\tilde n_{\ell}-1)}\g ,
\label{eq:6.8} \end{equation}
where $\de_{j_{\ell}}(\oo\cdot\nn_{\ell}^{0}+\s'\om_{j'}) =
|\oo\cdot\nn_{\ell}^{0}+\s'\om_{j'}-\s_{\ell}\om_{j_{\ell}}|
\ge \g |\nn_{\ell}^{0}|^{-\tau}$ by the 
Diophantine condition (\ref{eq:3.2b}). Therefore
$2^{\tilde n_{\ell}-1} \le |\nn_{\ell}^{0}|^{\tau} \le
(E_{1}k(\To))^{\tau} \le (E_{1}K_{0})^{\tau} 2^{n}=2^{n-8}$,
and hence $\tilde n_{\ell} \le n-7$. Since $| \oo\cdot\nn' -\s'\om_{j'}|
\le 2^{-n+2}\g$ by the inductive hypothesis, one has
\begin{eqnarray}
\de_{j_{\ell}}(\oo\cdot\nn_{\ell})
\!\!\! & = & \!\!\!
\left| \oo\cdot\nn_{\ell}^{0}+\oo\cdot\nn'-\s_{\ell}\om_{j_{\ell}} \right|
\nonumber \\
\!\!\! & \ge & \!\!\!
\left| \oo\cdot\nn_{\ell}^{0}+\s'\om_{j'}-\s_{\ell}\om_{j_{\ell}} \right|
- \left| \oo\cdot\nn' -\s'\om_{j'} \right| \ge
\frac{15}{16} \, \de_{j_{\ell}}(\oo\cdot\nn_{\ell}^{0}+\s'\om_{j'}) ,
\nonumber \end{eqnarray}
because $\de_{j_{\ell}}(\oo\cdot\nn_{\ell}^{0}+\s'\om_{j'}) \ge
2^{-(\tilde n_{\ell}+1)}\g \ge 2^{-n+6}\g \ge 2^{4}
\, | \oo\cdot\nn' -\s'\om_{j'}|$.
In the same way one can bound $\de_{j_{\ell}}(\oo\cdot\nn_{\ell}) \le
| \oo\cdot\nn_{\ell}^{0}+\s'\om_{j'}-\s_{\ell}\om_{j_{\ell}} | +
| \oo\cdot\nn' -\s'\om_{j'}|$, so that we conclude that
\begin{equation}
\frac{15}{16} \, \de_{j_{\ell}}(\oo\cdot\nn_{\ell}^{0}+\s'\om_{j'}) 
\le \de_{j_{\ell}}(\oo\cdot\nn_{\ell}) \le
\frac{17}{16} \, \de_{j_{\ell}}(\oo\cdot\nn_{\ell}^{0}+\s'\om_{j'}) .
\label{eq:6.9} \end{equation}
This yields the assertion.\\
2. Consider now the case that $\matO_{T}=\matR$. In that case
$\oo\cdot\nn_{\ell}(\underline{t}_{\ell})$ is given by (\ref{eq:6.6}).
Define $\tilde n_{\ell}$ as in (\ref{eq:6.8}),
with $\s'=\s_{\ell_{p}}$ and $j'=j_{\ell_{p}}$.
We want to prove that
\begin{equation}
\frac{7}{8} \, \de_{j_{\ell}}(\oo\cdot\nn_{\ell}^{0}+\s'\om_{j'}) 
\le \de_{j_{\ell}}(\oo\cdot\nn_{\ell}(\underline{t}_{\ell})) \le
\frac{9}{8} \, \de_{j_{\ell}}(\oo\cdot\nn_{\ell}^{0}+\s'\om_{j'}) .
\label{eq:6.10} \end{equation}
for all $\underline{t}_{\ell}=(t_{1},\ldots,t_{p})$, with
$t_{i}\in[0,1]$ for $i=1,\ldots,p$. This immediately implies
the assertion because, by using also (\ref{eq:6.9}), we obtain
\begin{equation}
\frac{14}{17} \, \de_{j_{\ell}}(\oo\cdot\nn_{\ell}) \le
\frac{7}{8} \, \de_{j_{\ell}}(\oo\cdot\nn_{\ell}^{0}+\s'\om_{j'}) \le
\de_{j_{\ell}}(\oo\cdot\nn_{\ell}(\underline{t}_{\ell})) \le
\frac{9}{8} \, \de_{j_{\ell}}(\oo\cdot\nn_{\ell}^{0}+\s'\om_{j'}) \le
\frac{18}{15} \, \de_{j_{\ell}}(\oo\cdot\nn_{\ell}) ,
\nonumber \end{equation}
and hence $4\de_{j_{\ell}}(\oo\cdot\nn_{\ell}) \le
5\de_{j_{\ell}}(\oo\cdot\nn_{\ell}(\underline{t}_{\ell})) \le
6\de_{j_{\ell}}(\oo\cdot\nn_{\ell})$.

By the inductive hypothesis and the discussion of the case 1,
in (\ref{eq:6.8}) we have 
\begin{equation}
\left| \oo\cdot\nn_{\ell_{i}}^{0} + \s_{\ell_{i-1}}\om_{j_{i-1}} -
\s_{\ell_{i}}\om_{j_{i}} \right| \le 2^{-n_{i}+2} \g ,
\qquad i=1,\ldots,p ,
\nonumber \end{equation}
where $n_{i}=n_{\ell_{i}}$. Moreover one has
$n_{i} \ge n_{i+1}$ for $i=1,\ldots,p-1$, so that we obtain
\begin{equation}
\de_{j_{\ell}}(\oo\cdot\nn_{\ell}(\underline{t}_{\ell})) \ge
\de_{j_{\ell}}(\oo\cdot\nn_{\ell}^{0}+\s'\om_{j'}) -
\sum_{i=1}^{p} 2^{-n_{i}+2} \g \ge
\de_{j_{\ell}}(\oo\cdot\nn_{\ell}^{0}+\s'\om_{j'}) - 2^{-n+3} \g .
\nonumber \end{equation}
Since $\de_{j_{\ell}}(\oo\cdot\nn_{\ell}^{0}+\s'\om_{j'}) \ge
2^{-(\tilde n_{\ell}+1)}\g$ and $\tilde n_{\ell} \le n -7$, one
finds $\de_{j_{\ell}}(\oo\cdot\nn_{\ell}(\underline{t}_{\ell})) \ge
(1-2^{-3}) \de_{j_{\ell}}(\oo\cdot\nn_{\ell}^{0}+\s'\om_{j'})$.
In the same way one has
$\de_{j_{\ell}}(\oo\cdot\nn_{\ell}(\underline{t}_{\ell})) \le
(1+2^{-3})\de_{j_{\ell}}(\oo\cdot\nn_{\ell}^{0}+\s'\om_{j'})$,
so that (\ref{eq:6.10}) follows.\EP
%%%%%%%%%%%%%%%%%%%%%%%%%%%%%%%%%%%%%%%%%%%%%%%%%%%%%%%%%%%%%%%%%%%%%%%%%%

%%%%%%%%%%%%%%%%%%%%%%%%%%%%%%%%%%%%%%%%%%%%%%%%%%%%%%%%%%%%%%%%%%%%%%%%%%
\begin{rmk} \label{rmk:6.5}
Given a renormalised tree $\theta$, with $\Val(\theta)\neq 0$,
if a line $\ell\in L(\theta)$ has scale $n_{\ell}$ then
$\Psi_{n_{\ell}}(\de_{j_{\ell}}(\oo\cdot\nn_{\ell})(\underline{t}_{\ell}))
\neq0$, and hence, by Lemma \ref{lem:6.4}, one has $2^{-(n_{\ell}+2)}\g
\le \de_{j_{\ell}}(\oo\cdot\nn_{\ell}) \le 2^{-(n_{\ell}-2)}\g$.
Therefore, Lemma \ref{lem:5.1} still holds for the renormalised
trees without any changes in the proof (see also Remark \ref{rmk:5.2}).
\end{rmk}
%%%%%%%%%%%%%%%%%%%%%%%%%%%%%%%%%%%%%%%%%%%%%%%%%%%%%%%%%%%%%%%%%%%%%%%%%%

%%%%%%%%%%%%%%%%%%%%%%%%%%%%%%%%%%%%%%%%%%%%%%%%%%%%%%%%%%%%%%%%%%%%%%%%%%
\begin{rmk} \label{rmk:6.6}
Another important consequence of Lemma \ref{lem:6.4} (and of
inequality (\ref{eq:4.8}) in Remark \ref{rmk:4.11})
is that the number of scale labels which can be
associated with each line of a renormalised tree is still at most $2$.
\end{rmk}
%%%%%%%%%%%%%%%%%%%%%%%%%%%%%%%%%%%%%%%%%%%%%%%%%%%%%%%%%%%%%%%%%%%%%%%%%%

%%%%%%%%%%%%%%%%%%%%%%%%%%%%%%%%%%%%%%%%%%%%%%%%%%%%%%%%%%%%%%%%%%%%%%%%%%
%%%%%%%%%%%%%%%%%%%%%%%%%%%%%%%%%%%%%%%%%%%%%%%%%%%%%%%%%%%%%%%%%%%%%%%%%%
\zerarcounters
\section{Symmetries and identities}
\label{sec:7}
%%%%%%%%%%%%%%%%%%%%%%%%%%%%%%%%%%%%%%%%%%%%%%%%%%%%%%%%%%%%%%%%%%%%%%%%%%
%%%%%%%%%%%%%%%%%%%%%%%%%%%%%%%%%%%%%%%%%%%%%%%%%%%%%%%%%%%%%%%%%%%%%%%%%%

Now we shall prove some symmetry properties on the localized value of the
self-energy clusters.

%%%%%%%%%%%%%%%%%%%%%%%%%%%%%%%%%%%%%%%%%%%%%%%%%%%%%%%%%%%%%%%%%%%%%%%%%%
\begin{lemma} \label{lem:7.1}
If $T\in\gotE^{k}_{j,\s,j,\s}(u,n)$ is such that  $\Tpru$ does not
contain any end node $v$ with $F_{v}=c_{j}^{-\s}$
then there exists $T'\in\overline\gotR^{k}_{j,\s,j,\s}(u,n)$
such that $-2\matL\Val(T)=\matL\Val(T',u)$.
\end{lemma}
%%%%%%%%%%%%%%%%%%%%%%%%%%%%%%%%%%%%%%%%%%%%%%%%%%%%%%%%%%%%%%%%%%%%%%%%%%

\begin{figure}[!ht]
\centering
{
\psfrag{njs}{$\nn\,j\,\s$}
\psfrag{se}{$\s\ee_{j}$}
\psfrag{v0}{$v_{0}$}
\psfrag{T}{$T$}
\psfrag{th}{$\theta_{T}=$}
\psfrag{T'}{$T'$}
\includegraphics[width=15cm]{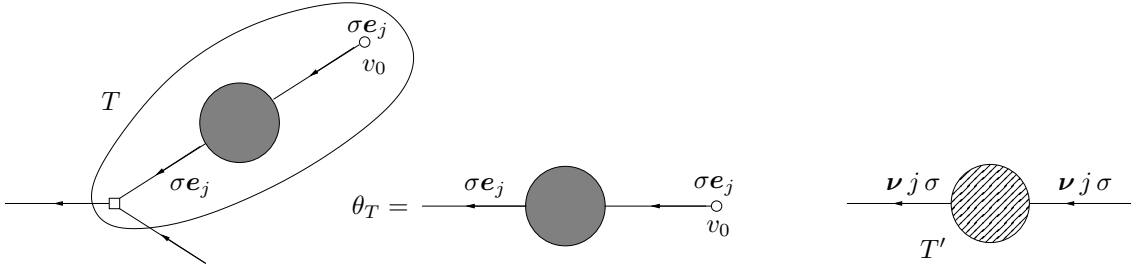}
}
\caption{\label{fig:15}
\footnotesize{The self-energy cluster $T$, the tree $\theta_{T}$,
and the self-energy cluster $T'$ in the proof of Lemma \ref{lem:7.1}.}}
\end{figure}

%%%%%%%%%%%%%%%%%%%%%%%%%%%%%%%%%%%%%%%%%%%%%%%%%%%%%%%%%%%%%%%%%%%%%%%%%%
\prova If $T\in\gotE^{k}_{j,\s,j,\s}(u,n)$ one has $|E_{j}^{\s}
(\Tpru)|=|E_{j}^{-\s}(\Tpru)|+1$ (see Remark \ref{rmk:4.24}),
so that if $|E_{j}^{-\s}(\Tpru)|=0$, then also
$|E_{j}^{\s}(\Tpru)|=1$. This means that $j_{v} \neq j$
for all $v\in E(\Tpru)\setminus\{v_{0}\}$,
if $E_{j}^{\s}(\Tpru)=\{v_{0}\}$. Consider the self-energy cluster
$T'\in\overline\gotR^{k}_{j,\s,j,\s}(u,n)$ obtained from
$\theta_{T}$ by replacing the line exiting $v_{0}$
with an entering line carrying a momentum $\nn$ such that
$\oo\cdot\nn=u$ and $n_{T'}=n_{T}$; see Figure \ref{fig:15}.
With the exception of $v_{0}$, the nodes of $\theta_{T}$
have the same node factors as $T'$; in particular
they have the same combinatorial factors. If we compute the propagators
$G_{\ell}$ of $\ell\in L(T')$, by setting $u=\s\om_{j}$, then they
are the same as the corresponding propagators of $\theta_{T}$.
Finally, as $n_{T'}=n_{T}$, one has
$\matL\Val(T)=0$ if and only if also $\matL\Val(T',u)=0$.
Thus, by recalling also Remark \ref{rmk:4.25},
one finds $-2\matL\Val(T)=\matL\Val(T',u)$.\EP
%%%%%%%%%%%%%%%%%%%%%%%%%%%%%%%%%%%%%%%%%%%%%%%%%%%%%%%%%%%%%%%%%%%%%%%%%%

%%%%%%%%%%%%%%%%%%%%%%%%%%%%%%%%%%%%%%%%%%%%%%%%%%%%%%%%%%%%%%%%%%%%%%%%%%
% figure 15
%%%%%%%%%%%%%%%%%%%%%%%%%%%%%%%%%%%%%%%%%%%%%%%%%%%%%%%%%%%%%%%%%%%%%%%%%%
%%%%%%%%%%%%%%%%%%%%%%%%%%%%%%%%%%%%%%%%%%%%%%%%%%%%%%%%%%%%%%%%%%%%%%%%%%

\begin{figure}[!ht]
\centering
{
\psfrag{n}{$\nn$}
\psfrag{n'}{$\nn'$}
\psfrag{se}{$\s\ee_{j}$}
\psfrag{se'}{$\s'\ee_{j'}$}
\psfrag{se-}{$-\s\ee_{j}$}
\psfrag{se'-}{$-\s'\ee_{j'}$}
\psfrag{v}{$v_{T}$}
\psfrag{T1}{$T_{1}$}
\psfrag{T2}{$T_{2}$}
\psfrag{T3}{$T_{3}$}
\includegraphics[width=15cm]{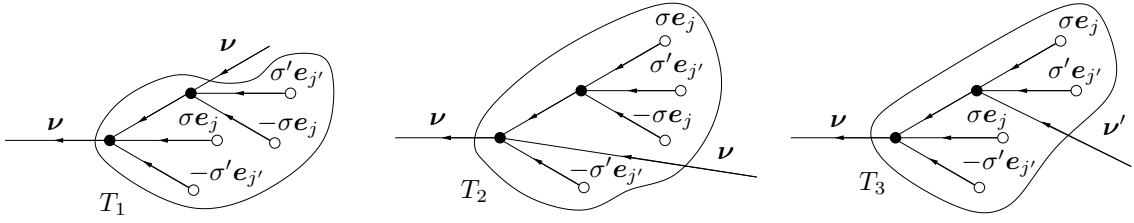}
}
\caption{\label{fig:16}
\footnotesize{The sets $\calF_{1}(T)=\{T_{1},
T_{2}\}$ and $\calF_{2}(T)=\{T_{3}\}$ corresponding
to the self-energy cluster $T$ in Figure \ref{fig:11}.}}
\end{figure}

For $T\in\gotE^{k}_{j,\s,j,\s}(u,n)$ let us call
$\calF_{1}(T)$ the set of all inequivalent self-energy clusters
$T'\in\overline\gotR^{k}_{j,\s,j,\s}(u,n)$ obtained from $\theta_{T}$
by replacing a line exiting an end node $v\in E_{j}^{\s}(\thetapru_{T})$
with an entering line carrying a momentum $\nn$ such that $\oo\cdot\nn=u$
and with $n_{T'}=n_{T}$. Call also $\calF_{2}(T)$ the set of all
inequivalent self-energy clusters $T'\in\gotR^{k}_{j,\s,j,-\s}(u',n)$,
with $u'=u-2\s\om_{j}$, obtained from $\theta_{T}$ by replacing a
line exiting an end node $v\in E_{j}^{-\s}(\thetapru_{T})$ (if any)
with an entering line carrying a momentum $\nn'$ such that
$\oo\cdot\nn'=u'$ and with $n_{T'}=n_{T}$; see Figure \ref{fig:16}.

%%%%%%%%%%%%%%%%%%%%%%%%%%%%%%%%%%%%%%%%%%%%%%%%%%%%%%%%%%%%%%%%%%%%%%%%%%
% figure 16
%%%%%%%%%%%%%%%%%%%%%%%%%%%%%%%%%%%%%%%%%%%%%%%%%%%%%%%%%%%%%%%%%%%%%%%%%%
%%%%%%%%%%%%%%%%%%%%%%%%%%%%%%%%%%%%%%%%%%%%%%%%%%%%%%%%%%%%%%%%%%%%%%%%%%

%%%%%%%%%%%%%%%%%%%%%%%%%%%%%%%%%%%%%%%%%%%%%%%%%%%%%%%%%%%%%%%%%%%%%%%%%%
\begin{lemma} \label{lem:7.2}
For all $T \in \gotE^{k}_{j,\s,j,\s}(u,n)$ one has
\begin{equation}
\Big( 2c_{j}^{\s} \matL\Val(T) +
c_{j}^{\s}\sum_{T'\in \calF_{1}(T)} \matL\Val(T',u) \Big) =
c_{j}^{-\s}\sum_{T'\in \calF_{2}(T)}\matL\Val(T',u'),
\nonumber \end{equation}
where $u'=u-2\s\om_{j}$ and the right hand side is meant as zero
if $\calF_{2}(T) = \emptyset$.
\end{lemma}
%%%%%%%%%%%%%%%%%%%%%%%%%%%%%%%%%%%%%%%%%%%%%%%%%%%%%%%%%%%%%%%%%%%%%%%%%%

%%%%%%%%%%%%%%%%%%%%%%%%%%%%%%%%%%%%%%%%%%%%%%%%%%%%%%%%%%%%%%%%%%%%%%%%%%
\prova The case $k(T)>K_{0}2^{n_{T}/\tau}$ is trivial so that we
consider only the case $k(T)\le K_{0}2^{n_{T}/\tau}$.
By construction any $T\in\gotE^{k}_{j,\s,j,\s}(u,n)$ is
such that $\Tpru$ contains at least an end node $v$
such that $F_{v}=c_{j}^{\s}$, hence $|E^{\s}_{j}(\Tpru)| \ge 1$.
By Lemma \ref{lem:7.1} either $|E^{-\s}_{j}(\Tpru)|\ge 1$ or there
exists $T'\in\overline\gotR^{k}_{j,\s,j,\s}(u,n)$ such that $2\matL\Val(T)
+ \matL\Val(T',u)=0$. Hence the assertion is proved
if $E_{j}^{-\s}(\Tpru)=\emptyset$.

So, let us consider the case $|E_{j}^{-\s}(\Tpru)| \ge 1$.
First of all note that there is a  1-to-1  correspondence between the
lines of $\theta_{T}$ and the lines and external lines,
respectively, of both $T'\in \calF_{1}(T)$ and $T'\in \calF_{2}(T)$;
the same holds for the internal nodes.
Moreover the propagators both of any $T'\in\calF_{1}(T)$ and of any
$T'\in\calF_{2}(T)$ are equal to the corresponding propagators of $T$
when setting $u=\s\om_{j}$ and $u=-\s\om_{j}$, respectively.
Also the node factors of the internal nodes of all
self-energy clusters $T'\in\calF_{1}(T) \cup \calF_{2}(T)$
are the same as those of $T$.
For $T'\in \calF_{1}(T)$ one has $|E_{i}^{+}(\Tpru')|=
|E_{i}^{-}(\Tpru')|$ for all $i=1,\ldots,d$, whereas for
$T''\in \calF_{2}(T)$ one has $|E_{i}^{+}(\Tpru'')|=|E_{i}^{-}(\Tpru'')|$
for all $i \neq j$ and $|E_{j}^{\s}(\Tpru'')|=|E_{j}^{-\s}(\Tpru'')|+2$;
thus, one has
\begin{equation}
\Big( \prod_{v \in E(\Tpru)} c_{j_{v}}^{\s_{v}} \Big) =
c_{j}^{\s} \Big( \prod_{v \in E(\Tpru')} c_{j_{v}}^{\s_{v}} \Big) =
c_{j}^{-\s} \Big( \prod_{v \in E(\Tpru'')} c_{j_{v}}^{\s_{v}} \Big)
\nonumber \end{equation}
for all $T'\in \calF_{1}(T)$ and all $T''\in \calF_{2}(T)$. 

Therefore, if we write
\begin{equation}
-2c_{j}^{\s} \matL\Val(T) =
\Val(\theta_{T}) =\matA(T)
\Big( \prod_{v \in E(\Tpru)} c_{j_{v}}^{\s_{v}} \Big) ,
\label{eq:7.1} \end{equation}
where $\matA(T)$ depends only on $T$, then one finds
\begin{equation}
\sum_{T'\in \calF_{1}(T)} \matL\Val(T',u) = \matA(T)
\frac{1}{c_{j}^{\s}} \Big( \prod_{v \in E(\Tpru)} c_{j_{v}}^{\s_{v}} \Big)
\sum_{v\in V(\Tpru)} r_{v,j,\s} ,
\nonumber \end{equation}
with the same factor $\matA(T)$ as in (\ref{eq:7.1}).
Analogously one has
\begin{equation}
\sum_{T'\in \calF_{2}(T)} \matL\Val(T',u') = \matA(T)
\frac{1}{c_{j}^{\s}} \Big( \prod_{v \in E(\Tpru)} c_{j_{v}}^{-\s_{v}} \Big)
\sum_{v\in V(\Tpru)} r_{v,j,-\s} ,
\nonumber \end{equation}
again with the same factor $\matA(T)$ as in (\ref{eq:7.1}),
so one can write
\begin{equation}
\Big( \!\! - 2c_{j}^{\s} \Val(T) +
c_{j}^{\s} \!\!\!\!\!\! \sum_{T'\in \calF_{1}(T)} \!\!\!
\matL\Val(T',u) \Big) -
c_{j}^{-\s} \!\!\!\!\!\!
\sum_{T'\in \calF_{2}(T)} \!\!\!
\matL\Val(T',u') = \matB(T) \Big( \!\! - 1 + \!\!\!\!
\sum_{v\in V(\Tpru)} \!\! \left( r_{v,j,\s} - r_{v,j,-\s} \right) \Big) ,
\label{eq:7.2} \end{equation}
where
\begin{equation}
\matB(T) = \matA(T)
\Big( \prod_{v \in E(\Tpru)} c_{j_{v}}^{\s_{v}} \Big) .
\nonumber \end{equation}
On the other hand one has
\begin{equation}
\sum_{v\in V(\Tpru)} r_{v,j,\s} = |E_{j}^{\s}(\Tpru)| ,
\nonumber \end{equation}
so that the term in the last parentheses of (\ref{eq:7.2}) gives
$-1+|E_{j}^{\s}(\Tpru)|-|E_{j}^{-\s}(\Tpru)|=0$.
Therefore the assertion is proved.\EP
%%%%%%%%%%%%%%%%%%%%%%%%%%%%%%%%%%%%%%%%%%%%%%%%%%%%%%%%%%%%%%%%%%%%%%%%%%

\begin{figure}[!ht]
\centering
{
\psfrag{n}{$\nn\,j\,\s$}
\psfrag{n'}{$\nn'\,j'\,\s'$}
\psfrag{n''}{$\nn''\,j\,\s$}
\psfrag{n'-}{$\nn'\,j'\,\!-\!\s'$}
\psfrag{l}{$\ell_{T}$}
\psfrag{l'}{$\ell_{T}'$}
\psfrag{n1}{$\nn_{1}$}
\psfrag{n2}{$\nn_{2}$}
\psfrag{n3}{$\nn_{3}$}
\psfrag{n4}{$-\s'\ee_{j'}$}
\psfrag{n1-}{$-\nn_{1}$}
\psfrag{n2-}{$-\nn_{2}$}
\psfrag{n3-}{$-\nn_{3}$}
\psfrag{n4-}{$\s'\ee_{j'}$}
\psfrag{se'}{$\s'\ee_{j'}$}
\psfrag{se-'}{$-\s'\ee_{j'}$}
\psfrag{v0}{$v_{0}$}
\psfrag{T}{$T$}
\psfrag{T1}{$T_{1}$}
\psfrag{T2}{$T_{2}$}
\psfrag{T3}{$T_{3}$}
\psfrag{T4}{$T_{4}$}
\psfrag{T5}{$T_{5}$}
\includegraphics[width=15cm]{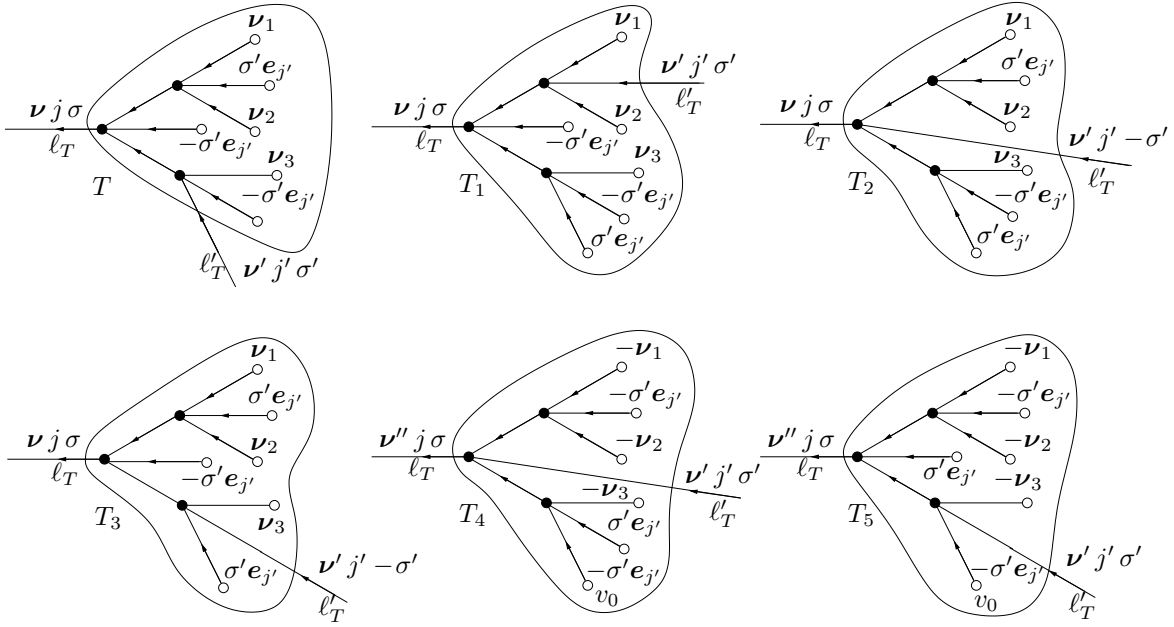}
}
\caption{\label{fig:17}
\footnotesize{A self-energy cluster $T$ and the
corresponding sets $\calG_{1}(T)=\{T,T_{1}\}$,
$\calG_{2}(T)=\{T_{2},T_{3}\}$, and $\calG_{3}(T)=\{T_{4},T_{5}\}$.}}
\end{figure}

For $T\in\gotR^{k}_{j,\s,j',\s'}(u,n)$ with $n_{\ell}\ge 0$ for all
$\ell \in \calP_{T}$, call $\calG_{1}(T)$ the set of
self-energy clusters $T'\in \gotR^{k}_{j,\s,j',\s'}(u,n)$ obtained
from $T$ by exchanging the entering line $\ell_{T}'$
with a line exiting an end node $v\in E_{j'}^{\s'}(\Tpru)$ (if any).
Call also $\calG_{2}(T)$ the set of self-energy clusters
$T'\in \gotR^{k}_{j,\s,j',-\s'}(u',n)$, with $u'=u-2\s\om_{j}$,
obtained from $T$ by (1) replacing
the momentum of $\ell_{T}'$ with a momentum $\nn'$ such that
$\oo\cdot\nn'=u'$, (2) changing the sign label of an end node
$v\in E_{j'}^{-\s'}(\Tpru)$ into $\s'$, and (3) exchanging
the lines $\ell_{T}'$ and $\ell_{v}$.
Finally call $\calG_{3}(T)$ the set of self-energy clusters
$T'\in \gotR^{k}_{j,-\s,j',\s'}(u,n)$, obtained from $T$ by (1) replacing
the entering line $\ell_{T}'$ with a line exiting a new end node
$v_{0}$ with $\s_{v_{0}}=\s'$ and $\nn_{v_{0}}=\s'\ee_{j'}$, (2) replacing
all the labels $\s_{v}$ of the nodes
$v\in N_{0}(T)\cup\{v_{0}\}$ with $-\s_{v}$ and (3)
replacing a line exiting an end node $v\in E^{\s'}_{j'}(\Tpru)$,
with the entering line $\ell_{T}'$; see Figure \ref{fig:17}.
Again we force $n_{T'}=n_{T}$ for all $T'\in\calG_{1}(T)\cup\calG_{2}(T)
\cup \calG_{3}(T)$.

%%%%%%%%%%%%%%%%%%%%%%%%%%%%%%%%%%%%%%%%%%%%%%%%%%%%%%%%%%%%%%%%%%%%%%%%%%
% figure 17
%%%%%%%%%%%%%%%%%%%%%%%%%%%%%%%%%%%%%%%%%%%%%%%%%%%%%%%%%%%%%%%%%%%%%%%%%%
%%%%%%%%%%%%%%%%%%%%%%%%%%%%%%%%%%%%%%%%%%%%%%%%%%%%%%%%%%%%%%%%%%%%%%%%%%

%%%%%%%%%%%%%%%%%%%%%%%%%%%%%%%%%%%%%%%%%%%%%%%%%%%%%%%%%%%%%%%%%%%%%%%%%%
\begin{lemma} \label{lem:7.3}
For all $T \in \gotR^{k}_{j,\s,j',\s'}(u,n)$, with $j\neq j'$ and
$n_{\ell}\ge 0$ for all $\ell \in \calP_{T}$, one has
\begin{equation}
\begin{aligned}
c_{j'}^{\s'} \!\!\! \sum_{T'\in\calG_{1}(T)} \!\!\! \matL\Val(T',u) &=
c_{j'}^{-\s'} \!\!\! \sum_{T'\in\calG_{2}(T)} \!\!\! \matL\Val(T',u')  \\
c_{j}^{-\s}c_{j'}^{\s'} \!\!\! \sum_{T'\in\calG_{1}(T)}\!\!\! \matL\Val(T',u)
&= c_{j}^{\s}c_{j'}^{-\s'} \!\!\! \sum_{T'\in\calG_{3}(T)} \!\!\!
\matL\Val(T',u).
\end{aligned} \nonumber \end{equation}
\end{lemma}
%%%%%%%%%%%%%%%%%%%%%%%%%%%%%%%%%%%%%%%%%%%%%%%%%%%%%%%%%%%%%%%%%%%%%%%%%%

%%%%%%%%%%%%%%%%%%%%%%%%%%%%%%%%%%%%%%%%%%%%%%%%%%%%%%%%%%%%%%%%%%%%%%%%%%
\prova Again we consider only the case $k(\To)\le K_{0}2^{n_{T}/\tau}$.
For fixed $T \in \gotR^{k}_{j,\s,j',\s'}(u,n)$, with $j\neq j'$,
let $\theta\in\gotT_{j,\s\ee_{j}}^{k}(n)$ be the tree obtained
from $T$ by replacing the entering line $\ell_{T}'$  with a line
exiting a new end node $v_{0}$ with $\s_{v_{0}}=\s'$ and $\nn_{v_{0}}=
\s'\ee_{j'}$. 
Note that in particular one has $|E_{j'}^{\s}(\thetapru)|=
|E_{j'}^{-\s}(\thetapru)|$. Any $T'\in\calG_{1}(T)$ can be obtained from
$\theta$ by replacing a line exiting an end node $v\in E_{j'}^{\s'}
(\thetapru)$ with entering line $\ell_{T'}'$,
with the same labels as $\ell_{T}$, so that
\begin{equation}
c_{j'}^{\s'} \!\!\! \sum_{T'\in\calG_{1}(T)} \!\!\! \matL\Val(T',u)=
|E_{j'}^{\s'}(\thetapru)|\Val(\theta).
\nonumber \end{equation}
On the other hand, any $T'\in\calG_{2}(T)$ can be obtained from $\theta$ 
by replacing a line exiting an end node $v \in E_{j'}^{-\s'}(\thetapru)$
with entering line $\ell_{T'}'$, with labels
$\nn'-2\s'\ee_{j'}$, $j'$, $-\s'$, hence
\begin{equation}
c_{j'}^{-\s'} \!\!\! \sum_{T'\in\calG_{2}(T)} \!\!\! \matL\Val(T',u)=
|E_{j'}^{-\s'}(\thetapru)|\Val(\theta),
\nonumber \end{equation}
so that the first equality is proved.

Now, let $\theta'\in\gotT_{j,-\s\ee_{j}}^{k}(n)$ be the tree obtained
from $\theta$ by replacing all the labels $\s_{v}$ of
the nodes $v\in N_{0}(\theta)$ with $-\s_{v}$.
Any $T'\in\calG_{3}(T)$ can be obtained from $\theta'$ by
replacing a line exiting an end node $v \in E_{j'}^{\s'}(\thetapru')$
with entering line $\ell_{T'}'$, carrying the same labels as $\ell_{T}$.
Hence, by Lemma \ref{lem:4.14},
\begin{equation}
c_{j}^{-\s}c_{j'}^{\s'} \!\!\! \sum_{T'\in\calG_{1}(T)} \!\!\!
\matL\Val(T',u)= c_{j}^{-\s}|E_{j'}^{\s'}(\thetapru)|
\Val(\theta)=c_{j}^{\s} |E_{j'}^{-\s'}(\thetapru')|
\Val(\theta')= c_{j}^{\s}c_{j'}^{-\s'} \!\!\!
\sum_{T'\in\calG_{3}(T)} \!\!\! \matL\Val(T',u) ,
\nonumber \end{equation}
which yields the second identity, and hence
completes the proof.\EP
%%%%%%%%%%%%%%%%%%%%%%%%%%%%%%%%%%%%%%%%%%%%%%%%%%%%%%%%%%%%%%%%%%%%%%%%%%

%%%%%%%%%%%%%%%%%%%%%%%%%%%%%%%%%%%%%%%%%%%%%%%%%%%%%%%%%%%%%%%%%%%%%%%%%%
\begin{lemma} \label{lem:7.4}
For all $k\in\ZZZ_{+}$, all $j,j'=1,\ldots,d$,
and all $\s,\s'\in\{\pm\}$, one has\\
(i) $\hhh^{(k)}=\hhh^{(k)}(|c_{1}|^{2},\ldots,|c_{d}|^{2})$, i.e.,
$\hhh^{(k)}$ depends on $\cc$ only through the quantities
$|c_{1}|^{2},\ldots,|c_{d}|^{2}$;\\
(ii) $\matL M^{(k)}_{j,\s,j',\s'}(u,n)=c_{j}^{-\s}c_{j'}^{\s'} 
M^{(k)}_{j,j'}(n)$, where $M^{(k)}_{j,j'}(n)$ does not depend on the
indices $\s,\s'$.
\end{lemma}
%%%%%%%%%%%%%%%%%%%%%%%%%%%%%%%%%%%%%%%%%%%%%%%%%%%%%%%%%%%%%%%%%%%%%%%%%%

%%%%%%%%%%%%%%%%%%%%%%%%%%%%%%%%%%%%%%%%%%%%%%%%%%%%%%%%%%%%%%%%%%%%%%%%%%
\prova One works on the single trees contributing
to $\matL M^{(k)}_{j,\s,j',\s'}(u,n)$.
Then the proof follows from Lemma \ref{lem:4.14} and the results above.\EP
%%%%%%%%%%%%%%%%%%%%%%%%%%%%%%%%%%%%%%%%%%%%%%%%%%%%%%%%%%%%%%%%%%%%%%%%%%

%%%%%%%%%%%%%%%%%%%%%%%%%%%%%%%%%%%%%%%%%%%%%%%%%%%%%%%%%%%%%%%%%%%%%%%%%%

%%%%%%%%%%%%%%%%%%%%%%%%%%%%%%%%%%%%%%%%%%%%%%%%%%%%%%%%%%%%%%%%%%%%%%%%%%
\begin{rmk} \label{rmk:7.5}
Note that Lemma \ref{lem:7.4} could be reformulated as
$$ \matL M^{(k)}_{j,\s,j',\s'}(u,n)=\partial_{c_{j'}^{\s'}}c_{j}^{\s}
\matL\widetilde{M}_{j,\s,j,\s}^{(k)}(n) , $$
with $\widetilde{M}_{j,\s,j,\s}^{(k)}(n)$ defined after (\ref{eq:4.13}).
We omit the proof of the identity, since it will not be used.
\end{rmk}
%%%%%%%%%%%%%%%%%%%%%%%%%%%%%%%%%%%%%%%%%%%%%%%%%%%%%%%%%%%%%%%%%%%%%%%%%%

%%%%%%%%%%%%%%%%%%%%%%%%%%%%%%%%%%%%%%%%%%%%%%%%%%%%%%%%%%%%%%%%%%%%%%%%%%
%%%%%%%%%%%%%%%%%%%%%%%%%%%%%%%%%%%%%%%%%%%%%%%%%%%%%%%%%%%%%%%%%%%%%%%%%%
\zerarcounters
\section{Cancellations and bounds}
\label{sec:8}
%%%%%%%%%%%%%%%%%%%%%%%%%%%%%%%%%%%%%%%%%%%%%%%%%%%%%%%%%%%%%%%%%%%%%%%%%%
%%%%%%%%%%%%%%%%%%%%%%%%%%%%%%%%%%%%%%%%%%%%%%%%%%%%%%%%%%%%%%%%%%%%%%%%%%

We have seen in Section \ref{sec:5} that, as far as resonant lines
are not considered, no problems arise in obtaining `good bounds',
i.e., bounds on the tree values of order $k$ proportional to some constant
to the power $k$ (see Lemma \ref{lem:5.3}). For the same bound to hold
for all tree values we need a gain  factor  proportional to $2^{-n_{\ell}}$
for each resonant line $\ell$ on scale $n_{\ell}\ge 1$.

Let us consider a tree $\theta$, and write its value as in (\ref{eq:6.4}).
Let $\ell$ be a resonant line. Then $\ell$ exits a self-energy cluster
$T_{2}$ and enters a self-energy cluster $T_{1}$; see Figure \ref{fig:9}.
By construction $T_{1}\in\gotR^{k_{1}}_{j_{1},\s_{1},j_{1}',\s_{1}'}
(\oo\cdot\nn_{\ell_{T_{1}}'},n_{1})$ and
$T_{2}\in\gotR^{k_{2}}_{j_{2},\s_{2},j_{2}',\s_{2}'}
(\oo\cdot\nn_{\ell_{T_{2}}'},n_{2})$, for suitable values
of the labels, with the constraint $j_{1}=j_{2}'=j_{\ell}$ and
$\s_{1}=\s_{2}'=\s_{\ell}$.

\begin{figure}[!ht]
\begin{center}{
\psfrag{l}{$\ell$}
\psfrag{js1}{$\nn_{1}\,j_{1}\,\s_{1}$}
\psfrag{js2'}{$\nn_{2}'\,j_{2}'\,\s_{2}'$}
\psfrag{jsl}{$\nn_{\ell}\,j_{\ell}\,\s_{\ell}$}
\psfrag{jsl'}{$\nn_{\ell}' \, j_{\ell}\, -\s_{\ell}$}
\psfrag{sum1}{$\displaystyle\sum_{\substack{T_{1}'\in
\gotR^{k_{1}}_{j_{1},\s_{1},j_{\ell},\s_{\ell}}(\oo\cdot\nn_{\ell},n_{1}),\\
T_{2}'\in\gotR^{k_{2}}_{j_{\ell},\s_{\ell},j_{2}',\s_{2}'}
(\oo\cdot\nn_{2}',n_{2})}}$}
\psfrag{sum2}{$\displaystyle\sum_{\substack{T_{1}'\in
\gotR^{k_{1}}_{j_{1},\s_{1},j_{\ell},-\s_{\ell}}(\oo\cdot\nn_{\ell}',n_{1}),\\
T_{2}'\in\gotR^{k_{2}}_{j_{\ell},-\s_{\ell},j_{2}',\s_{2}'}
(\oo\cdot\nn_{2}',n_{2})}}$}
\psfrag{T1}{$T_{1}'$}
\psfrag{T2}{$T_{2}'$}
\psfrag{L}{$\matL$}
\includegraphics[width=14cm]{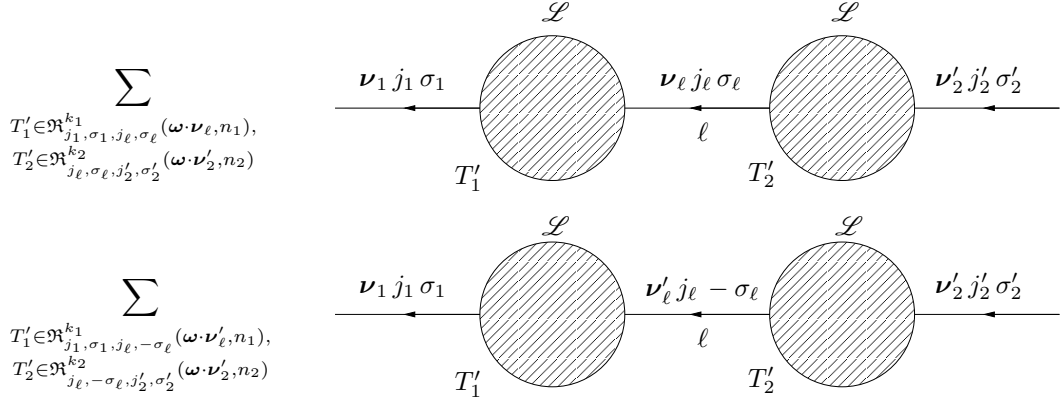}
}
\caption{\label{fig:18}
\footnotesize{Graphical representation of the cancellation mechanism
discussed in the text: $\nn_{\ell}'=\nn_{\ell}-2\s_{\ell}\ee_{j_{\ell}}$.
If we sum the two contributions we obtain a gain factor
$O(2^{-n_{\ell}})$.}}
\end{center}
\end{figure}

If $\matO_{T_{1}}=\matO_{T_{2}}=\matL$,
we consider also all trees obtained from $\theta$
by replacing $T_{1}$ and $T_{2}$ with
other clusters $T_{1}'\in\gotR^{k_{1}}_{j_{1},\s_{1},j_{1}',\s_{1}'}
(\oo\cdot\nn_{\ell_{T_{1}}'},n_{1})$ and
$T_{2}'\in\gotR^{k_{2}}_{j_{2},\s_{2},j_{2}',\s_{2}'}
(\oo\cdot\nn_{\ell_{T_{2}}'},n_{2})$, respectively,
with $\matO_{T_{1}'}=\matO_{T_{2}'}=\matL$. In this way
\begin{equation}
\matL\Val(T_{1},\oo\cdot\nn_{\ell_{T_{1}}'}) \,
G^{[n_{\ell}]}_{j_{\ell}}(\oo\cdot\nn_{\ell}) \,
\matL\Val(T_{2},\oo\cdot\nn_{\ell_{T_{2}}}) 
\nonumber \end{equation}
is replaced with
\begin{equation}
\matL M^{(k_{1})}_{j_{1},\s_{1},j_{\ell},\s_{\ell}}
(\oo\cdot\nn_{\ell_{T_{1}}'},n_{1}) \,
G^{[n_{\ell}]}_{j_{\ell}}(\oo\cdot\nn_{\ell}) \, \matL
M^{(k_{2})}_{j_{\ell},\s_{\ell},j_{2}',\s_{2}'}
(\oo\cdot\nn_{\ell_{T_{2}}'},n_{2}) .
\label{eq:8.1} \end{equation}
Then consider also all trees in which the factor (\ref{eq:8.1})
is replaced with
\begin{equation}
\matL M^{(k_{1})}_{j_{1},\s_{1},j_{\ell},-\s_{\ell}}
(\oo\cdot\nn_{\ell_{T_{1}}'},n_{1}) \,
G^{[n_{\ell}]}_{j_{\ell}}(\oo\cdot\nn_{\ell}') \, \matL
M^{(k_{2})}_{j_{\ell},-\s_{\ell},j_{2}',\s_{2}}
(\oo\cdot\nn_{\ell_{T_{2}}'},n_{2}) ,
\label{eq:8.2} \end{equation}
with $\nn_{\ell}'$ such that $\oo\cdot\nn_{\ell}-\s_{\ell}\om_{j_{\ell}}=
\oo\cdot\nn_{\ell}'+\s_{\ell}\om_{j_{\ell}}$; see Figure \ref{fig:18}.
Because of Lemmas \ref{lem:7.2} and \ref{lem:7.3} the
sum of the two contributions (\ref{eq:8.1}) and (\ref{eq:8.2}) gives
\begin{equation}
\matL M^{(k_{1})}_{j_{1},\s_{1},j_{\ell},\s_{\ell}}
(\oo\cdot\nn_{\ell_{T_{1}}'},n_{1}) \left(
G^{[n_{\ell}]}_{j_{\ell}}(\oo\cdot\nn_{\ell}) +
G^{[n_{\ell}]}_{j_{\ell}}(\oo\cdot\nn_{\ell}') \right) \matL
M^{(k_{2})}_{j_{\ell},\s_{\ell},j_{2}',\s_{2}'}
(\oo\cdot\nn_{\ell_{T_{2}}'},n_{2}) ,
\nonumber \end{equation}
where
\begin{eqnarray}
G^{[n_{\ell}]}_{j_{\ell}}(\oo\cdot\nn_{\ell}) +
G^{[n_{\ell}]}_{j_{\ell}}(\oo\cdot\nn_{\ell}')
& \!\!\! = \!\!\! &
\frac{\Psi_{n_{\ell}}(\de_{j_{\ell}}(\oo\cdot\nn_{\ell}))}{
(\oo\cdot\nn_{\ell}-\s_{\ell}\om_{j_{\ell}})} \left(
\frac{1}{\oo\cdot\nn_{\ell}+\s_{\ell}\om_{j_{\ell}}} +
\frac{1}{\oo\cdot\nn_{\ell}'-\s_{\ell}\om_{j_{\ell}}} \right)
\nonumber \\
& \!\!\! = \!\!\! &
\frac{2\Psi_{n_{\ell}}(\de_{j_{\ell}}(\oo\cdot\nn_{\ell}))}{
(\oo\cdot\nn_{\ell}+\s_{\ell}\om_{j_{\ell}})
(\oo\cdot\nn_{\ell}'-\s_{\ell}\om_{j_{\ell}})} ,
\label{eq:cancel} \end{eqnarray}
and hence $| G^{[n_{\ell}]}_{j_{\ell}}(\oo\cdot\nn_{\ell}) +
G^{[n_{\ell}]}_{j_{\ell}}(\oo\cdot\nn_{\ell}')|\le 2\om_{j_{\ell}}^{-2}$.
This provides the gain factor $O(2^{-n_{\ell}})$ we were
looking for, with respect to the original bound $C_{0}2^{n_{\ell}}$
on the propagator $G_{\ell}$.

%%%%%%%%%%%%%%%%%%%%%%%%%%%%%%%%%%%%%%%%%%%%%%%%%%%%%%%%%%%%%%%%%%%%%%%%%%
% figure 18
%%%%%%%%%%%%%%%%%%%%%%%%%%%%%%%%%%%%%%%%%%%%%%%%%%%%%%%%%%%%%%%%%%%%%%%%%%
%%%%%%%%%%%%%%%%%%%%%%%%%%%%%%%%%%%%%%%%%%%%%%%%%%%%%%%%%%%%%%%%%%%%%%%%%%

If $\matO_{T_{1}}=\matR$ then if $k(\To_{1}) > K_{0}2^{n_{T_{1}}/\tau}$ we can
extract a factor $C^{k(\To_{1})}$ from $\Val(T_{1},\oo\cdot\nn_{\ell_{T_{1}}'})$
($C$ is the constant appearing in Lemma \ref{lem:5.3}),
and, after writing
$C^{k(\To_{1})}=C^{2 k(\To_{1})} C^{-k(\To_{1})}$, use that
$C^{-k(\To_{1})}\le C^{-K_{0}2^{n_{T_{1}}/\tau}} \le {\rm const.}2^{-n_{T_{1}}}$
in order  to obtain a gain factor $O(2^{-n_{\ell}})$.

If $k(\To_{1}) \le K_{0}2^{n_{T_{1}}/\tau}$ and 
$n_{\ell}\ge 0$ for all $\ell \in \calP_{T}$,
we obtain a gain factor proportional to $2^{-n_{\ell}}$ because of
the first line of (\ref{eq:6.2}). 
Of course whenever one has such a case,
then one has a derivative acting on $\Val(T,u)$ -- see (\ref{eq:6.2}).
Therefore one needs to control derivatives like
\begin{equation}
\partial_{u} \Val(T,u) =
\sum_{\ell\in \calP_{T}} \partial_{u} G_{\ell}
\Big( \prod_{\ell'\in L(T)\setminus \{\ell\} } G_{\ell'} \Big)
\Big( \prod_{v\in N(T)} F_{v} \Big) ,
\label{eq:8.3} \end{equation}
where
\begin{equation}
\partial_{u} G_{\ell} =
\frac{\partial_{u} \Psi_{n_{\ell}}(\de(\oo\cdot\nn_{\ell}))}
{(\oo\cdot\nn_{\ell})^{2}-\om_{j_{\ell}}^{2}} -
2 \oo\cdot\nn_{\ell} \frac{\Psi_{n_{\ell}}(\de(\oo\cdot\nn_{\ell}))}{
((\oo\cdot\nn_{\ell})^{2}-\om_{j_{\ell}}^{2})^{2}} .
\label{eq:8.4}\end{equation}
The derived propagator (\ref{eq:8.4}) can be easily bounded by
\begin{equation}
\left| \partial_{u} G_{\ell} \right| \le
C_{1} 2^{2n_{\ell}} ,
\label{eq:8.5}\end{equation}
for some positive constant $C_{1}$.

In principle, given a line $\ell$, one could have one derivative
of $G_{\ell}$ for each self-energy cluster containing $\ell$.
This should be a problem, because in a tree of order $k$,
a propagator $G_{\ell}$ could be derived up to $O(k)$ times,
and no bound proportional to some constant to the power $k$ can
be expected to hold to order $k$. In fact, it happens that
no propagator has to be derived more than once.
This can be seen by reasoning as follows.

Let $T$ be a self-energy cluster of depth $D(T)=1$. If $\matO_{T}=\matR$
then a gain factor $O(2^{-n_{\ell_{T}}})$ is obtained. When writing
$\partial_{u}\Val(T,u)$ according to (\ref{eq:8.3}) one obtains
$|\calP_{T}|$ terms, one for each line $\ell\in\calP_{T}$.
Then we can bound the derivative of $G_{\ell}$ according
to (\ref{eq:8.5}). By collecting together the gain factor
and the bound (\ref{eq:8.5}) we obtain $2^{2n_{\ell}}2^{-n_{\ell_{T}'}}$.
We can interpret such a bound by saying that, at the cost
of replacing the bound $2^{n_{\ell}}$ of the propagator $G_{\ell}$
with its square $2^{2n_{\ell}}$, we have a gain factor
$2^{-n_{\ell_{T}'}}$ for the self-energy cluster $T$.

Suppose that $\ell$ is contained inside other self-energy cluster
besides $T$, say $T_{p} \subset T_{p-1} \subset \ldots \subset T_{1}$
(hence $T_{p}$ is that with largest depth, and $D(T_{p})=p+1$).
Then, when taking the contribution to (\ref{eq:8.3}) with
the derivative $\partial_{u}$ acting on the propagator $G_{\ell}$,
we consider together the labels $\matO_{T_{i}}=\matR$ and
$\matO_{T_{i}}=\matL$ for all $i=1,\ldots,p$ (in other words
we do not distinguish between localised and regularised values
for such self-energy clusters), because we do not want to produce
further derivatives on the propagator $G_{\ell}$.
Of course we have obtained no gain factor corresponding to
the entering lines of the self-energy clusters $T_{1},\ldots,T_{p}$,
and all these lines can be resonant lines. So, eventually
we shall have to keep track of this.

Then we can iterate the procedure. If the self-energy cluster
$T$ does not contain any line whose propagator is derived,
we split its value into the sum of the localised value plus
the regularised value. On the contrary, if a line
along the path $\calP_{T}$ of $T$ is derived we do not
separate the localised value of $T$ from its regularised value.
Note that, if $T$ is contained inside a regularised self-energy cluster,
then both $\oo\cdot\nn_{\ell}$ and $\oo\cdot\nn_{\ell}'$
in (\ref{eq:8.1}) and (\ref{eq:8.2}) must be replaced with
$\oo\cdot\nn_{\ell}(\underline{t}_{\ell})$ and
$\oo\cdot\nn_{\ell}'(\underline{t}_{\ell})$, respectively,
but still $\oo\cdot\nn_{\ell}(\underline{t}_{\ell})-\s_{\ell}\om_{j_{\ell}}
=\oo\cdot\nn_{\ell}'(\underline{t}_{\ell})+\s_{\ell}\om_{j_{\ell}}$,
so that the cancellation (\ref{eq:cancel}) still holds.

Let us call \emph{ghost line} a resonant line $\ell$
such that (1) $\ell$ is along the path $\calP_{T}$ of
a regularised self-energy cluster $T$ and either
(2a) $\ell$ enters or exits a self-energy cluster $T'\subset T$
containing a line whose propagator is derived or
(2b) the propagator of $\ell$ is derived.
Then, eventually one obtains a gain
$2^{-n_{\ell}}$ for all resonant lines $\ell$, except for the ghost lines.
In other words we can say that there is an overall factor
proportional to
\begin{equation}
\Big( \prod_{\ell\in L_{\rm R}(\theta)} 2^{-n_{\ell}} \Big)
\Big( \prod_{\ell\in L_{\rm G}(\theta)} 2^{n_{\ell}}\Big) ,
\label{eq:8.6} \end{equation}
where $L_{\rm G}(\theta)$ is the set of ghost lines.
Indeed, in case (2a) there is no gain corresponding to the line
$\ell$, so that we can insert a `good' factor
$2^{-n_{\ell}}$ provided we allow also a compensating
`bad' factor $2^{n_{\ell}}$. In case (2b) one can reason as follows.
Call (with some abuse of notation) $T_{1}$ and $T_{2}$ the self-energy
clusters which $\ell$ enters and exits, respectively. If
$\matO_{T_{1}}=\matO_{T_{2}}=\matL$, we consider
\begin{equation}
\matL\Val(T_{1},\oo\cdot\nn_{\ell_{T_{1}}'}) \,
\partial_{u} G^{[n_{\ell}]}_{j_{\ell}}
(\oo\cdot\nn_{\ell}(\underline{t}_{\ell})) \,
\matL\Val(T_{2},\oo\cdot\nn_{\ell_{T_{2}}}) ,
\nonumber \end{equation}
and, by summing over all possible self-energy clusters
as done in (\ref{eq:8.1}), we obtain
\begin{equation}
\matL M^{(k_{1})}_{j_{1},\s_{1},j_{\ell},\s_{\ell}}
(\oo\cdot\nn_{\ell_{T_{1}}'},n_{1}) \,
\partial_{u} G^{[n_{\ell}]}_{j_{\ell}}
(\oo\cdot\nn_{\ell}(\underline{t}_{\ell})) \, \matL
M^{(k_{2})}_{j_{\ell},\s_{\ell},j_{2}',\s_{2}'}
(\oo\cdot\nn_{\ell_{T_{2}}'},n_{2}) ;
\nonumber \end{equation}
then we sum this contribution with
\begin{equation}
\matL M^{(k_{1})}_{j_{1},\s_{1},j_{\ell},-\s_{\ell}}
(\oo\cdot\nn_{\ell_{T_{1}}'},n_{1}) \,
\partial_{u} G^{[n_{\ell}]}_{j_{\ell}}
(\oo\cdot\nn_{\ell}'(\underline{t}_{\ell})) \, \matL
M^{(k_{2})}_{j_{\ell},-\s_{\ell},j_{2}',\s_{2}}
(\oo\cdot\nn_{\ell_{T_{2}}'},n_{2}) ,
\nonumber\end{equation}
where $\nn_{\ell}'=\nn_{\ell}-2\s_{\ell}\ee_{j_{\ell}}$; again we can use
Lemmas \ref{lem:7.2} and \ref{lem:7.3} to obtain
\begin{equation}
\matL M^{(k_{1})}_{j_{1},\s_{1},j_{\ell},\s_{\ell}}
(\oo\cdot\nn_{\ell_{T_{1}}'},n_{1}) \,
\big(\partial_{u} G^{[n_{\ell}]}_{j_{\ell}}
(\oo\cdot\nn_{\ell}(\underline{t}_{\ell})) +
\partial_{u} G^{[n_{\ell}]}_{j_{\ell}}
(\oo\cdot\nn_{\ell}'(\underline{t}_{\ell}))\big) \,
\matL M^{(k_{2})}_{j_{\ell},\s_{\ell},j_{2}',\s_{2}'}
(\oo\cdot\nn_{\ell_{T_{2}}'},n_{2}) ,
\nonumber \end{equation}
where
\begin{eqnarray}
\partial_{u} G^{[n_{\ell}]}_{j_{\ell}}
(\oo\cdot\nn_{\ell}(\underline{t}_{\ell})) +
\partial_{u} G^{[n_{\ell}]}_{j_{\ell}}
(\oo\cdot\nn_{\ell}'(\underline{t}_{\ell}))
& \!\!\! = \!\!\! &
\frac{2 \partial_{u} \Psi_{n_{\ell}}
(\delta(\oo\cdot\nn_{\ell}(\underline{t}_{\ell})))}
{(\oo\cdot\nn_{\ell}(\underline{t}_{\ell})+\s_{\ell}\om_{j_{\ell}})
(\oo\cdot\nn_{\ell}'(\underline{t}_{\ell})-\s_{\ell}\om_{j_{\ell}})} 
\nonumber \\
& \!\!\! - \!\!\! &
\frac{4(\oo\cdot\nn_{\ell}(\underline{t}_{\ell})-\s_{\ell}\om_{j_{\ell}})
\Psi_{n_{\ell}}(\delta(\oo\cdot\nn_{\ell}(\underline{t}_{\ell})))}
{(\oo\cdot\nn_{\ell}(\underline{t}_{\ell})+\s_{\ell}\om_{j_{\ell}})^{2}
(\oo\cdot\nn_{\ell}'(\underline{t}_{\ell})-\s_{\ell}\om_{j_{\ell}})^{2}} ,
\nonumber \end{eqnarray}
so that we have not only the gain factor $2^{-n_{\ell}}$ due to the
cancellation, but also a factor $2^{n_{\ell}}$ because of the term
$\partial_{u} \Psi_{n_{\ell}}(\delta(\oo\cdot\nn_{\ell}))$.

A trivial but important remark is that all the ghost lines
contained inside the same self-energy cluster have
different scales: in particular there is at most one ghost line
on a given scale $n$. Therefore we can rely upon Lemma \ref{lem:5.4}
and Lemma \ref{lem:6.4}, to ensure that
for each of such lines there is also at least one non-resonant line
on scale $\ge n-3$ (inside the same self-energy cluster).
Therefore we can bound the second product in (\ref{eq:8.6}) as
\begin{equation}
\Big( \prod_{\ell\in L_{\rm G}(\theta)} 2^{n_{\ell}} \Big) \le
\prod_{n=1}^{\io} 2^{n\gotN_{n-3}(\theta)} ,
\nonumber \end{equation}
which in turn is bounded as a constant to the power $k=k(\theta)$,
as argued in the proof of Lemma \ref{lem:5.3}.

Finally if $k(\To_{1})\le K_{0}2^{n_{T_{1}}/\tau}$ and $T_{1}$ contains
at least one line $\ell\in \calP_{T_{1}}$ with $n_{\ell'}=-1$,
in general there are $p\ge 1$
self-energy clusters $T_{p}' \subset T_{p-1}' \subset \ldots \subset 
T_{1}'=T_{1}$
such that $\ell\in\calP_{T_{i}'}$ for $i=1,\ldots,p$, and $T_{p}'$ is
the one with largest depth containing $\ell$. For $i=1,\ldots,p$
call $\ell_{i}$ the exiting line of the self-energy cluster $T_{i}'$
and $\theta_{i}=\theta_{\ell_{i}}$. Denote also,
for $i=1,\ldots,p-1$, by $\Ga_{i}=\Ga_{\ell_{i+1}}(\theta_{i})$ 
(recall Notation \ref{def:4.29}). By Lemma \ref{lem:4.30} one has
$|\nn_{\ell_{i}}-\nn_{\ell_{i+1}}| \le E_{2} k(\Gao_{i})$ for
$i=1,\ldots,p-1$. Moreover one has $|\nn_{\ell_{1}}-\s_{\ell}\ee_{j_{\ell}}|
\le E_{2} (k(\Gao_{1})+\ldots+k(\Gao_{p-1}))$. On the other hand
one has
\begin{eqnarray}
& \displaystyle{
\frac{\g}{|\nn_{\ell_{i}}-\nn_{\ell_{i+1}}|^{\tau}} \le
\de_{j_{i}}(\oo\cdot\nn_{\ell_{i}})+
\de_{j_{i+1}}(\oo\cdot\nn_{\ell_{i+1}}) \le 2^{-n_{T_{i+1}'}+2} \g } , 
\nonumber \\
& \displaystyle{
\frac{\g}{|\nn_{\ell_{1}}-\s_{\ell}\ee_{j_{\ell}}|^{\tau}} \le
\de_{j_{1}}(\oo\cdot\nn_{\ell_{1}}) \le 2^{-n_{T_{1}'}} \g } ,
\nonumber \end{eqnarray}
so that one can write
\begin{equation}
C^{k(\Gao_{1})+\ldots+k(\Gao_{p-1}))} \le
C^{3k(\Gao_{1})+\ldots+k(\Gao_{p-1}))}
2^{-n_{T_{1}'}} \prod_{i=2}^{p} 2^{-n_{T_{i}'}} ,
\end{equation}
which assures the gain factors for all self-energy clusters
$T_{1}',\ldots,T_{p}'$.

To conclude the analysis, if $\matO_{T_{1}}=\matL$ but
$\matO_{T_{2}}=\matR$, one can reason
in the same way by noting that $|n_{\ell_{T_{2}}'}-n_{\ell}|\le1$.

%%%%%%%%%%%%%%%%%%%%%%%%%%%%%%%%%%%%%%%%%%%%%%%%%%%%%%%%%%%%%%%%%%%%%%%%%%
\begin{lemma} \label{lem:8.1}
Set $\Gamma(\cc)=\max\{|c_{1}|,\ldots,|c_{d}|,1\}$.
There exists a positive constant $C$ such that for $k\in\NNN$,
$j\in\{1,\ldots,d\}$ and $\nn\in\ZZZ^{d}$ one has
$\displaystyle | \sum_{\theta\in\gotT^{k}_{j,\nn}}\Val(\theta) | \le C^{k}
\Gamma^{3k}(\cc)$.
\end{lemma}
%%%%%%%%%%%%%%%%%%%%%%%%%%%%%%%%%%%%%%%%%%%%%%%%%%%%%%%%%%%%%%%%%%%%%%%%%%

%%%%%%%%%%%%%%%%%%%%%%%%%%%%%%%%%%%%%%%%%%%%%%%%%%%%%%%%%%%%%%%%%%%%%%%%%%
\prova Each time one has a resonant line $\ell$,
when summing together the values of all self-energy clusters, a gain
$B_{1}2^{-n_{\ell}}$ is obtained (either by the cancellation mechanism
described at the beginning of this Section or as an effect of
the regularisation operator $\matR$). The number of trees of
order $k$ is bounded by $B_{2}^{k}$ for some constant $B_{2}$;
see Remark \ref{rmk:4.12}. The derived propagators can be bounded
by (\ref{eq:8.5}). By taking into account also the bound of
Lemma \ref{lem:5.3}, setting $B_{3}=C_{0}\Phi$, and
bounding by $B^{k}_{4}$, with
\begin{equation}
B_{4} = \exp \left( 3c \log 2 \sum_{n=0}^{\io} 2^{-n/\tau} n \right) ,
\nonumber \end{equation}
the product of the propagators (both derived and non-derived) of the
non-resonant lines times the derived propagators of the resonant lines,
we obtain the assertion with $C=B_{1}B_{2}B_{3}B_{4}$.\EP
%%%%%%%%%%%%%%%%%%%%%%%%%%%%%%%%%%%%%%%%%%%%%%%%%%%%%%%%%%%%%%%%%%%%%%%%%%

%%%%%%%%%%%%%%%%%%%%%%%%%%%%%%%%%%%%%%%%%%%%%%%%%%%%%%%%%%%%%%%%%%%%%%%%%%
\begin{lemma} \label{lem:8.2}
The function (\ref{eq:2.4}), with $x_{j,\nn}$ as in (\ref{eq:4.10}),
and the counterterms $\h_{j}$ defined in (\ref{eq:4.11}) are analytic
in $\e$ and $\cc$, for $|\e| \Gamma^{3}(\cc) \le \eta_{0}$
with $\eta_{0}$ small enough and $\Gamma(\cc)=\max\{|c_{1}|,\ldots,
|c_{d}|,1\}$. Therefore the solution $\xx(t,\e,\cc)$
is analytic in $t,\e,\cc$ for $|\e| \Gamma^{3}(\cc)
{\rm e}^{3|\oo|\,|{\rm Im}\, t|} \le \eta_{0}$,
with $\eta_{0}$ small enough.
\end{lemma}
%%%%%%%%%%%%%%%%%%%%%%%%%%%%%%%%%%%%%%%%%%%%%%%%%%%%%%%%%%%%%%%%%%%%%%%%%%

%%%%%%%%%%%%%%%%%%%%%%%%%%%%%%%%%%%%%%%%%%%%%%%%%%%%%%%%%%%%%%%%%%%%%%%%%%
\prova Just collect together all the results above, in order to obtain
the convergence of the series for $\eta_{0}$ small enough
and $|\e| \Gamma^{\xi}(\cc) \le \eta_{0}$, for some
constant $\xi$. Moreover $x^{(k)}_{j,\nn}=0$ for $|\nn| > \xi k$,
for the same constant $\xi$. Lemma \ref{lem:4.10} gives $\xi=3$.\EP
%%%%%%%%%%%%%%%%%%%%%%%%%%%%%%%%%%%%%%%%%%%%%%%%%%%%%%%%%%%%%%%%%%%%%%%%%%

%%%%%%%%%%%%%%%%%%%%%%%%%%%%%%%%%%%%%%%%%%%%%%%%%%%%%%%%%%%%%%%%%%%%%%%%%
\appendix
%%%%%%%%%%%%%%%%%%%%%%%%%%%%%%%%%%%%%%%%%%%%%%%%%%%%%%%%%%%%%%%%%%%%%%%%%

%%%%%%%%%%%%%%%%%%%%%%%%%%%%%%%%%%%%%%%%%%%%%%%%%%%%%%%%%%%%%%%%%%%%%%%%%
%%%%%%%%%%%%%%%%%%%%%%%%%%%%%%%%%%%%%%%%%%%%%%%%%%%%%%%%%%%%%%%%%%%%%%%%%
\zerarcounters
\section{Momentum-depending perturbation}
\label{app:a}
%%%%%%%%%%%%%%%%%%%%%%%%%%%%%%%%%%%%%%%%%%%%%%%%%%%%%%%%%%%%%%%%%%%%%%%%%
%%%%%%%%%%%%%%%%%%%%%%%%%%%%%%%%%%%%%%%%%%%%%%%%%%%%%%%%%%%%%%%%%%%%%%%%%

Here we discuss the Hamiltonian case in which the
perturbation depends also on the coordinates $y_{1},\ldots,y_{d}$,
as in (\ref{eq:2.10}). As we shall see, differently
from the $\yy$-independent case, here the
Hamiltonian structure of the system is fundamental.

It is more convenient to work in complex variables $\zz,\,\ww=\zz^{*}$,
with $z_{j}=(y_{j}+\ii\om_{j}x_{j})/\sqrt{2\om_{j}}$,
where the Hamilton equations are of the form
\begin{equation}\left\{\begin{aligned}
-\ii\dot{z}_{j}&=\om_{j}z_{j}+\e\partial_{w_{j}}F(\zz,\ww,\e)+\h_{j}z_{j},\\
\ii\dot{w}_{j}&=\om_{j}w_{j}+\e\partial_{z_{j}}F(\zz,\ww,\e)+\h_{j}w_{j},
\end{aligned}\right.\label{eq:a.1}\end{equation}
with
\begin{equation}
F(\zz,\ww,\e)=\sum_{p=0}^{\infty}\e^{p} \!\!\!\!\!\!\!\!\!\!
\sum_{\substack{
s^{+}_{1},\ldots,s^{+}_{d},s^{-}_{1},\ldots,s^{-}_{d} \ge 0
\\ s^{+}_{1}+\ldots+s^{+}_{d}+s^{-}_{1}+\ldots+s^{-}_{d}=p+3}}
\!\!\!\!\!\!\!\!\!\!
a_{s^{+}_{1},\ldots,s^{+}_{d},s^{-}_{1},\ldots,s^{-}_{d} } \,
z_{1}^{s^{+}_{1}} \ldots z_{d}^{s^{+}_{d}}w_{1}^{s^{-}_{1}}
\ldots w_{d}^{s^{-}_{d}} .
\label{eq:a.2} \end{equation}
Note that, since the Hamiltonian (\ref{eq:2.8}) is real, one has
\begin{equation}
a_{\bs^{+},\bs^{-}}=a^{*}_{\bs^{-},\bs^{+}},\quad 
\bs^{\pm}=(s_{1}^{\pm},\ldots,s_{d}^{\pm})\in\ZZZ_{+}^{d}.
\label{eq:a.3}\end{equation}
Let us write
\begin{equation}
f^{+}_{j}(\zz,\ww,\e)=\e\partial_{w_{j}}F(\zz,\ww,\e),\qquad
f^{-}_{j}(\zz,\ww,\e)=\e\partial_{z_{j}}F(\zz,\ww,\e)
\nonumber\end{equation}
so that
\begin{equation}
f^{\s}_{j}(\zz,\ww,\e)=\sum_{p=1}^{\infty}\e^{p} \!\!\!\!\!\!\!\!\!\!
\sum_{\substack{
\bs^{+},\bs^{-} \in \ZZZ_{+}^{d}
\\ s^{+}_{1}+\ldots+s^{+}_{d}+s^{-}_{1}+\ldots+s^{-}_{d}=p+1}}
\!\!\!\!\!\!\!\!\!\!
f^{\s}_{j,\bs^{+},\bs^{-} } \,
z_{1}^{s^{+}_{1}} \ldots z_{d}^{s^{+}_{d}}w_{1}^{s^{-}_{1}}
\ldots w_{d}^{s^{-}_{d}}, \qquad \s=\pm ,
\nonumber \end{equation}
with $f^{+}_{j,\bs^{+},\bs^{-}}=(s_{j}^{-}+1)a_{\bs^{+},\bs^{-}+\ee_{j}}$
and $f^{-}_{j,\bs^{+},\bs^{-}}=(s_{j}^{+}+1)a_{\bs^{+}+\ee_{j},\bs^{-}}$,
and hence
\begin{subequations}
\begin{align}
& f_{j,\bs^{+},\bs^{-}}^{-}=\left(f_{j,\bs^{-},\bs^{+}}^{+}\right)^{*}, 
\qquad j=1,\ldots,d,\quad \bs^{+},\bs^{-}\in\ZZZ^{d},
\label{eq:a.4a} \\
& (s_{j_{2}}^{+}+1)f_{j_{1},\bs^{+}+\ee_{j_{2}},\bs^{-}}^{+}=
(s_{j_{1}}^{-}+1)f_{j_{2},\bs^{+},\bs^{-}+\ee_{j_{1}}}^{-},\qquad
j_{1},j_{2}=1,\ldots,d,\quad \bs^{+},\bs^{-}\in\ZZZ^{d}, 
\label{eq:a.4b}\\
& (s_{j_{2}}^{-}+1)f_{j_{1},\bs^{+},\bs^{-}+\ee_{j_{2}}}^{+}=
(s_{j_{1}}^{-}+1)f_{j_{2},\bs^{+},\bs^{-}+\ee_{j_{1}}}^{+},
\qquad
j_{1},j_{2}=1,\ldots,d,\quad \bs^{+},\bs^{-}\in\ZZZ^{d},
\label{eq:a.4c}\\
& (s_{j_{2}}^{+}+1)f_{j_{1},\bs^{+}+\ee_{j_{2}},\bs^{-}}^{-}=
(s_{j_{1}}^{+}+1)f_{j_{2},\bs^{+}+\ee_{j_{1}},\bs^{-}}^{-},
\qquad
j_{1},j_{2}=1,\ldots,d,\quad \bs^{+},\bs^{-}\in\ZZZ^{d}.
\label{eq:a.4d}\end{align}
\label{eq:a.4}\end{subequations}
Expanding the solution $(\zz(t),\ww(t))$ in Fourier series with
frequency vector $\oo$, (\ref{eq:a.1}) gives
\begin{equation}\begin{cases}
(\oo\cdot\nn-\om_{j})z_{j,\nn} =
\h_{j}z_{j,\nn}+f^{+}_{j,\nn}(\zz,\ww,\e), & \\ 
(-\oo\cdot\nn-\om_{j})w_{j,\nn} =
\h_{j}w_{j,\nn}+f^{-}_{j,\nn}(\zz,\ww,\e) . &
\end{cases}
\label{eq:a.5}\end{equation}
We write the unperturbed solutions as
\begin{equation}
z_{j}^{(0)}(t)=c_{j}^{+}e^{\ii\om_{j}t}, \qquad
w_{j}^{(0)}(t)=c_{j}^{-}e^{-\ii\om_{j}t}, \qquad
j=1,\ldots,d,
\nonumber\end{equation}
with $c_{j}=c_{j}^{-} \in \CCC$ and $c_{j}^{+}=c_{j}^{*}$.
As in Section \ref{sec:2} we can split (\ref{eq:a.5}) into
\begin{subequations}
\begin{align}
& f^{+}_{j,\ee_{j}}(\zz,\ww,\e) + \h_{j} \, z_{j,\ee_{j}} = 0 ,
\qquad j=1,\ldots,d,
\label{eq:a.6a} \\
& f^{-}_{j,-\ee_{j}}(\zz,\ww,\e) + \h_{j} \, w_{j,-\ee_{j}} = 0 ,
\qquad j=1,\ldots,d,
\label{eq:a.6b} \\
& \left[ (\oo\cdot\nn) - \om_{j} \right] z_{j,\nn} =
f^{+}_{j,\nn}(\zz,\ww,\e) + \h_{j} \, z_{j,\nn} ,
\qquad j=1,\ldots,d , \quad \nn \neq  \ee_{j} ,
\label{eq:a.6c} \\
& \left[ -(\oo\cdot\nn) - \om_{j} \right] w_{j,\nn} =
f^{-}_{j,\nn}(\zz,\ww,\e) + \h_{j} \, w_{j,\nn} ,
\qquad j=1,\ldots,d , \quad \nn \neq  -\ee_{j} ,
\label{eq:a.6d}
\end{align}
\end{subequations}
so that first of all one has to show that the same choice of
$\h_{j}$ makes both (\ref{eq:a.6a}) and (\ref{eq:a.6b})
to hold simultaneously, and that such $\h_{j}$ is real.

We consider a tree expansion very close to the one performed in Section
\ref{sec:4}: we simply drop (3) in Constraint \ref{sec:4.4}.
We denote by $\gotT^{k}_{j,\nn,\s}$ the set of inequivalent trees
of order $k$, tree component $j$, tree momentum $\nn$ and
{\it tree sign} $\s$ that is, the sign label of the root line is $\s$.

We introduce $\thetapru$ and $\thetao$ as in Notation \ref{def:4.5}
and \ref{def:4.27} respectively, and we define the value of a tree
as follows.

The node factors are defined as in (\ref{eq:4.2}) for the end nodes,
while for the internal nodes $v\in V(\theta)$ we define
\begin{equation}
F_{v} = \begin{cases}
\displaystyle{ \frac{s_{v,1}^{+}!\ldots s_{v,d}^{+}!
s_{v,1}^{-}!\ldots s_{v,d}^{-}!}{s_{v}!}
f^{\s_{\ell_{v}}}_{j_{v},\bs_{v}^{+},\bs_{v}^{-}} } ,
& \qquad k_{v} \ge 1 , \\
& \\
\displaystyle{ - \frac{1}{2c_{j_{v}}^{\s_{v}}} } ,
& \qquad k_{v} = 0 .
\end{cases}
\label{eq:a.7}\end{equation}

The propagators are defined as $G_{\ell}=1$ if $\nn_{\ell}=
\s_{\ell}\ee_{j_{\ell}}$
and
\begin{equation}
G_{\ell} = G_{j_{\ell}}^{[n_{\ell}]}(\s_{\ell}\,\oo\cdot\nn_{\ell}), \qquad
G_{j}^{[n]}(u)=\frac{\Psi_{n}(|u-\om_{j}|)}{u-\om_{j}},
\label{eq:a.8}\end{equation}
otherwise, and we define $\Val(\theta)$ as in (\ref{eq:4.9}).

Finally we set $z_{j,\ee_{j}}=w_{j,-\ee_{j}}^{*}=c_{j}$,
and formally define
\begin{equation}
\begin{aligned}
& {z}_{j,\nn}=\sum_{k=1}^{\io}\e^{k}z^{(k)}_{j,\nn},
\qquad
z^{(k)}_{j,\nn}=\sum_{\theta \in \gotT^{k}_{j,\nn,+}}\Val(\theta),
\qquad\nn\neq\ee_{j},\\
& {w}_{j,\nn}=\sum_{k=1}^{\io}\e^{k}w^{(k)}_{j,\nn},
\qquad
w^{(k)}_{j,\nn}=\sum_{\theta \in \gotT^{k}_{j,\nn,-}}\Val(\theta),
\qquad \nn\neq-\ee_{j},
\end{aligned}\label{eq:a.9}\end{equation}
and
\begin{equation}
{\h}_{j,\s}=\sum_{k=1}^{\io}\e^{k}\h^{(k)}_{j,\s}, \qquad
\h^{(k)}_{j,\s}=-\frac{1}{c_{j}^{\s}}
\sum_{\theta \in \gotT^{k}_{j,\s\ee_{j},\s}}
\Val(\theta).
\label{eq:a.10}\end{equation}
Note that Remarks \ref{rmk:4.9}, \ref{rmk:4.13} and \ref{rmk:4.17}
still hold.

%%%%%%%%%%%%%%%%%%%%%%%%%%%%%%%%%%%%%%%%%%%%%%%%%%%%%%%%%%%%%%%%%%%%%%%%%%
\begin{lemma} \label{lem:a.1}
With the notations introduced above, one has ${\h}_{j,+}^{*}=
{\h}_{j,-}$ and $z_{j,\nn}^{*}=w_{j,-\nn}$.
\end{lemma}
%%%%%%%%%%%%%%%%%%%%%%%%%%%%%%%%%%%%%%%%%%%%%%%%%%%%%%%%%%%%%%%%%%%%%%%%%%

%%%%%%%%%%%%%%%%%%%%%%%%%%%%%%%%%%%%%%%%%%%%%%%%%%%%%%%%%%%%%%%%%%%%%%%%%%
\prova By definition of we only have to prove
that for any $\theta\in\gotT^{k}_{j,\nn,+}$ there exists $\theta'\in
\gotT^{k}_{j,-\nn,-}$ such that $\Val(\theta)^{*}=\Val(\theta')$.
The proof is by induction on the order of the tree.
Given $\theta\in\gotT^{k}_{j,\nn,+}$, let us consider the tree
$\theta'$ obtained from $\theta$ by replacing the labels $\s_{v}$ of
all the nodes $v\in N_{0}(\theta)$ with $-\s_{v}$ and the labels
$\s_{\ell}$ of all the lines $\ell\in L(\theta)$ with $-\s_{\ell}$.
Call $\ell_{1},\ldots,\ell_{p}$ the lines on scale $-1$ (if any) closest
to the root of $\theta$, and denote by $v_{i}$ the node $\ell_{i}$ enters
and by $\theta_{i}$ the tree with root line $\ell_{i}$.
Each tree $\theta_{i}$ is then replaced with a tree $\theta_{i}'$ such
that $\Val(\theta_{i})^{*}=\Val(\theta_{i}')$ by the inductive hypothesis.
Moreover, as for any internal line in $\theta$ the momentum becomes
$-\nn_{\ell}$, the propagators do not change. Finally, for any
$v\in V(\thetapru)$ the node factor is changed into
\begin{equation}
F_{v}' = \begin{cases}
\displaystyle{  \frac{s_{v,1}^{-}!\ldots s_{v,d}^{-}!
s_{v,1}^{+}!\ldots s_{v,d}^{+}!}{s_{v}!}
f^{-\s_{\ell_{v}}}_{j_{v},\bs_{v}^{-},\bs_{v}^{+}} } ,
& k_{v} \ge 1 , \\
& \\
\displaystyle{ - \frac{1}{2c_{j_{v}}^{-\s_{v}}} } ,
& k_{v} = 0 .
\end{cases}
\label{eq:a.11}\end{equation}
Hence by (\ref{eq:a.4a}) one has $\Val(\theta)^{*}=\Val(\theta')$.\EP
%%%%%%%%%%%%%%%%%%%%%%%%%%%%%%%%%%%%%%%%%%%%%%%%%%%%%%%%%%%%%%%%%%%%%%%%%%

%%%%%%%%%%%%%%%%%%%%%%%%%%%%%%%%%%%%%%%%%%%%%%%%%%%%%%%%%%%%%%%%%%%%%%%%%%
\begin{lemma} \label{lem:a.2}
With the notations introduced above, one has ${\h}_{j,+}\in\RRR$.
\end{lemma}
%%%%%%%%%%%%%%%%%%%%%%%%%%%%%%%%%%%%%%%%%%%%%%%%%%%%%%%%%%%%%%%%%%%%%%%%%%

%%%%%%%%%%%%%%%%%%%%%%%%%%%%%%%%%%%%%%%%%%%%%%%%%%%%%%%%%%%%%%%%%%%%%%%%%%
\prova We only have to prove that for any $\theta\in\gotT^{k}_{j,\ee_{j},+}$
there exists $\theta''\in\gotT^{k}_{j,\ee_{j},+}$ such that
$$
c_{j}^{+}\Val(\theta)^{*}=c_{j}^{-}\Val(\theta'').
$$
Let $v_{0}\in E_{j}^{+}(\thetapru)$ (existing by Remark \ref{rmk:4.9})
and let us consider the tree $\theta''$ obtained from $\theta$ by
(1) exchanging the root line $\ell_{0}$ with $\ell_{v_{0}}$, (2) replacing
all the labels $\s_{v}$ of all the end-nodes $v\in N_{0}(\theta)
\setminus\{v_{0}\}$ with $-\s_{v}$, and (3) replacing all the labels
$\s_{\ell}$ of all the internal lines with $-\s_{\ell}$, except for those
in $\calP(\ell_{v_{0}},\ell_{0})$ which remain the same. The propagators
do not change; this is trivial for the lines outside
$\calP(\ell_{v_{0}},\ell_{0})$, while for $\ell\in\calP(\ell_{v_{0}},
\ell_{0})$ one can reason as follows.
The line $\ell$ divides $E(\thetapru)\setminus\{v_{0}\}$ into two disjoint
sets of end nodes $E(\thetapru,p)$ and $E(\thetapru,s)$ such that if
$\ell=\ell_{w}$ one has $E(\thetapru,p)=\{v\in E(\thetapru)\setminus
\{v_{0}\}:v\prec w\}$ and $E(\thetapru,s)=(E(\thetapru)\setminus\{v_{0}\})
\setminus E(\thetapru,p)$. %See Figure \ref{???} mettiamo?
If
$$
\nn^{(p)}=\sum_{v\in E(\thetapru,p)}\nn_{v},\qquad
\nn^{(s)}=\sum_{v\in E(\thetapru,s)}\nn_{v},
$$
one has $\nn^{(p)}+\nn^{(s)}=\vzero$. When considering $\ell$ as a line in
$\theta$ one has $\nn_{\ell}=\nn^{(p)}+\ee_{j}$ while in $\theta''$ one has
$\nn_{\ell}=-\nn^{(s)}+\ee_{j}$. Hence, as we have not changed the sign
label $\s_{\ell}$, also $G_{\ell}$ does not change.
The node factors of the internal nodes are changed into their complex
conjugated; this can be
obtained as in Lemma \ref{lem:a.1} for the internal nodes $w$ such that
$\ell_{w}\notin \calP(\ell_{v_{0}},\ell_{0})$ while for the other nodes one
can reason as follows.

First of all if $v$ is such that $\ell_{v}\in\calP(\ell_{v_{0}},\ell_{0})$,
there is a line $\ell_{v}'\in\calP(\ell_{v_{0}},\ell_{0})$
entering $v$.
We shall denote $j_{\ell_{v}}=j_{1}$, $\s_{\ell_{v}}=\s$, $j_{\ell_{v}'}=j_{2}$,
and $\s_{\ell_{v}'}=\s'$. 
Moreover we call $s_{i}^{\s''}$ the number of
lines outside $\calP(\ell_{v_{0}},\ell_{0})$ with component label $i$
and sign label $\s''$ entering $v$.
Let us consider first the case $\s=\s'=+$.
When considering $v$ as node of $\theta$ one has
$$
F_{v}^{*}=\left(\frac{s_{1}^{+}!\ldots s_{d}^{+}! s_{1}^{-}!\ldots
s_{d}^{-}! (s_{j_{2}}^{+}+1)}{s_{v}!}
f_{j_{1},\bs^{+}+\ee_{j_{2}},\bs^{-}}^{+}\right)^{*}
=\frac{s_{1}^{+}!\ldots s_{d}^{+}! s_{1}^{-}!\ldots
s_{d}^{-}! (s_{j_{2}}^{+}+1)}{s_{v}!}
f_{j_{1},\bs^{-},\bs^{+}+\ee_{j_{2}}}^{-}.
$$
When considering $v$ as node of $\theta''$ one has
$\bs_{v}^{+}=\bs^{-}+\ee_{j_{1}}$ and $\bs_{v}^{-}=\bs^{+}$, so that
$$
F_{v}''=\frac{s_{1}^{+}!\ldots s_{d}^{+}! s_{1}^{-}!\ldots
s_{d}^{-}! (s_{j_{1}}^{-}+1)}{s_{v}!}
f_{j_{2},\bs^{-}+\ee_{j_{1}},\bs^{+}}^{+},
$$
and hence by (\ref{eq:a.4b}) $F_{v}^{*}=F_{v}''$.
Reasoning analogously one obtains $F_{v}^{*}=F_{v}''$
also in the cases $\s=\s'=-$ and $\s\neq\s'$, using again
(\ref{eq:a.4b}) when $\s=\s'=-$, and (\ref{eq:a.4c}) and (\ref{eq:a.4d})
for $\s=-$, $\s'=+$ and $\s=+$, $\s'=-$ respectively.
Hence the assertion is proved.\EP
%%%%%%%%%%%%%%%%%%%%%%%%%%%%%%%%%%%%%%%%%%%%%%%%%%%%%%%%%%%%%%%%%%%%%%%%%%

We define the self-energy clusters as in Section \ref{sec:4.6},
but replacing the constraint (3) with
(3$'$) one has $|\nn_{\ell_{T}}-\nn_{\ell'_{T}}|\le 2$ and 
$|\s_{\ell_{T}}\,\oo\cdot\nn_{\ell_{T}}-\om_{j_{\ell_{T}}}| =
|\s_{\ell'_{T}}\,\oo\cdot\nn_{\ell'_{T}}-\om_{j_{\ell'_{T}}}|$.
We introduce $\Tpru$ and $\To$ as in Notation \ref{def:4.23} and
\ref{def:4.27} respectively, and we can define
$\Val(T)$ as in (\ref{eq:4.12}) and the localisation and the
regularisation operators as in Section \ref{sec:6}.

Note that the main difference with the $\yy$-independent case is in the
role of the sign label $\s_{\ell}$. In fact, here the sign label of a line
does not depend on its momentum and component labels, and the small
divisor is given by $\delta_{j,\s}(\oo\cdot\nn)=|\s\oo\cdot\nn-\om_{j}|$.

Hence the dimensional bounds of Section \ref{sec:5} and the symmetries
discussed in Section \ref{sec:7} and summarised
in Lemma \ref{lem:7.1} can be proved word by word as in the
$\yy$-independent case, except for the second equality in Lemma
\ref{lem:7.3} where one has to take into account a change of signs.
More precisely for $T \in \gotR^{k}_{j,\s,j',\s'}(u,n)$, with $j\neq j'$
and $n_{\ell}\ge 0$ for all $\ell \in \calP_{T}$, we define $\calG_{1}(T)$
as in Section \ref{sec:7} and $\calG_{3}(T)$ as in Section \ref{sec:7}
but replacing also the sign labels $\s_{\ell}$ of the lines
$\ell \in L(T)$ with $-\s_{\ell}$.

%%%%%%%%%%%%%%%%%%%%%%%%%%%%%%%%%%%%%%%%%%%%%%%%%%%%%%%%%%%%%%%%%%%%%%%%%%
\begin{lemma} \label{lem:a.3}
For all $T \in \gotR^{k}_{j,\s,j',\s'}(u,n)$, with $j\neq j'$ and
$n_{\ell}\ge 0$ for all $\ell \in \calP_{T}$, one has
\begin{equation}\label{eq:a.12}
c_{j}^{-\s}c_{j'}^{\s'} \!\!\! \sum_{T'\in\calG_{1}(T)}\!\!\! \matL\Val(T',u)
= c_{j}^{\s}c_{j'}^{-\s'} \!\!\! \sum_{T'\in\calG_{3}(T)} \!\!\!
\matL\Val(T',u).
\end{equation}
\end{lemma}
%%%%%%%%%%%%%%%%%%%%%%%%%%%%%%%%%%%%%%%%%%%%%%%%%%%%%%%%%%%%%%%%%%%%%%%%%%

%%%%%%%%%%%%%%%%%%%%%%%%%%%%%%%%%%%%%%%%%%%%%%%%%%%%%%%%%%%%%%%%%%%%%%%%%%
\prova We consider only the case $k(\To)\le K_{0}2^{n_{T}/\tau}$.
For fixed $T \in \gotR^{k}_{j,\s,j',\s'}(u,n)$, with $j\neq j'$,
let $\theta\in\gotT_{j,\s\ee_{j},\s}^{k}(n)$ be the tree obtained
from $T$ by replacing the entering line $\ell_{T}'$  with a line
exiting a new end node $v_{0}$ with $\s_{v_{0}}=\s'$ and $\nn_{v_{0}}=
\s'\ee_{j'}$. 
As in the proof of Lemma \ref{lem:7.3} one has
\begin{equation}
c_{j'}^{\s'} \!\!\! \sum_{T'\in\calG_{1}(T)} \!\!\! \matL\Val(T',u)=
|E_{j'}^{\s'}(\thetapru)|\Val(\theta).
\nonumber \end{equation}
Now, let $\theta'\in\gotT_{j,-\s\ee_{j},-\s}^{k}(n)$ be the tree obtained
from $\theta$ by replacing all the labels $\s_{v}$ of
the nodes $v\in N_{0}(\theta)$ with $-\s_{v}$, and the labels $\s_{\ell}$
of all the lines $\ell\in L(\theta)$ with $-\s_{\ell}$.
Any $T'\in\calG_{3}(T)$ can be obtained from $\theta'$ by
replacing a line exiting an end node $v \in E_{j'}^{\s'}(\thetapru')$
with entering line $\ell_{T'}'$, carrying the same labels as $\ell_{T}$.
Hence, by Lemma \ref{lem:a.1},
\begin{equation}
c_{j}^{-\s}c_{j'}^{\s'} \!\!\! \sum_{T'\in\calG_{1}(T)} \!\!\!
\matL\Val(T',u)= c_{j}^{-\s}|E_{j'}^{\s'}(\thetapru)|
\Val(\theta)=c_{j}^{-\s} |E_{j'}^{\s'}(\thetapru')|
\Val(\theta')^{*}= c_{j}^{-\s}c_{j'}^{-\s'} \!\!\!
\sum_{T'\in\calG_{3}(T)} \!\!\! (\matL\Val(T',u))^{*} .
\nonumber \end{equation}
On the other hand, exactly as in Lemma \ref{lem:a.2} one can prove that
for any $T'\in\calG_{3}(T)$ there exists $T''\in\calG_{3}(T)$ such that
\begin{equation}
c_{j}^{-\s} (\matL\Val(T',u))^{*}= c_{j}^{\s}\matL\Val(T'',u),
\nonumber\end{equation}
and hence the assertion follows.\EP
%%%%%%%%%%%%%%%%%%%%%%%%%%%%%%%%%%%%%%%%%%%%%%%%%%%%%%%%%%%%%%%%%%%%%%%%%%

The cancellation mechanism and the bounds
proved in Section \ref{sec:8} follows by the same reasoning
(in fact it is even simpler); see the next appendix for details.

%%%%%%%%%%%%%%%%%%%%%%%%%%%%%%%%%%%%%%%%%%%%%%%%%%%%%%%%%%%%%%%%%%%%%%%%%
%%%%%%%%%%%%%%%%%%%%%%%%%%%%%%%%%%%%%%%%%%%%%%%%%%%%%%%%%%%%%%%%%%%%%%%%%
\zerarcounters
\section{Matrix representation of the cancellations}
\label{app:b}
%%%%%%%%%%%%%%%%%%%%%%%%%%%%%%%%%%%%%%%%%%%%%%%%%%%%%%%%%%%%%%%%%%%%%%%%%
%%%%%%%%%%%%%%%%%%%%%%%%%%%%%%%%%%%%%%%%%%%%%%%%%%%%%%%%%%%%%%%%%%%%%%%%%

As we have discussed in Section \ref{sec:6} the only obstacle to
convergence of the formal power series of the solution
is given by the accumulation of resonant lines;
see Figure \ref{fig:14}. 

The cancellation mechanism described in Section \ref{sec:8} can be
expressed in matrix notation. This is particularly helpful in the
$\yy$-dependent case. For this reason, and for the fact that
the formalism introduced in Appendix \ref{app:a} include
the $\yy$-independent case, we prefer to work here
with the variables $(\zz,\ww)$.

We first develop a convenient notation.
Given $\nn$ such that $\s(\nn,1)=+$ and $\de_{1,+}(\om\cdot\nn)<\gamma$
let us group together, in an ordered set $S(\nn)$, all the $\nn'$
such that $\nn' = \nn'(j,\s) := \nn-\ee_{1}+\s \ee_j$, $\s=\pm 1$
and $j=1,\dots,d$; see Remark \ref{rmk:4.19}.
By definition one has  $\de_{1,+}(\om\cdot\nn)= 
\de_{j,\s}(\om\cdot\nn'(j,\s))$ for all $j=1,\ldots,d$ and $\s=\pm$.
Then we construct a $2d\times 2d$ \emph{localised self-energy matrix}
$\matL M^{(k)}(\oo\cdot\nn,n)$ with entries $\matL
M^{(k)}_{j,\s,j',\s'}(\oo\cdot\nn'(j',\s'),n)$.
We also define the $2d\times 2d$ diagonal propagator matrix 
$\matG^{[n]}(\oo\cdot\nn)$ with entries $\matG^{[n]}_{j,\s,j,\s}
(\oo\cdot \nn) = \de_{j,j'}\de_{\s,\s'}G^{[n]}_{j'}(\oo\cdot \nn'(j,\s))$,
with $G^{[n]}_{j'}(u)$ defined according to (\ref{eq:a.8}).

As in Section \ref{sec:8} let us consider a chain  of two
self-energy clusters; see figure \ref{fig:9}.
By definition its value is
\begin{equation}
\matL\Val(T_{1},\om\cdot \nn_{1}) \,
G^{[n_{\ell}]}_{j_{\ell}}(\om\cdot \nn_{\ell}) \,
\matL\Val(T_{2},\om\cdot \nn_{2}),
\nonumber \end{equation}
with $\nn_{1}=\nn_{\ell_{T_{1}}'}$ and $\nn_{2}=\nn_{\ell_{T_{2}}}$.

Notice that, if one sets also for sake of simplicity,
$\s_{1}=\s_{\ell_{T_{1}}'}$, $j_{1}=j_{\ell_{T_{1}'}}$,
$\s_{2}=\s_{\ell_{T_{2}}}$, and $j_{2}=j_{\ell_{T_{2}}}$,
by the constraint (3$'$) in the definition of
self-energy clusters given in Appendix \ref{app:a},
one has $\nn_{1}-\nn_{\ell} = \s_{1}\ee_{j_{1}}-\s_{\ell}\ee_{j_{\ell}}$
and $\nn_{\ell}-\nn_{2} = \s_{\ell}\ee_{j_{\ell}}-\s_{2}\ee_{j_{2}}$;
moreover $\nn_{1},\nn_{\ell},\nn_{2}$ all belong to a single
set $S(\nn)$ for some $\nn$.

As done in Section \ref{sec:8} let us sum together the values of
all the possible self-energy clusters $T_{1}$ and $T_{2}$
with fixed labels associated with the external lines, and
of fixed orders $k_{1}$ and $k_{2}$, respectively.
We obtain
\begin{equation}
\matL M^{(k_{1})}_{j_{1},\s_{1},j_{\ell},\s_{\ell}}
(\oo\cdot\nn'(j_{\ell},\s_{\ell}),n_{T_{1}})\,
\matG^{[n_{\ell}]}_{j_{\ell},\s_{\ell},j_{\ell},\s_{\ell}}(\oo\cdot\nn)\,
\matL M^{(k_{2})}_{j_{\ell},\s_{\ell},j_{2},\s_{2}}
(\oo\cdot\nn'(j_{2},\s_{2}),n_{T_{2}}) .
\nonumber \end{equation}
If we also sum over all possible values of the labels
$j_{\ell},\s_{\ell}$ we get 
\begin{eqnarray}
& \displaystyle{
\sum_{\s_{\ell}=\pm}\sum_{j_{\ell}=1}^{d}
\matL M^{(k_{1})}_{j_{1},\s_{1},j_{\ell},\s_{\ell}}
(\oo\cdot\nn'(j_{\ell},\s_{\ell}),n_{T_{1}})\, 
\matG^{[n_{\ell}]}_{j_{\ell},\s_{\ell},j_{\ell},\s_{\ell}}(\oo\cdot\nn)
\matL M^{(k_{2})}_{j_{\ell},\s_{\ell},j_{2},\s_{2}}
(\oo\cdot\nn'(j_{2},\s_{2}),n_{T_{2}}) }
\nonumber \\
& \qquad \qquad \qquad \qquad
= \displaystyle{
\left[ \matL M^{(k_{1})} (\oo\cdot\nn,n_{T_{1}})\,
\matG^{[n_{\ell}]} (\oo\cdot\nn) \,\matL M^{(k_{2})}
(\oo\cdot\nn,n_{T_{2}}) \right]_{j_{1},\s_{1},j_{2},\s_{2}} } ,
\nonumber \end{eqnarray} 
(i.e. the entry $j_{1},\s_{1},j_{2},\s_{2}$ of the matrix in square
brackets).

By the definition (\ref{eq:a.8}) of the propagators and by the
symmetries of Lemma \ref{lem:7.1}, $\matG^{[n]}(\oo\cdot\nn)$
and $\matL M^{(k)}(\oo\cdot\nn,n)$ have the form
\begin{equation}
\matG^{[n]}(\oo\cdot\nn)=
\frac{\Psi_{n}(|\oo\cdot\nn-\om_{1}|)}{\oo\cdot\nn-\om_{1}}
\begin{pmatrix}\begin{pmatrix} 1 & 0 \\
0 & -1 \end{pmatrix}& 0 & \cdots &0 \\
0 & \ddots& \ddots &\vdots \\
\vdots &\ddots &\ddots &0\\
0 &\cdots&0&\begin{pmatrix}
1 & 0 \\ 0 & -1 \end{pmatrix}\end{pmatrix} ,
\label{eq:b.1} \end{equation}
and
\begin{equation}
\matL M^{(k)}(\oo\cdot\nn,n)= 
\begin{pmatrix} M^{(k)}_{1,1}(n) \begin{pmatrix}
c^{*}_{1}c_{1} & c^{*}_{1}c^{*}_{1} \\
c_{1}c_{1} & c_{1}c_{1}^{*} \end{pmatrix}
& \cdots & M^{(k)}_{1,d}(n) \begin{pmatrix}
c^{*}_{1}c_{d} & c^{*}_{1}c^{*}_{d} \\
c_{1}c_{d} & c_{1}c_{d}^{*} \end{pmatrix}\\
\vdots & M^{(k)}_{j,j'}(n) \begin{pmatrix}
c^{*}_{j}c_{j'} & c^{*}_{j}c^{*}_{j'} \\
c_{j}c_{j'} & c_{j} c_{j'}^{*} \end{pmatrix}
& \vdots \\M^{(k)}_{d,1}(n)
\begin{pmatrix} c^{*}_{d}c_{1} & c^{*}_{d}c^{*}_{1} \\
c_{d}c_{1} & c_{d}c_{1}^{*} \end{pmatrix}
& \cdots & M^{(k)}_{d,d}(n) \begin{pmatrix}
c^{*}_{d}c_{d}& c^{*}_{d}c^{*}_{d} \\
c_{d}c_{d} & c_{d}c_{d}^{*} \end{pmatrix} \end{pmatrix} ,
\nonumber \end{equation}
respectively. A direct computation gives
\begin{eqnarray}
& \hskip-4truecm \displaystyle{
\left[ \matL M^{(k_{1})}(\oo\cdot\nn,n_{T_{1}})\,
\matG^{[n_{\ell}]}(\oo\cdot\nn)\, \matL M^{(k_{2})}(\oo\cdot\nn,n_{T_{2}})
\right]_{j_{1},\s_{1},j_{2},\s_{2}} } \nonumber \\
& %\qquad\qquad\qquad\qquad
\displaystyle{ =
\frac{\Psi_{n_{\ell}}(|\oo\cdot\nn-\om_{1}|)}{\oo\cdot\nn-\om_{1}}
c_{j_{1}}^{-\s_{1}} c_{j_{2}}^{\s_{2}}
\sum_{j=1}^{d} M_{j_{1},j}(n_{T_{1}}) \, M_{j,j_{2}}(n_{T_{2}})\,
|c_{j}|^{2} \sum_{\s=\pm} (-1)^{1+\s 1} = 0 } ,
\label{eq:b.2} \end{eqnarray}
for all choices of the scales $n_{\ell},n_{T_{1}},n_{T_{2}}$
and of the orders $k_{1},k_{2}$.

This proves the necessary cancellation. Note that this is
an exact cancellation in terms of the variables $(\zz,\ww)$:
all chains of localised self-energy clusters of length
$p \geq 2$  can be ignored as their values sum up to zero.
In the $\yy$-independent case, and in terms of the variables $\xx$,
the cancellation is only partial, and one only finds
$\matL M^{(k_{1})} G^{[n]} \matL M^{(k_{2})} =O( 2^{-n})$,
as discussed in Section \ref{sec:8}.

%%%%%%%%%%%%%%%%%%%%%%%%%%%%%%%%%%%%%%%%%%%%%%%%%%%%%%%%%%%%%%%%%%%%%%%%%
%%%%%%%%%%%%%%%%%%%%%%%%%%%%%%%%%%%%%%%%%%%%%%%%%%%%%%%%%%%%%%%%%%%%%%%%%
\zerarcounters
\section{Resummation of the perturbation series}
\label{app:c}
%%%%%%%%%%%%%%%%%%%%%%%%%%%%%%%%%%%%%%%%%%%%%%%%%%%%%%%%%%%%%%%%%%%%%%%%%
%%%%%%%%%%%%%%%%%%%%%%%%%%%%%%%%%%%%%%%%%%%%%%%%%%%%%%%%%%%%%%%%%%%%%%%%%

The fact that the series obtained by systematically eliminating
the self-energy clusters converges, as seen in Section \ref{sec:5},
suggests that one may follow another approach, alternative to
that we have described so far, and leading to the same result.
Indeed, one can consider a \emph{resummed expansion},
where one really gets rid of the
self-energy clusters at the price of changing the propagators into new
\emph{dressed propagators} -- again terminology is
borrowed from quantum field theory.
This is a standard procedure, already exploited in the case
of KAM tori \cite{GBG}, lower-dimensional tori \cite{GBG,G2},
skew-product systems \cite{G1}, etc. The convergence of the
perturbation series reflects the fact that the dressed propagators
can be bounded proportionally to (a power of) the original ones
for all values of the perturbation parameter $\e$.
In our case, the latter property can be seen as a consequence of the
cancellation mechanism just described. In a few words --
and oversimplifying the strategy -- the
dressed propagators are obtained starting from a tree expansion
where no self-energy clusters are allowed, and then `inserting
arbitrary chains of self-energy clusters': this means that
each propagator $\matG^{[n]}=\matG^{[n]}(\oo\cdot\nn)$
is replaced by a dressed propagator
\begin{equation}
\Gamma^{[n]} = \matG^{[n]} +
\matG^{[n]} M \matG^{[n]} +
\matG^{[n]} M \matG^{[n]} M
\matG^{[n]} + \ldots ,
\label{eq:b.3} \end{equation}
where $M=M(\oo\cdot\nn)$ denotes the insertion of all possible
self-energy clusters compatible with the labels of the propagators
of the external lines ($M$ is is the matrix with entries
$M_{j,\s,j'\s'}(\oo\cdot\nn'(j',\s'))$ formally defined
in Remark \ref{rmk:4.26}). Then, formally, one can sum together
all possible contributions in (\ref{eq:b.3}), so as to obtain
\begin{equation}
\Gamma^{[n]} = \matG^{[n]} \left(
\uno - M \matG^{[n]} \right)^{-1} =
\left( A^{-1} - B \right)^{-1} , \qquad A:=\matG^{[n]} ,
\quad B:= M .
\label{eq:b.4} \end{equation}
For sake of simplicity, let us also identify the self-energy values
with their localised parts, so as to replace in (\ref{eq:b.3}),
and hence in (\ref{eq:b.4}), $M$ with $\matL M$,
if $\matL$ is the localisation operator.
Then, in the notations we are using, the cancellation (\ref{eq:b.2})
reads $BAB=0$, which implies
\begin{equation}
\Gamma^{[n]} = A + A B A .
\nonumber \end{equation}
Therefore one finds $\|\Gamma^{[n]}\| \le
\|A\| + \| A\|^{2} \| B \| = O(2^{2n})$. 
So the values of the trees appearing in the resummed expansion
can be bounded as done in Section \ref{sec:5}, with the only
difference that now, instead of the propagators $G_{\ell}$ bounded
proportionally to $2^{n_{\ell}}$, one has the dressed propagators
$\Gamma^{[n_{\ell}]}$ bounded proportionally to $2^{2n_{\ell}}$.

Of course, the argument above should be made more precise.
First of all one should have to take into account also the regularised
values of the self-energy clusters. Moreover, the dressed propagators
should be defined recursively, by starting from the lower scales:
indeed, the dressed propagator of a line on scale $n$
is defined in terms of the values of the self-energy clusters
on scales $<n$, as in (\ref{eq:b.4}),
and the latter in turn are defined in terms of
(dressed) propagators on scales $<n$, according to (\ref{eq:4.13}).
As a consequence, the cancellation mechanism becomes more involved
because the propagators are no longer of the form (\ref{eq:b.1});
in particular the symmetry properties of the self-energy values
should be proved inductively on the scale label.
In conclusion, really proceeding by following the strategy outlined
above requires some work (essentially the same amount as performed
in this paper). We do not push forward the analysis, which in principle
could be worked out by reasoning as done in the papers quoted above.

%%%%%%%%%%%%%%%%%%%%%%%%%%%%%%%%%%%%%%%%%%%%%%%%%%%%%%%%%%%%%%%%%%%%%%%%%%
%%%%%%%%%%%%%%%%%%%%%%%%%%%%%%%%%%%%%%%%%%%%%%%%%%%%%%%%%%%%%%%%%%%%%%%%%%
% References
%%%%%%%%%%%%%%%%%%%%%%%%%%%%%%%%%%%%%%%%%%%%%%%%%%%%%%%%%%%%%%%%%%%%%%%%%%
%%%%%%%%%%%%%%%%%%%%%%%%%%%%%%%%%%%%%%%%%%%%%%%%%%%%%%%%%%%%%%%%%%%%%%%%%%

\end{document}